\documentclass[12pt]{article} 
\usepackage[round,colon,authoryear]{natbib}

\RequirePackage[OT1]{fontenc}
\RequirePackage{amsthm,amsmath,amsfonts,amssymb}
\RequirePackage[colorlinks]{hyperref}
\hypersetup{
 colorlinks=false,
 citecolor=Violet,
 linkcolor=Red,
 urlcolor=Blue}
\RequirePackage{epsfig, graphicx, color}
\RequirePackage{latexsym}
\RequirePackage{psfrag}

\usepackage{fullpage}
\usepackage{epsf,epsfig,subfigure}
\usepackage{graphics,graphicx}
\usepackage{enumerate}
\usepackage{rotating}
\usepackage{epsfig}

\theoremstyle{plain}

\newtheorem{theo}{Theorem}[section]

\newtheorem{lem}{Lemma}[section]
\newtheorem{prop}{Proposition}[section]
\newtheorem{cor}{Corollary}[section]

\theoremstyle{definition} 

\newtheorem{nota}{Notation}[section]
\newtheorem{de}{Definition}[section]
\newtheorem{exa}{Example}[section]
\newtheorem{as}{Assumption}[section]
\newtheorem{alg}{Algorithm}[section]

\newcommand{\btheo}{\begin{theo}}
\newcommand{\bde}{\begin{de}}
\newcommand{\ble}{\begin{lem}}
\newcommand{\bpr}{\begin{prop}}
\newcommand{\bno}{\begin{nota}}
\newcommand{\bex}{\begin{exa}}
\newcommand{\bcor}{\begin{cor}}
\newcommand{\spro}{\begin{proof}}
\newcommand{\bas}{\begin{as}}
\newcommand{\balg}{\begin{alg}}

\newcommand{\etheo}{\end{theo}}
\newcommand{\ede}{\end{de}}
\newcommand{\ele}{\end{lem}}
\newcommand{\epr}{\end{prop}}
\newcommand{\eno}{\end{nota}}
\newcommand{\eex}{\end{exa}}
\newcommand{\ecor}{\end{cor}}
\newcommand{\fpro}{\end{proof}}
\newcommand{\eas}{\end{as}}
\newcommand{\ealg}{\end{alg}}

\theoremstyle{plain}

\newtheorem{theos}{Theorem}
\newtheorem{props}{Proposition}
\newtheorem{lems}{Lemma}
\newtheorem{cors}{Corollary}

\theoremstyle{definition}
\newtheorem{exas}{Example}
\newtheorem{algs}{Algorithm}
\newtheorem{asss}{Asumption}
\newtheorem{defns}{Definition}

\newcommand{\btheos}{\begin{theos}}
\newcommand{\etheos}{\end{theos}}
\newcommand{\bprops}{\begin{props}}
\newcommand{\eprops}{\end{props}}
\newcommand{\bdes}{\begin{defns}}
\newcommand{\edes}{\end{defns}}
\newcommand{\blems}{\begin{lems}}
\newcommand{\elems}{\end{lems}}
\newcommand{\bcors}{\begin{cors}}
\newcommand{\ecors}{\end{cors}}
\newcommand{\bexs}{\begin{exas}}
\newcommand{\eexs}{\end{exas}}
\newcommand{\balgs}{\begin{algs}}
\newcommand{\ealgs}{\end{algs}}
\newcommand{\bass}{\begin{asss}}
\newcommand{\eass}{\end{asss}}

\def\RR{\mathbb{R}}
\def\calA{\mathcal{A}}
\def\calB{\mathcal{B}}
\def\calR{\mathcal{R}}
\def\calM{\mathcal{M}}

\def\calT{\mathcal{T}}

\def\one{\mathbf{1}}
\def\zero{\mathbf{0}}
\DeclareMathOperator*{\argmin}{arg\,min}
\DeclareMathOperator*{\argmax}{arg\,max}
\DeclareMathOperator*{\essinf}{ess\,inf}
\DeclareMathOperator*{\esssup}{ess\,sup}

\begin{document}
\title{Convex Regularization for High-Dimensional Multi-Response Tensor Regression}

\author{Garvesh Raskutti$^\ast,$ Ming Yuan$^\dag$ and Han Chen$^\dag$\\
University of Wisconsin-Madison}

\date{}

\footnotetext[1]{
Departments of Statistics and Computer Science, and Optimization Group at Wisconsin Institute for Discovery, University of Wisconsin-Madison, 1300 University Avenue, Madison, WI 53706. The research of Garvesh Raskutti is supported in part by NSF Grant DMS-1407028}
\footnotetext[2]{
Morgridge Institute for Research and Department of Statistics, University of Wisconsin-Madison, 1300 University Avenue, Madison, WI 53706. The research of Ming Yuan and Han Chen was supported in part by NSF FRG Grant DMS-1265202, and NIH Grant 1-U54AI117924-01.}
\maketitle

\begin{abstract}
In this paper we present a general convex optimization approach for solving high-dimensional multiple response tensor regression problems under low-dimensional structural assumptions. We consider using convex and weakly decomposable regularizers assuming that the underlying tensor lies in an unknown low-dimensional subspace. Within our framework, we derive general risk bounds of the resulting estimate under fairly general dependence structure among covariates. Our framework leads to upper bounds in terms of two very simple quantities, the \emph{Gaussian width} of a convex set in tensor space and the \emph{intrinsic dimension} of the low-dimensional tensor subspace. To the best of our knowledge, this is the first general framework that applies to multiple response problems. These general bounds provide useful upper bounds on rates of convergence for a number of fundamental statistical models of interest including multi-response regression, vector auto-regressive models, low-rank tensor models and pairwise interaction models. Moreover, in many of these settings we prove that the resulting estimates are minimax optimal. We also provide a numerical study that both validates our theoretical guarantees and demonstrates the breadth of our framework.
\end{abstract}

\newpage

\section{Introduction}

Many modern scientific problems involve solving high-dimensional statistical problems where the sample size is small relative to the ambient dimension of the underlying parameter to be estimated. Over the past few decades there has been a large amount of work on solving such problems by imposing low-dimensional structure on the parameter of interest. In particular sparsity, low-rankness and other low-dimensional subspace assumptions have been studied extensively both in terms of the development of fast algorithms and theoretical guarantees. See, e.g., \cite{BuhlmannVDGBook} and \cite{HastieTibshiraniWainwrightBook}, for an overview. Most of the prior work has focussed on scenarios in which the parameter of interest is a vector or matrix. Increasingly common in practice, however, the parameter or object to be estimated naturally has a higher order tensor structure. Examples include hyper-spectral image analysis \citep{LiLi10}, multi-energy computed tomography \citep{Semerci14}, radar signal processing \citep{SidNion10}, audio classification \citep{Mesgarani06} and text mining \citep{CohenCollins12} among numerous others. It is much less clear how the low dimensional structures inherent to these problems can be effectively accounted for. The main purpose of this article is to fill in this void and provide a general and unifying framework for doing so.

Consider a general tensor regression problem where covariate tensors $X^{(i)}\in \RR^{d_1\times\cdots\times d_M}$ and response tensors $Y^{(i)}\in \RR^{d_{M+1}\times\cdots\times d_N}$ are related through:
\begin{equation}
\label{eq:model}
Y^{(i)}=\langle X^{(i)}, T\rangle+\epsilon^{(i)}, \qquad i=1,2,\ldots,n.
\end{equation}
Here $T\in \RR^{d_1\times \cdots \times d_N}$ is an unknown parameter of interest, and $\epsilon^{(i)}$s are independent and identically distributed noise tensors whose entries are independent and identically distributed centred normal random variables with variance $\sigma^2$. Further, for simplicity we assume the covariates $(X^{(i)})_{i=1}^n$ are Gaussian, but with fairly general dependence assumptions. The notation $\langle\cdot,\cdot \rangle$ will refer throughout this paper to the standard inner product taken over appropriate Euclidean spaces. Hence, for $A\in \RR^{d_1\times\cdots\times d_M}$ and $B\in \RR^{d_1\times\cdots\times d_N}$:
$$\langle A , B \rangle=\sum_{j_1 = 1}^{d_1}\cdots\sum_{j_M = 1}^{d_M}A_{j_1,\ldots,j_M}B_{j_1,\ldots,j_M}\in \RR$$
is the usual inner product if $M=N$; and if $M<N$, then $\langle A , B \rangle\in \RR^{d_{M+1}\times\cdots\times d_N}$ such that its $(j_{M+1},\ldots, j_N)$ entry is given by
$$\left(\langle A , B \rangle\right)_{j_{M+1},\ldots, j_N}=\sum_{j_1 = 1}^{d_1}\cdots\sum_{j_M = 1}^{d_M}{A_{j_1,\ldots,j_M}B_{j_1,\ldots,j_M,j_{M+1},\ldots,j_N}}.$$
The goal of tensor regression is to estimate the coefficient tensor $T$ based on observations $\{(X^{(i)},Y^{(i)}): 1\le i\le n\}$. In addition to the canonical example of tensor regression with $Y$ a scalar response (i.e., $M=N$), many other commonly encountered regression problems are also special cases of the general tensor regression model (\ref{eq:model}). Multi-response regression \citep[see, e.g.,][]{AndersonStat}, vector autoregressive model \citep[see, e.g.,][]{Lut06}, and pairwise interaction tensor model \citep[see, e.g.,][]{RendleSchmidt10} are some of the notable examples. In this article, we provide a general treatment to these seemingly different problems.

Our main focus here is on situations where the dimensionality $d_k$s are large when compared with the sample size $n$. In many practical settings, the true regression coefficient tensor $T$ may have certain types of low-dimensional structure. Because of the high ambient dimension of a regression coefficient tensor, it is essential to account for such a low-dimensional structure when estimating it. Sparsity and low-rankness are the most common examples of such low dimensional structures. In the case of tensors, sparsity could occur at the entry-wise level, fiber-wise level, or slice-wise level, depending on the context and leading to different interpretations. There are also multiple ways in which low-rankness may be present when it comes to higher order tensors, either at the original tensor level or at the \emph{matricized} tensor level. 

In this article, we consider a general class of convex regularization techniques to exploit either type of low-dimensional structure. In particular, we consider the standard convex regularization framework:
\begin{equation}
\label{EqnGeneral}
\widehat{T} \in \argmin_{A\in \RR^{d_1\times \cdots\times d_N}}\left\{\frac{1}{2n} \sum_{i=1}^{n} \|Y^{(i)} - \langle A, X^{(i)} \rangle \|_{\rm F}^2 + \lambda \mathcal{R}(A)\right\},
\end{equation}
where the regularizer $\mathcal{R}(\cdot)$ is a norm on $\RR^{d_1\times \cdots\times d_N}$, and $\lambda>0$ is a tuning parameter. Hereafter, for a tensor $A$, $\|A\|_{\rm F}=\langle A, A\rangle^{1/2}$. We derive general risk bounds for a family of so-called \emph{weakly decomposable} regularizers under fairly general dependence structure among the covariates. These general upper bounds apply to a number of concrete statistical inference problems including the aforementioned multi-response regression, high-dimensional vector auto-regressive models, low-rank tensor models, and pairwise interaction tensors where we show that they are typically optimal in the minimax sense.

In developing these general results, we make several contributions to a fast growing literature on high dimensional tensor estimation. First of all, we provide a unified and principled approach to exploit the low dimensional structure in these tensor problems. In doing so, we incorporate an extension of the notion of decomposability originally introduced by \cite{Neg10} for vector and matrix models to \emph{weak decomposability} previously introduced in~\cite{vandeGeer14} which allows us to handle more delicate tensor models such as the nuclear norm regularization for low-rank tensor models. Moreover, we provide, for the regularized least squared estimate given by (\ref{EqnGeneral}), a general risk bound under an easily interpretable condition on the design tensor. The risk bound we derive is presented in terms of merely two geometric quantities, the \emph{Gaussian width} which depends on the choice of regularization and the \emph{intrinsic dimension} of the subspace that the tensor $T$ lies in. We believe this is the first general framework that applies to multiple responses and general dependence structure for the covariate tensor $X$. Finally, our general results lead to novel upper bounds for several important regression problems involving high-dimensional tensors: multi-response regression, multi-variate auto-regressive models and pairwise interaction models, for which we also prove that the resulting estimates are minimiax rate optimal with appropriate choices of regularizers.

Our framework incorporates both tensor structure and multiple responses which present a number of challenges compared to previous approaches. These challenges manifest themselves both in terms of the choice of regularizer $\mathcal{R}$ and the technical challenges in the proof of the main result. Firstly since the notion of low-dimensional is more generic for tensors meaning there are a number of choices of convex regularizer $\mathcal{R}$ and these must satisfy a form of weak decomposability and provide optimal rates. Multiple responses and the flexible dependence structure among the covariates also present significant technical challenges for proving restricted strong convexity, a key technical tool for establishing rates of convergence. In particular, a one-sided uniform law (Lemma~\ref{LemEmpLowerBound}) is required instead of classical techniques as developed in e.g.~\cite{NegWain12,RasWaiYu10b} that only apply to univariate responses.

The remainder of the paper is organized as follows: In Section~\ref{SecProbSetup} we introduce the general framework of using weakly decomposable regularizers for exploiting low-dimensional structures in high dimensional tensor regression. In Section~\ref{SecBounds} we present a general upper bound for weakly decomposable regularizers and discuss specific risk bounds for commonly used sparsity or low-rankness regularizers for tensors. In Section~\ref{SecExamples} we apply our general result to three specific statistical problems, namely, multi-response regression, multivariate autoregressive model, and the pairwise interaction model. We show that in each of the three examples appropriately chosen weakly decomposable regularizers leads to minimax optimal estimation of the unknown parameters. Numerical experiments are presented in Section~\ref{SecSim} to further demonstrate the merits and breadth of our approach. Proofs are provided in Section~\ref{SecProofs}.

\section{Methodology}
\label{SecProbSetup}

Recall that the regularized least-squares estimate is given by
\begin{equation*}
\widehat{T}=\argmin_{A\in \RR^{d_1\times \cdots\times d_N}}\left\{\frac{1}{2n} \sum_{i=1}^{n} \|Y^{(i)} - \langle A, X^{(i)} \rangle \|_{\rm F}^2 + \lambda \mathcal{R}(A)\right\}.
\end{equation*}
For brevity, we assume implicitly hereafter that the minimizer on its left hand side is uniquely defined. Our development here actually applies to the more general case where $\widehat{T}$ can be taken as an arbitrary element from the set of the minimizers. Of particularly interest here is the so-called \emph{weakly decomposable} convex regularizers, extending a similar concept introduced by \cite{Neg10} for vectors and matrices.

Let $\mathcal{A}$ be an arbitrary linear subspace of $\RR^{d_1 \times \cdots \times d_N}$ and $\mathcal{A}^{\perp}$ its orthogonal complement:
$$\mathcal{A}^{\perp} := \{A \in \mathbb{R}^{d_1 \times \cdots\times d_N}\;| \; \langle A, B \rangle = 0\; \mbox{for all}\; B \in \mathcal{A}\}.$$ 
We call a regularizer $\mathcal{R}(\cdot)$ weakly decomposable with respect to a pair $(\mathcal{A}, \calB)$ where $\calB\subseteq \calA$ if there exist a constant $0<c_\calR\le 1$ such that for any $A \in \mathcal{A}^\perp$ and $B \in \mathcal{B}$,
\begin{equation}
\label{eq:decomcond}
\mathcal{R}(A + B)\ge \mathcal{R}(A) + c_{\calR}\mathcal{R}(B).
\end{equation}
In particular, if (\ref{eq:decomcond}) holds for any $B\in \calB=\calA$, we say $\mathcal{R}(\cdot)$ is weakly decomposable with respect to $\calA$. A more general version of this concept was first introduced in~\cite{vandeGeer14}. Because $\mathcal{R}$ is a norm, by triangular inequality, we also have
$$
\mathcal{R}(A + B) \le \mathcal{R}(A) + \mathcal{R}(B). 
$$
Many of the commonly used regularizers for tensors are weakly decomposable, or decomposable for short. When $c_\calR=1$, our definition of decomposability naturally extends from similar notion for vectors ($N=1$) and matrices ($N=2$) introduced by \cite{Neg10}. We also allow for more general choices of $c_\calR$ here to ensure a wider applicability. For example as we shall see the popular tensor nuclear norm regularizer is decomposable with respect to appropriate linear subspaces with $c_\calR=1/2$, but not decomposable if $c_\calR=1$. 

We have now described a catalogue of commonly used regularizers for tensors and argue that they are all decomposable with respect to appropriately chosen subspaces of $\RR^{d_1\times \cdots\times d_N}$. To fix ideas, we shall focus in what follows on estimating a third-order tensor $T$, that is $N=3$, although our discussion can be straightforwardly extended to higher-order tensors. 

\subsection{Sparsity regularizers}

An obvious way to encourage entry-wise sparsity is to impose the vector $\ell_1$ penalty on the entries of $A$:
\begin{equation}
\label{eq:lasso}
\calR(A):=\sum_{j_1=1}^{d_1}\sum_{j_2=1}^{d_2}\sum_{j_3=1}^{d_3} |A_{j_1j_2j_3}|,
\end{equation}
following the same idea as the Lasso for linear regression \citep[see, e.g.,][]{Tibshirani96}. This is a canonical example of decomposable regularizers. For any fixed $I\subset [d_1]\times[d_2]\times[d_3]$ where $[d]=\{1,2,\ldots,d\}$, write
\begin{equation}
\label{eq:defA1}
\mathcal{A}(I)=\mathcal{B}(I)=\left\{A\in \mathbb{R}^{d_1 \times d_2\times d_3}: A_{j_1j_2j_3}=0 {\rm \ for\ all\ } (j_1,j_2,j_3)\notin I\right\}.
\end{equation}
It is clear that
$$
\mathcal{A}^\perp(I)=\left\{A\in \mathbb{R}^{d_1 \times d_2\times d_3}: A_{j_1j_2j_3}=0 {\rm \ for\ all\ } (j_1,j_2,j_3)\in I\right\},
$$
and $\calR(A)$ defined by (\ref{eq:lasso}) is decomposable with respect to $\calA$ with $c_{\calR}=1$.

In many applications, sparsity arises with a more structured fashion for tensors. For example, a fiber or a slice of a tensor is likely to be zero simultaneously. Mode-$1$ fibers of a tensor $A\in \calR^{d_1\times d_2\times d_3}$ are the collection of $d_1$-dimensional vectors
$$\left\{A_{\cdot j_2j_3}=(A_{1j_2j_3},\ldots,A_{d_1j_2j_3})^\top: 1\le j_2\le d_2, 1\le j_3\le d_3\right\}.$$
Mode-$2$ and -$3$ fibers can be defined in the same fashion. To fix ideas, we focus on mode-$1$ fibers. Sparsity among mode-$1$ fibers can be exploited using the group-based $\ell_1$ regularizer:
\begin{equation}
\label{eq:groupLasso}
\mathcal{R}(A) = \sum_{j_2=1}^{d_2}\sum_{j_3=1}^{d_3} \|A_{\cdot j_2j_3}\|_{\ell_2},
\end{equation} 
similar to the group Lasso \citep[see, e.g.,][]{YuaLi06}, where $\|\cdot\|_{\ell_2}$ stands for the usual vector $\ell_2$ norm. Similar to the vector $\ell_1$ regularizer, the group $\ell_1$-based regularizer is also decomposable. For any fixed $I\subset [d_2]\times[d_3]$, write
\begin{equation}
\label{eq:defA2}
\mathcal{A}(I)=\mathcal{B}(I)=\left\{A\in \mathbb{R}^{d_1 \times d_2\times d_3}: A_{j_1j_2j_3}=0 {\rm \ for\ all\ } (j_2,j_3)\notin I\right\}.
\end{equation}
It is clear that
$$
\mathcal{A}^\perp(I)=\left\{A\in \mathbb{R}^{d_1 \times d_2\times d_3}: A_{j_1j_2j_3}=0 {\rm \ for\ all\ } (j_2,j_3)\in I\right\},
$$
and $\calR(A)$ defined by (\ref{eq:groupLasso}) is decomposable with respect to $\calA$ with $c_{\calR}=1$. Note that in defining the regularizer in (\ref{eq:groupLasso}), instead of vector $\ell_2$ norm, other $\ell_q$ ($q>1$) norms could also be used. See, e.g., \cite{Turlach05}.

Sparsity could also occur at the slice level. The $(1,2)$ slices of a tensor $A\in \RR^{d_1\times d_2\times d_3}$ are the collection of $d_1\times d_2$ matrices
$$\left\{A_{\cdot \cdot j_3}=(A_{j_1j_2j_3})_{1\le j_1\le d_1,1\le j_2\le d_2}: 1\le j_3\le d_3\right\}.$$
Let $\|\cdot\|$ be an arbitrary norm on $d_1\times d_2$ matrices. Then the following group regularizer can be considered:
\begin{equation}
\label{eq:groupLassoMat}
\calR(A)=\sum_{j_3=1}^{d_3}\|A_{\cdot \cdot j_3}\|.
\end{equation}
Typical examples of the matrix norm that can be used in (\ref{eq:groupLassoMat}) include Frobenius norm and nuclear norm among others. In the case when $\|\cdot\|_{\rm F}$ is used, $\calR(\cdot)$ is again a decomposable regularizer with respect to
\begin{equation}
\label{eq:defA3}
\mathcal{A}(I)=\mathcal{B}(I)=\left\{A\in \mathbb{R}^{d_1 \times d_2\times d_3}: A_{j_1j_2j_3}=0 {\rm \ for\ all\ } j_3\notin I\right\}.
\end{equation}
for any $I\subset[d_3]$.

Now consider the case when we use the matrix nuclear norm $\|\cdot\|_\ast$ in (\ref{eq:groupLassoMat}). Let $P_{1j}$ and $P_{2j}$, $j=1,\ldots, d_3$ be two sequences of projection matrices on $\RR^{d_1}$ and $\RR^{d_2}$ respectively. Let
\begin{equation}
\label{eq:defA4}
\calA(P_{1j},P_{2j}: 1\le j\le d_3)=\left\{A\in \RR^{d_1\times d_2\times d_3}: P_{1j}^\perp A_{\cdot\cdot j}P_{2j}^\perp=0, j=1,\ldots, d_3\right\},
\end{equation}
and
\begin{equation}
\label{eq:defB4}
\calB(P_{1j},P_{2j}: 1\le j\le d_3)=\left\{A\in \RR^{d_1\times d_2\times d_3}: A_{\cdot\cdot j}=P_{1j}A_{\cdot\cdot j}P_{2j}, j=1,\ldots, d_3\right\}.
\end{equation}
By pinching inequality \citep[see, e.g.,][]{Bhatia97}, it can be derived that $\calR(\cdot)$ is decomposable with respect to $\calA(P_{1j},P_{2j}: 1\le j\le d_3)$ and $\calB(P_{1j},P_{2j}: 1\le j\le d_3)$.

\subsection{Low-rankness regularizers}
In addition to sparsity, one may also consider tensors with low-rank. There are multiple notions of rank for higher-order tensors. See, e.g., \cite{KoldarBader}, for a recent review. In particular, the so-called CP rank is defined as the smallest number $r$ of rank-one tensors needed to represent a tensor $A\in \RR^{d_1\times d_2\times d_3}$:
\begin{equation}
\label{eq:cpdecomp}
A=\sum_{k=1}^r u_k\otimes v_k\otimes w_k
\end{equation}
where $u_k\in \RR^{d_1}$, $v_k\in \RR^{d_2}$ and $w_k\in \RR^{d_3}$. To encourage a low rank estimate, we can consider the nuclear norm regularization. Following \cite{YuanZhang14}, we define the nuclear norm of $A$ through its dual norm. More specifically, let the spectral norm of $A$ be given by
$$
\|A\|_{s}=\max_{\|u\|_{\ell_2},\|v\|_{\ell_2},\|w\|_{\ell_2}\le 1}\langle A, u\otimes v\otimes w\rangle.
$$
Then its nuclear norm is defined as
$$
\|A\|_\ast=\max_{\|B\|_s\le 1}\langle A, B\rangle.
$$
We shall then consider the regularizer:
\begin{equation}
\label{eq:nuclear}
\mathcal{R}(A) = \|A\|_\ast.
\end{equation}
We now show this is also a weakly decomposable regularizer.

Let $P_k$ be a projection matrix in $\RR^{d_k}$. Define
$$
(P_1\otimes P_2\otimes P_3) A=\sum_{k=1}^r P_1u_k\otimes P_2v_k\otimes P_3w_k.
$$
Write
$$
Q=P_1\otimes P_2\otimes P_3+P_1^\perp\otimes P_2\otimes P_3+P_1\otimes P_2^\perp\otimes P_3+P_1\otimes P_2\otimes P_3^\perp,
$$
and
$$Q^\perp =P_1^\perp\otimes P_2^\perp\otimes P_3^\perp+P_1^\perp\otimes P_2^\perp\otimes P_3+P_1\otimes P_2^\perp\otimes P_3^\perp+P_1^\perp\otimes P_2\otimes P_3^\perp,$$
where $P_k^\perp=I-P_k$.
\blems
\label{le:pinch}
For any $A\in \RR^{d_1\times d_2\times d_3}$ and projection matrices $P_k$ in $\RR^{d_k}$, $k=1,2,3$, we have
$$
\|A\|_{\ast}\ge \|(P_1\otimes P_2\otimes P_3) A\|_\ast+{1\over 2}\|Q^\perp A\|_\ast.
$$
\elems
Lemma \ref{le:pinch} is a direct consequence from the characterization of sub-differential for tensor nuclear norm given by \cite{YuanZhang14}, and can be viewed as a tensor version of the pinching inequality for matrices.

Write
\begin{equation}
\label{eq:defA5}
\calA(P_1,P_2,P_3)=\left\{A\in \RR^{d_1\times d_2\times d_3}: Q A=A\right\},
\end{equation}
and
\begin{equation}
\label{eq:defB5}
\calB(P_1,P_2,P_3)=\left\{A\in \RR^{d_1\times d_2\times d_3}: (P_1\otimes P_2\otimes P_3) A=A\right\}.
\end{equation}
By Lemma \ref{le:pinch}, $\calR(\cdot)$ defined by (\ref{eq:nuclear}) is weakly decomposable with respect to $\calA(P_1,P_2,P_3)$ and $\calB(P_1,P_2,P_3)$ with $c_{\calR}=1/2$. We note that a counterexample is also given by \cite{YuanZhang14} which shows that, for the tensor nuclear norm, we cannot take $c_\calR=1$.

Another popular way to define tensor rank is through the so-called Tucker decomposition. Recall that the Tucker decomposition of a tensor $A\in \RR^{d_1\times d_2\times d_3}$ is of the form:
\begin{equation}
\label{eq:tucker}
A_{j_1j_2j_3}=\sum_{k_1=1}^{r_1}\sum_{k_2=1}^{r_2}\sum_{k_3=1}^{r_3} S_{k_1k_2k_3}U_{j_1k_1}V_{j_2k_2}W_{j_3k_3}
\end{equation}
so that $U$, $V$ and $W$ are orthogonal matrices, and the so-called core tensor $S=(S_{k_1k_2k_3})_{k_1,k_2,k_3}$ is such that any two slices of $S$ are orthogonal. The triplet $(r_1,r_2,r_3)$ are referred to as the Tucker ranks of $A$. It is not hard to see that if (\ref{eq:cpdecomp}) holds, then the Tucker ranks $(r_1,r_2,r_3)$ can be equivalently interpreted as the dimensionality of the linear spaces spanned by $\{u_k: 1\le k\le r\}$, $\{v_k: 1\le k\le r\}$, and $\{w_k: 1\le k\le r\}$ respectively. The following relationship holds between CP rank and Tucker ranks:
$$
\max\{r_1,r_2,r_3\}\le r\le \min\{r_1r_2,r_2r_3,r_1r_3\}.
$$

A convenient way to encourage low Tucker ranks in a tensor is through matricization. Let $\calM_1(\cdot)$ denote the mode-$1$ matricization of a tensor. That is $\calM_1(A)$ is a $d_1\times (d_2d_3)$ matrix whose column vectors are the the mode-$1$ fibers of $A\in \RR^{d_1\times d_2\times d_3}$. $\calM_2(\cdot)$ and $\calM_3(\cdot)$ can also be defined in the same fashion. It is clear
$$
{\rm rank}(\calM_k(A))=r_k(A).
$$
A natural way to encourage low-rankness is therefore through nuclear norm regularization:
\begin{equation}
\label{eq:matricize}
\mathcal{R}(A) = {1\over 3}\sum_{k=1}^3 \|\calM_k(A)\|_{\ast}.
\end{equation}
By the pinching inequality for matrices, $\calR(\cdot)$ defined by (\ref{eq:matricize}) is also decomposable with respect to $\calA(P_1,P_2,P_3)$ and $\calB(P_1,P_2,P_3)$ with $c_{\calR}=1$.

\section{Risk Bounds for Decomposable Regularizers}
\label{SecBounds}

We now establish risk bounds for general decomposable regularizers. In particular, our bounds are given in terms of the \emph{Gaussian width} of a suitable set of tensors. Recall that the Gaussian width of a set $S \subset \mathbb{R}^{d_1 \times d_2 \times...\times d_N}$ is given by
\begin{equation*}
w_G(S) := \mathbb{E}\left(\sup_{A \in S} \langle A, G \rangle \right),
\end{equation*}
where $G \in \mathbb{R}^{d_1 \times d_2 \times ... \times d_N}$ is a tensor whose entries are independent $\mathcal{N}(0,1)$ random variables. See, e.g.~\cite{Gordon88} for more details on Gaussian width.
 
Note that the Gaussian width is a geometric measure of the volume of the set $S$ and can be related to other volumetric characterizations \citep[see, e.g.,][]{Pisier89}.We also define the unit ball for the norm-regularizer $\mathcal{R}(.)$ as follows:
$$\mathbb{B}_{\mathcal{R}}(1) := \{A \in \mathbb{R}^{d_1 \times d_2 \times...\times d_N}\;|\; \mathcal{R}(A) \leq 1 \}.$$ We impose the mild assumption that $\|A\|_{\rm F} \leq \mathcal{R}(A)$ which ensures that the regularizer $\mathcal{R}(\cdot)$ encourages low-dimensional structure.

Now we define a quantity that relates the size of the norm $\mathcal{R}(A)$ to the Frobenius norm $\|A\|_{\rm F}$ over the the low-dimensional subspace $\calA$. Following \cite{Neg10}, for a subspace $\calA$ of $\RR^{d_1\times \cdots\times d_N}$, define its compatibility constant $s(\mathcal{A})$ as
\begin{equation*}
s(\mathcal{A}) := \sup_{A \in \mathcal{A}/\{0\}} \frac{\mathcal{R}^2(A)}{\|A\|_{\rm F}^2},
\end{equation*}
which can be interpreted as a notion of intrinsic dimensionality of $\calA$.

Now we turn our attention to the covariate tensor. Denote by $X^{(i)}={\rm vec}(X^{(i)})$ the vectorized covariate from the $i$th sample. With slight abuse of notation, write
$$X={\rm vec}((X^{(1)})^\top,\ldots, (X^{(n)})^\top) \in \mathbb{R}^{n.d_1d_2\cdots d_M}$$
the concatenated covariates from all $n$ samples. For convenience let $D_M = d_1d_2\cdots d_M$. Further for brevity we assume a Gaussian design so that
$$X \sim \mathcal{N}(0,\Sigma)$$ where $$\Sigma = {\rm cov}(X, X) \in \mathbb{R}^{nD_M \times nD_M}.$$
With more technical work our results may be extended beyond Gaussian designs. We note that we do not require that the sample tensors $X^{(i)}$ be independent.

We shall assume that $\Sigma$ has bounded eigenvalues which we later verify for a number of statistical examples. Let $\lambda_{\min}(\cdot)$ and $\lambda_{\max}(\cdot)$ represent the smallest and largest eigenvalues of a matrix, respectively. In what follows, we shall assume that
\begin{equation}
\label{AssCov}
c_{\ell}^2 \leq \lambda_{\min}(\Sigma) \leq \lambda_{\max}(\Sigma) \leq c_u^2,
\end{equation}
for some constants $0< c_\ell\le c_u<\infty$.

Note that in particular if all covariates $\{X^{(i)}: i=1,\ldots,n\}$ are independent and identically distributed, then $\Sigma$ has a block diagonal structure, and (\ref{AssCov}) boils down to similar conditions on ${\rm cov}(X^{(i)},X^{(i)})$. However (\ref{AssCov}) is more general and applicable to settings in which the $X^{(i)}$'s may be dependent such as time-series models, which we shall discuss in further detail in Section~\ref{SecExamples}.

We are now in position to state our main result on the risk bounds in terms of both Frobenius norm $\|\cdot\|_{\rm F}$ and the empirical norm $\|\cdot\|_n$ where for a tensor $A\in \RR^{d_1\times \cdots \times d_N}$, which we define as:
$$
\|A\|_n^2 := {1\over n}\sum_{i=1}^n \|\langle A, X^{(i)}\rangle\|_{F}^2.
$$
The main reason we focus on random Gaussian design is so that we can prove a one-sided uniform law that relates the empirical norm defined above to the Frobenius norm of a tensor in $\mathcal{A}$ (see Lemma~\ref{LemEmpLowerBound}). Lemma~\ref{LemEmpLowerBound} is analogous to restricted strong convexity defined in~\cite{Neg10} but since we are dealing with multiple responses a more refined technique is required to prove Lemma~\ref{LemEmpLowerBound}. 

\btheos
\label{ThmUpper}
Suppose that (\ref{eq:model}) holds for a tensor $T$ from a linear subspace $\mathcal{A}_0\subset \RR^{d_1\times \cdots\times d_N}$ where (\ref{AssCov}) holds. Let $\widehat{T}$ be defined by~\eqref{EqnGeneral} where the regularizer $\calR(\cdot)$ is decomposable with respect to $\calA$ and $\calA_0$ for some linear subspace $\calA\supseteq\calA_0$. If
\begin{equation}
\lambda \geq \frac{2 \sigma c_u (3+c_\calR)}{c_\calR\sqrt{n}}w_G[\mathbb{B}_{\calR}(1)],
\end{equation}
then there exists a constant $c>0$ such that with probability at least $1 - \exp\{-c w_G^2[\mathbb{B}_{\calR}(1)]\}$,
\begin{equation}
\label{eq:risk}
\max\left\{\|\widehat{T}-T\|_n^2, \|\widehat{T}-T\|_{\rm F}^2\right\} \leq \frac{6(1+c_\calR)}{3+c_\calR}\frac{9 c_u^2}{c_{\ell}^2} s(\mathcal{A}) \lambda^2,
\end{equation}
when $n$ is sufficiently large, assuming that the right hand side converges to zero as $n$ increases.
\etheos
As stated in Theorem \ref{ThmUpper}, our upper bound boils down to bounding two quantities, $s(\mathcal{A})$ and $w_G[\mathbb{B}_{\calR}(1)]$ which are both purely geometric quantities. To provide some intuition, $w_G[\mathbb{B}_{\calR}(1)]$ captures how large the $\calR(\cdot)$ norm is relative to the $\|\cdot\|_{\rm F}$ norm and $s(\mathcal{A})$ captures the low dimension of the subspace $\mathcal{A}$.

Several technical remarks are in order. Note that $w_G[\mathbb{B}_{\calR}(1)]$ can be expressed as expectation of the \emph{dual norm} of $G$. According to  $\mathcal{R}$ \citep[see, e.g.,][for details]{Rockafellar}, the dual norm $\mathcal{R}^*(\cdot)$ is given by:
\begin{equation*}
\mathcal{R}^*(B) := \sup_{A \in \mathbb{B}_{\calR}(1)} \langle A, B \rangle,
\end{equation*}
where the supremum is taken over tensors of the same dimensions as $B$. 
It is straightforward to see that $w_G[\mathbb{B}_{\calR}(1)] = \mathbb{E}[\mathcal{R}^*(G)]$.

To the best of our knowledge, this is the first general result that applies to multiple responses. As mentioned earlier, incorporating multiple responses presents a technical challenge (see Lemma~\ref{LemEmpLowerBound}) which is a one-sided uniform law analogous to restricted strong convexity. While Theorem~\ref{ThmUpper} focusses on Gaussian design, results can be extended to random  sub-Gaussian design using more sophisticated techniques \citep[see e.g.,][]{Mendelson15,Zhou09} or for fixed design by assuming covariates deterministically satisfy the condition in Lemma~\ref{LemEmpLowerBound}. Since the focus of this paper is on general dependence structure, we assume random Gaussian design. 

One important practical challenge is that $\sigma^2$, $c_u$ and $c_{\ell}$ are typically unknown and these clearly influence the choice of $\lambda$. This is a common challenge for high-dimensional statistical inference and we don't address this issue in this paper. In practice, $\lambda$ is typically chosen through cross-validation. A more sophisticated choice of $\lambda$ based on estimation of $\sigma^2$ and other constants remains an open question. Another important and open question is for what choices of  $\mathcal{A}_0$ is the upper bound optimal (up to a constant). In Section~\ref{SecExamples} we provide specific examples in which we provide minimax lower bounds which match the upper bounds up to constant. However as we see for low-rank tensor regression for low-rank tensor regression discussed in Section 3.2, we are not aware of a convex regularizer that matches the minimax lower bound. 

Now we develop upper bounds on both quantities in different scenarios. As in the previous section, we shall focus on third order tensor in the rest of the section for the ease of exposition.

\subsection{Sparsity regularizers}
We first consider sparsity regularizers described in the previous section.

\subsubsection{Entry-wise and fiber-wise sparsity}
Recall that vectorized $\ell_1$ regularizer:
$$
\calR_1(A)=\sum_{j_1=1}^{d_1}\sum_{j_2=1}^{d_2}\sum_{j_3=1}^{d_3}|A_{j_1j_2j_3}|,
$$
could be used to exploit entry-wise sparsity. Clearly,
$$
\calR^\ast_1(A)=\max_{j_1,j_2,j_3}|A_{j_1j_2j_3}|.
$$
It can then be shown that:
\blems
\label{LemSparsity}
There exists a constant $0 < c < \infty$ such that
\begin{equation}
\label{eq:lassog}
w_G[\mathbb{B}_{\calR_1}(1)] \leq c\sqrt{\log(d_1d_2d_3)}.
\end{equation}
\elems
Let
$$
\Theta_1(s)=\left\{A\in \RR^{d_1\times d_2\times d_3}: \sum_{j_1=1}^{d_1}\sum_{j_2=1}^{d_2}\sum_{j_3=1}^{d_3}{\mathbb I}(A_{j_1j_2j_3}\neq 0)\le s\right\}.
$$
For an arbitrary $A\in \Theta_1(s)$, write
$$
I(A)=\left\{(j_1,j_2,j_3)\in [d_1]\times [d_2]\times [d_3]: A_{j_1j_2j_3}\neq 0\right\}.
$$
Then $\calR_1(\cdot)$ is decomposable with respect to $\calA(I(A))$ as defined by (\ref{eq:defA1}). It is easy to verify that for any $A\in \Theta_1(s)$,
\begin{equation}
\label{eq:lassos}
s_1(\calA(I))=\sup_{B \in \mathcal{A}(I(A))/\{0\}} \frac{\mathcal{R}_1^2(B)}{\|B\|_{\rm F}^2}\le s.
\end{equation}

In light of (\ref{eq:lassos}) and (\ref{eq:lassog}), Theorem \ref{ThmUpper} implies that
$$
\sup_{T\in \Theta_1(s)}\max\left\{\|\widehat{T}_1-T\|_n^2, \|\widehat{T}_1-T\|_{\rm F}^2\right\} \lesssim {s\log(d_1d_2d_3)\over n},
$$
with high probability by taking
$$
\lambda\asymp \sqrt{\log(d_1d_2d_3)\over n},
$$
where $\widehat{T}_1$ is the regularized least squares estimate defined by (\ref{EqnGeneral}) when using regularizer $\calR_1(\cdot)$.

A similar argument can also be applied to fiber-wise sparsity. To fix ideas, we consider here only sparsity among mode-1 fibers. In this case, we use a group Lasso type of regularizer:
$$
\calR_2(A)=\sum_{j_2=1}^{d_2}\sum_{j_3=1}^{d_3}\|A_{\cdot j_2j_3}\|_{\ell_2}.
$$
Then
$$
\calR^\ast_2(A)=\max_{j_2,j_3}\|A_{\cdot j_2j_3}\|_{\ell_2}.
$$

\blems
\label{LemSparsityfiber}
There exists a constant $0 < c < \infty$ such that
\begin{equation}
\label{eq:glassog}
w_G[\mathbb{B}_{\calR_2}(1)]\leq c \sqrt{ \max\{d_1,\log(d_2d_3)\}}.
\end{equation}
\elems

Let
$$
\Theta_2(s)=\left\{A\in \RR^{d_1\times d_2\times d_3}: \sum_{j_2=1}^{d_2}\sum_{j_3=1}^{d_3}{\mathbb I}(A_{\cdot j_2j_3}\neq \zero)\le s\right\}.
$$
Similar to the previous case, for an arbitrary $A\in \Theta_1(s)$, write
$$
I(A)=\{(j_2,j_3)\in [d_2]\times [d_3]: A_{\cdot j_2j_3}\neq \zero\}.
$$
Then $\calR_2(\cdot)$ is decomposable with respect to $\calA(I(A))$ as defined by (\ref{eq:defA2}). It is easy to verify that for any $A\in \Theta_2(s)$,
\begin{equation}
\label{eq:glassos}
s_2(\calA(I))=\sup_{B \in \mathcal{A}(I(A))/\{0\}} \frac{\mathcal{R}_2^2(B)}{\|B\|_{\rm F}^2}\le s.
\end{equation}

In light of (\ref{eq:glassos}) and (\ref{eq:glassog}), Theorem \ref{ThmUpper} implies that
$$
\sup_{T\in \Theta_2(s)}\max\left\{\|\widehat{T}_2-T\|_n^2, \|\widehat{T}_2-T\|_{\rm F}^2\right\} \lesssim {s\max\{d_1,\log(d_2d_3)\}\over n},
$$
with high probability by taking
$$
\lambda\asymp \sqrt{\max\{d_1,\log(d_2d_3)\}\over n},
$$
where $\widehat{T}_2$ is the regularized least squares estimate defined by (\ref{EqnGeneral}) when using regularizer $\calR_2(\cdot)$.

Comparing with the rates for entry-wise and fiber-wise sparsity regularization, we can see the benefit of using group Lasso type of regularizer $\calR_2$ when sparsity is likely to occur at the fiber level. More specifically, consider the case when there are a total of $s_1$ nonzero entries from $s_2$ nonzero fibers. If an entry-wise $\ell_1$ regularization is applied, we can achieve the risk bound:
$$
\|\widehat{T}_1-T\|_{\rm F}^2\lesssim {s_1\log(d_1d_2d_3)\over n}.
$$
On the other hand, if fiber-wise group $\ell_1$ regularization is applied, then the risk bound becomes:
$$
\|\widehat{T}_2-T\|_{\rm F}^2\lesssim {s_2\max\{d_1,\log(d_2d_3)\}\over n}.
$$
When nonzero entries are clustered in fibers, we may expect $s_1\asymp s_2d_1$. In this case, $\widehat{T}_2$ enjoys performance superior to that of $\widehat{T}_1$ since $s_2 d_1 \log(d_1d_2d_3)$ is larger than $s_2\max\{d_1,\log(d_2d_3)\}$.

\subsubsection{Slice-wise sparsity and low-rank structure}

Now we consider slice-wise sparsity and low-rank structure. Again, to fix ideas, we consider here only sparsity among $(1,2)$ slices. As discussed in the previous section, two specific types of regularizers could be employed:
$$
\calR_3(A)=\sum_{j_3=1}^{d_3}\|A_{\cdot \cdot j_3}\|_{\rm F},
$$
and
$$
\calR_4(A)=\sum_{j_3=1}^{d_3}\|A_{\cdot \cdot j_3}\|_\ast,
$$
where recall that $\|.\|_{\ast}$ denotes the nuclear norm of a matrix, that is the sum of all singular values.

Note that
$$
\calR^\ast_3(A)=\max_{1\le j_3\le d_3}\|A_{\cdot \cdot j_3}\|_{\rm F}.
$$
Then we have the following result:
\blems
\label{LemSparsityslice}
There exists a constant $0 < c < \infty$ such that
\begin{equation}
\label{eq:sglassog}
w_G[\mathbb{B}_{\calR_3}(1)] \leq c \sqrt{ \max\{d_1d_2,\log(d_3)\}}.
\end{equation}
\elems

Let
$$
\Theta_3(s)=\left\{A\in \RR^{d_1\times d_2\times d_3}: \sum_{j_3=1}^{d_3}{\mathbb I}(A_{\cdot \cdot j_3}\neq \zero)\le s\right\}.
$$
For an arbitrary $A\in \Theta_1(s)$, write
$$
I(A)=\{j_3\in [d_3]: A_{\cdot \cdot j_3}\neq \zero\}.
$$
Then $\calR_3(\cdot)$ is decomposable with respect to $\calA(I(A))$ as defined by (\ref{eq:defA3}). It is easy to verify that for any $A\in \Theta_3(s)$,
\begin{equation}
\label{eq:sglassos}
s_3(\calA(I(A)))=\sup_{B \in \mathcal{A}(I(A))/\{0\}} \frac{\mathcal{R}_3^2(B)}{\|B\|_{\rm F}^2}\le s.
\end{equation}

Based on (\ref{eq:sglassos}) and (\ref{eq:sglassog}), Theorem \ref{ThmUpper} implies that
$$
\sup_{T\in \Theta_3(s)}\max\left\{\|\widehat{T}_3-T\|_n^2, \|\widehat{T}_3-T\|_{\rm F}^2\right\} \lesssim {s\max\{d_1d_2,\log(d_3)\}\over n},
$$
with high probability by taking
$$
\lambda\asymp \sqrt{\max\{d_1d_2,\log(d_3)\}\over n},
$$
where $\widehat{T}_3$ is the regularized least squares estimate defined by (\ref{EqnGeneral}) when using regularizer $\calR_3(\cdot)$.

Alternatively, for $\calR_4(\cdot)$,
$$
\calR^\ast_4(A)=\max_{j_3}\|A_{\cdot \cdot j_3}\|_{s},
$$
we have the following:
\blems
\label{LemSparsityslice2}
There exists a constant $0 < c < \infty$ such that
\begin{equation}
\label{eq:sglassog2}
w_G[\mathbb{B}_{\calR_4}(1)] \leq c \sqrt{ \max\{d_1,d_2,\log(d_3)\}}.
\end{equation}
\elems

Now consider
$$
\Theta_4(r)=\left\{A\in \RR^{d_1\times d_2\times d_3}: \sum_{j_3=1}^{d_3}{\rm rank}(A_{\cdot \cdot j_3})\le r\right\}.
$$
For an arbitrary $A\in \Theta_4(r)$, denote by $P_{1j}$ and $P_{2j}$ the projection onto the row and column space of $A_{\cdot \cdot j}$ respectively. It is clear that $A\in \calB(P_{1j}, P_{2j}: 1\le j\le d_3)$ as defined by (\ref{eq:defB4}). In addition, recall that $\calR_4$ is decomposable with respect to $\calB(P_{1j}, P_{2j}: 1\le j\le d_3)$ and $\calA(P_{1j}, P_{2j}: 1\le j\le d_3)$ as defined by (\ref{eq:defA4}). It is not hard to see that for any $A\in \Theta_4(r)$, $\calA(P_{1j}, P_{2j}: 1\le j\le d_3)\subset \Theta_4(2r)$, from which we can derive that:
\blems
\label{le:sglasso2}
For any $A\in \Theta_4(r)$,
\begin{equation}
\label{eq:sglassos2}
s_4(\calA(P_{1j}, P_{2j}: 1\le j\le d_3))\le\sup_{B\in \calA/\{0\}} \frac{\mathcal{R}_4^2(B)}{\|B\|_{\rm F}^2}\le 2r.
\end{equation}
\elems

In light of (\ref{eq:sglassos2}) and (\ref{eq:sglassog2}), Theorem \ref{ThmUpper} implies that
$$
\sup_{T\in \Theta_4(r)}\max\left\{\|\widehat{T}_4-T\|_n^2, \|\widehat{T}_4-T\|_{\rm F}^2\right\} \lesssim {r\max\{d_1,d_2,\log(d_3)\}\over n},
$$
with high probability by taking
$$
\lambda\asymp \sqrt{\max\{d_1,d_2,\log(d_3)\}\over n},
$$
where $\widehat{T}_4$ is the regularized least squares estimate defined by (\ref{EqnGeneral}) when using regularizer $\calR_4(\cdot)$.

Comparing with the rates for estimates with regularizers $\calR_3$ and $\calR_4$, we can see the benefit of using $\calR_4$ when the nonzero slices are likely to be of low-rank. In particular, consider the case when there are $s_1$ nonzero slices and each nonzero slice has rank up to $r$. Then applying $\calR_3$ leads to risk bound:
$$
\|\widehat{T}_{3}-T\|_{\rm F}^2\lesssim {s_1\max\{d_1d_2,\log(d_3)\}\over n},
$$
whereas applying $\calR_4$ leads to:
$$
\|\widehat{T}_{4}-T\|_{\rm F}^2\lesssim {s_1 r\max\{d_1,d_2,\log(d_3)\}\over n}.
$$
It is clear that $\widehat{T}_{4}$ is a better estimator when $r\ll d_1=d_2=d_3$.

\subsection{Low-rankness regularizers}

We now consider regularizers that encourages low rank estimates. We begin with the tensor nuclear norm regularization:
$$
\calR_5(A)=\|A\|_\ast.
$$
Recall that $\calR^\ast_5(A)=\|A\|_s$.
\blems
\label{LemNuclearNorm}
There exists a constant $0 < c < \infty$ such that
\begin{equation}
\label{eq:glassog}
w_G[\mathbb{B}_{\calR_5}(1)] \leq c \sqrt{(d_1 + d_2 + d_3)}.
\end{equation}
\elems

Now let
$$
\Theta_5(r)=\left\{A\in \RR^{d_1\times d_2\times d_3}: \max\{r_1(A), r_2(A), r_3(A)\} \le r\right\}.
$$
For an arbitrary $A\in \Theta_5(r)$, denote by $P_1$, $P_2$, $P_3$ the projection onto the linear space spanned by the mode-1, -2 and -3 fibers respectively. As we argued in the previous section, $\calR_5(\cdot)$ is weakly decomposable with respect to $\calA(P_1, P_2,P_3)$ and $\calB(P_1, P_2,P_3)$, and $A\in \calB(P_1, P_2,P_3)$ where $\calA(P_1, P_2,P_3)$ and $\calB(P_1, P_2,P_3)$ are defined by (\ref{eq:defA5}) and (\ref{eq:defB5}) respectively. 
\blems
For any $A\in \Theta_5(r)$,
\label{LemNuclearNorm1}
$$s_5(\calA(P_1, P_2,P_3))=\sup_{B \in \calA(P_1, P_2,P_3)/\{0\}} \frac{\mathcal{R}_5^2(B)}{\|B\|_{\rm F}^2}\le r^2.$$
\elems

Lemmas~\ref{LemNuclearNorm} and \ref{LemNuclearNorm1} show that
$$
\sup_{T\in \Theta_5(r)}\max\left\{\|\widehat{T}_5-T\|_n^2, \|\widehat{T}_5-T\|_{\rm F}^2\right\} \lesssim {r^2(d_1+d_2+d_3)\over n},
$$
with high probability by taking
$$
\lambda\asymp \sqrt{d_1+d_2+d_3\over n},
$$
where $\widehat{T}_5$ is the regularized least squares estimate defined by (\ref{EqnGeneral}) when using regularizer $\calR_5(\cdot)$.

Next we consider the low-rankness regularization via matricization:
$$
\calR_6(A)={1\over 3}\left(\|\calM_1(A)\|_\ast+\|\calM_2(A)\|_\ast+\|\calM_3(A)\|_\ast\right).
$$
It is not hard to see that
$$
\calR_6^\ast(A)=3\max\left\{\|\calM_1(A)\|_s, \|\calM_2(A)\|_s, \|\calM_3(A)\|_s\right\}.
$$

\blems
\label{LemLowRank}
There exists a constant $0 < c < \infty$ such that
\begin{equation}
\label{eq:glassog}
w_G[\mathbb{B}_{\calR_6}(1)] \leq c \sqrt{\max\{d_1d_2,d_2d_3,d_1d_3\}}.
\end{equation}
\elems

On the other hand,
\blems
\label{LemLowRank1}
For any $A\in \Theta_5(r)$,
$$s_6(\calA(P_1,P_2,P_3))=\sup_{B \in \mathcal{A}(P_1,P_2,P_3)/\{0\}} \frac{\mathcal{R}_6^2(B)}{\|B\|_{\rm F}^2}\le r.$$
\elems

Lemmas \ref{LemLowRank} and \ref{LemLowRank1} suggest that
$$
\sup_{T\in \Theta_5(r)}\max\left\{\|\widehat{T}_6-T\|_n^2, \|\widehat{T}_6-T\|_{\rm F}^2\right\} \lesssim {r\max\{d_1d_2, d_2d_3, d_1d_3\}\over n},
$$
with high probability by taking
$$
\lambda\asymp \sqrt{\max\{d_1d_2, d_2d_3, d_1d_3\}\over n}.
$$
where $\widehat{T}_6$ is the regularized least squares estimate defined by (\ref{EqnGeneral}) when using regularizer $\calR_6(\cdot)$.

Comparing with the rates for estimates with regularizers $\calR_5$ and $\calR_6$, we can see the benefit of using $\calR_5$. For any $T\in \Theta_5(r)$, If we apply regularizer $\calR_5$, then
$$
\|\widehat{T}_{5}-T\|_{\rm F}^2\lesssim {r^2(d_1+d_2+d_3)\over n}.
$$
This is to be compared with the risk bound for matricized regularization:
$$
\|\widehat{T}_{6}-T\|_{\rm F}^2\lesssim {r\max\{d_1d_2,d_2d_3,d_1d_3\}\over n}.
$$
Obviously $\widehat{T}_{5}$ always outperform $\widehat{T}_6$ since $r\le \min\{d_1,d_2,d_3\}$. The advantage of $\widehat{T}_5$ is typically rather significant since in general $r\ll \min\{d_1,d_2,d_3\}$. On the other hand, $\widehat{T}_6$ is more amenable for computation. 

Both upper bounds on Frobenius error on $\widehat{T}_{5}$ and $\widehat{T}_{6}$ are novel results and complement the existing results on tensor completion~\cite{GandRect11, MuGoldfarb14} and \cite{YuanZhang14}. Neither $\widehat{T}_{5}$ nor $\widehat{T}_6$ is minimax optimal and remains an interesting question as to whether there exists a convex regularization approach that is minimax optimal.

\section{Specific Statistical Problems}
\label{SecExamples}

In this section, we apply our results to several concrete examples where we are attempting to estimate a tensor under certain sparse or low rank constraints, and show that the regularized least squares estimate $\widehat{T}$ is typically minimiax rate optimal with appropriate choices of regularizers. In particular we focus on the multi-response aspect of the general framework to provide novel upper bounds and matching minimax lower bounds.

\subsection{Multi-Response regression with large $p$}

The first example we consider is the multi-response regression model:
\begin{equation*}
Y_{k}^{(i)} = \sum_{j=1}^p \sum_{\ell=1}^m {X^{(i)}_{j \ell} T_{j \ell k}} + \epsilon^{(i)}_{k},
\end{equation*}
where $1 \leq i \leq n$ represents the index for each sample, $1 \leq k \leq m$ represents the index for each response and $1 \leq j \leq p$ represents the index for each feature. For the multi-response regression problem we have $N = 3$, $M = 2$, $d_1 = d_2 = m$ which represents the total number of responses and $d_3 = p$, which represent the total number of parameters.

Since we are in the setting where $p$ is large but only a small number $s$ are relevant, we define the subspace:
\begin{equation*}
\mathcal{T}_1 = \left\{A \in \mathbb{R}^{m \times m \times p}\;|\; \sum_{j=1}^p \mathbb{I}({\|A_{\cdot\cdot j}\|_{\rm F}} \neq 0) \leq s \right\}.	
\end{equation*}
Furthermore for each $i$ we assume $X^{(i)} \in \mathbb{R}^{m \times p}$ where each entry of $X^{(i)}$, $[X^{(i)}]_{k,j}$, corresponds to the $j^{th}$ feature for the $k^{th}$ response. For simplicity, we assume the $X^{(i)}$'s are independent Gaussian with covariance $\widetilde{\Sigma} \in \mathbb{R}^{mp \times mp}$. The penalty function we are considering is:
\begin{equation}
\label{eq:multipen1}
\mathcal{R}(A) = \sum_{j = 1}^{p} \| A_{\cdot\cdot j} \|_{\rm F},
\end{equation}
and the corresponding dual function applied to the i.i.d. Gaussian tensor $G$  is:
\begin{equation*}
\mathcal{R}^*(G) = \max_{1 \leq j \leq p} \| G_{..j}\|_{\rm F}.
\end{equation*}

\btheos
\label{ThmUpperMultiReg}
Under the multi-response regression model with $T\in \calT_1$ and independent Gaussian design where $c_\ell^2 \leq \lambda_{min}(\widetilde{\Sigma}) \leq \lambda_{max}(\widetilde{\Sigma}) \leq c_u^2$, if
$$\lambda \geq 3 \sigma c_u \sqrt{\frac{\max\{m^2,\log p\}}{n}},$$
such that $\sqrt{s} \lambda$ converges to zero as $n$ increases, then there exist some constants $c_1,c_2>0$ such that with probability at least $1 - p^{-c_1}$
\begin{equation*}
\max\left\{\|\widehat{T} - T\|_n^2, \|\widehat{T} - T\|_{\rm F}^2\right\} \leq \frac{c_2c_u^2}{c_{\ell}^2} s \lambda^2,
\end{equation*}
when $n$ is sufficiently large, where $\widehat{T}$ is the regularized least squares estimate defined by (\ref{EqnGeneral}) with regularizer given by (\ref{eq:multipen1}). In addition,
\begin{equation*}
\min_{\widetilde{T}} \max_{T \in  \mathcal{T}_1} \|\widetilde{T} - T\|_{\rm F}^2 \geq \frac{c_3  \sigma^2 s \max\{m^2, \log p/s\}}{c_u^{2} n},
\end{equation*}
for some constant $c_3>0$, with probability at least $1/2$, where the minimum is taken over all estimators $\widetilde{T}$ based on data $\{(X^{(i)},Y^{(i)}): 1\le i\le n\}$.
\etheos

Theorem \ref{ThmUpperMultiReg} shows that when taking
$$
\lambda\asymp \sqrt{\frac{\max\{m^2,\log p\}}{n}},
$$
the regularized least squares estimate defined by (\ref{EqnGeneral}) with regularizer given by (\ref{eq:multipen1}) achieves minimax optimal rate of convergence over the parameter space $\calT_1$.

Alternatively, there are settings where the effect of covariates on the multiple tasks may be of low rank structure. In such a situation, we may consider
\begin{equation*}
\mathcal{T}_2 = \left\{A \in \mathbb{R}^{m \times m \times p}\;|\; \sum_{j=1}^p {\rm rank}({A_{..j}}) \leq r\right\}.	
\end{equation*}
An appropriate penalty function in this case is:
\begin{equation}
\label{eq:multipen2}
\mathcal{R}(A) = \sum_{j = 1}^{p} \| A_{..j} \|_{\ast},
\end{equation}
and the corresponding dual function applied to $G$ is:
\begin{equation*}
\mathcal{R}^*(G) = \max_{1 \leq j \leq p} \| G_{..j}\|_{s}.
\end{equation*}

\btheos
\label{ThmUpperMultiReg1}
Under the multi-response regression model with $T\in \calT_2$ and independent Gaussian design where $c_\ell^2 \leq \lambda_{min}(\widetilde{\Sigma}) \leq \lambda_{max}(\widetilde{\Sigma}) \leq c_u^2$, if 
$$\lambda \geq 3 \sigma c_u\sqrt{\frac{\max\{m,\log p\}}{n}},$$
such that $\sqrt{r} \lambda$ converges to zero as $n$ increases, then there exist some constants $c_1,c_2>0$ such that with probability at least $1 - p^{-c_1}$,
\begin{equation*}
\max\left\{\|\widehat{T} - T\|_n^2, \|\widehat{T} - T\|_{\rm F}^2\right\} \leq \frac{c_2 c_u^2}{c_{\ell}^2} r \lambda^2
\end{equation*}
when $n$ is sufficiently large, where $\widehat{T}$ is the regularized least squares estimate defined by (\ref{EqnGeneral}) with regularizer given by (\ref{eq:multipen2}). In addition,
\begin{equation*}
\min_{\widetilde{T}} \max_{T \in  \mathcal{T}_2} \|\widetilde{T} - T\|_{\rm F}^2 \geq \frac{c_3 \sigma^2 r \max\{m, \log (p/r)\}}{c_u^{2} n},
\end{equation*}
for some constant $c_3>0$, with probability at least $1/2$, where the minimum is taken over all estimators $\widetilde{T}$ based on data $\{(X^{(i)},Y^{(i)}): 1\le i\le n\}$.
\etheos

Again Theorem \ref{ThmUpperMultiReg1} shows that by taking
$$
\lambda\asymp \sqrt{\frac{\max\{m,\log p\}}{n}},
$$
the regularized least squares estimate defined by (\ref{EqnGeneral}) with regularizer given by (\ref{eq:multipen2}) achieves minimax optimal rate of convergence over the parameter space $\calT_2$. Comparing with optimal rates for estimating a tensor from $\calT_1$, one can see the benefit and importance to take advantage of the extra low rankness if the true coefficient tensor is indeed from $\calT_2$. As far as we are aware, these are the first results that provide upper bounds and matching minimax lower bounds for high-dimensional multi-response regression with sparse or low-rank slices. As pointed out earlier, the challenge in going from scalar to multiple response is proving Lemma~\ref{LemEmpLowerBound} which is an analog of restricted strong convexity.

\subsection{Multivariate sparse auto-regressive models}

Now we consider the setting of vector auto-regressive models. In this case, our generative model is:
\begin{equation}
\label{EqnVAR}
X^{(t+p)} = \sum_{j=1}^p {A_j X^{(t+p-j)}} + \epsilon^{(t)}, 
\end{equation}
where $1 \leq t \leq n$ represents the time index, $1 \leq j \leq p$ represents the lag index, $\{X^{(t)}\}_{t=0}^{n+p}$ is an $m$-dimensional vector, $\epsilon^{(t)} \sim \mathcal{N}(0, \sigma^2 I_{m \times m})$ represents the additive noise. Note that the parameter tensor $T$ is an $m \times m \times p$ tensor so that $T_{\cdot\cdot j}=A_j$, and $T_{k\ell j}$ represents the co-efficient of the $k^{th}$ variable on the $\ell^{th}$ variable at lag $j$. This model is studied by \cite{BasuMichail15} where $p$ is relatively small (to avoid introducing long-range dependence) and $m$ is large. Our main results allow more general structure and regularization schemes than those considered in \cite{BasuMichail15}.

Since we assume the number of series $m$ is large, and there are $m^2$ possible interactions between the series we assume there are only $s \ll m^2$ interactions in total.
\begin{equation}
\label{EqnVARClass}
\mathcal{T}_3 = \left\{A \in \mathbb{R}^{m \times m \times p}\;|\; \sum_{k=1}^m \sum_{\ell=1}^m \mathbb{I}({A_{k\ell\cdot}} \neq \zero) \leq s \right\}.	
\end{equation}
The penalty function we are considering is:
\begin{equation}
\label{eq:varpen}
\mathcal{R}(A) = \sum_{k = 1}^{m} \sum_{\ell=1}^m \| A_{k\ell\cdot}\|_{\ell_2},
\end{equation}
and the corresponding dual function applied to $G$ is:
\begin{equation*}
\mathcal{R}^*(G) = \max_{1 \leq k,\ell \leq m} \| G_{k,\ell,.}\|_{\ell_2}.
\end{equation*}
The challenge in this setting is that the $X$'s are highly dependent and we use the results developed in \cite{BasuMichail15} to prove that (\ref{AssCov}) is satisfied. 

Prior to presenting the main results, we introduce concepts developed in \cite{BasuMichail15} that play a role in determining the constants $c_{u}^2$ and $c_{\ell}^2$ which relate to the stability of the auto-regressive processes. A $p$-variate Gaussian time series is defined by its auto-covariance matrix function
$$\Gamma_X(h) = \mbox{Cov}(X^{(t)}, X^{(t+h)}),$$
for all $t, h \in \mathbb{Z}$. Further, we define the spectral density function:
\begin{equation*}
f_X(\theta) := \frac{1}{2 \pi} \sum_{\ell = -\infty}^{\infty} {\Gamma_X(\ell) e^{-i \ell \theta}},\;\; \theta \in [-\pi, \pi].
\end{equation*}
To ensure the spectral density is bounded, we make the following assumption:
\begin{equation*}
\mathcal{M}(f_X) := \esssup_{\theta} \Lambda_{\max}(f_X(\theta)) < \infty.
\end{equation*}
Further, we define the matrix polynomial
$$\mathcal{A}(z) = I_{m \times m} - \sum_{j=1}^p {A_j z^j}$$
where $\{A_j\}_{j=1}^p$ denote the back-shift matrices, and $z$ represents any point on the complex plane. Note that for a stable, invertible AR($p$) process,
\begin{equation*}
f_X(\theta) = \frac{1}{2 \pi} \mathcal{A}^{-1}(e^{-i \theta}) \overline{\mathcal{A}^{-1}(e^{-i \theta})}. 
\end{equation*}
We also define the lower extremum of the spectral density:
\begin{equation*}
m(f_X) := \essinf_{\theta} \Lambda_{\min}(f_X(\theta)).
\end{equation*}
Note that $m(f_X)$ and $\mathcal{M}(f_X)$ satisfy the following bounds:
\begin{equation*}
m(f_X) \geq \frac{1}{2 \pi \mu_{\max}(\mathcal{A})},\qquad {\rm and}\qquad \mathcal{M}(f_X) \leq \frac{1}{2 \pi \mu_{\min}(\mathcal{A})},
\end{equation*}
where
$$\mu_{\min}(\mathcal{A}) :=  \min_{|z| = 1} {\Lambda_{\min}(\overline{\mathcal{A}(z)}\mathcal{A}(z)) }$$
and
$$\mu_{\max}(\mathcal{A}) :=  \max_{|z| = 1} {\Lambda_{\max}(\overline{\mathcal{A}(z)}\mathcal{A}(z)) }.$$

From a straightforward calculation, we have that for any fixed $\Delta$:
\begin{equation}
\label{EqnCondVAR}
\frac{1}{\mu_{\max}} \|\Delta\|_{\rm F}^2 \leq \mathbb{E}\left[\|\Delta\|_n^2\right] \leq \frac{1}{\mu_{\min}} \|\Delta\|_{\ell_2}^2.
\end{equation}
Hence $c_u^2 = 1/\mu_{\min}$ and $c_{\ell}^2 = 1/\mu_{\max}$. Now we state our main result for auto-regressive models.

\btheos
\label{ThmUpperVAR}
Under the vector auto-regressive model defined by~\eqref{EqnVAR} with $T\in \calT_3$, if
$$\lambda \geq 3 \sigma \sqrt{\frac{\max\{p,2\log m\}}{n\mu_{\min}}},$$
such that $\sqrt{s}\lambda$ converges to zero as $n$ increases, then there exist some constants $c_1,c_2>0$ such that with probability at least $1-m^{-c_1}$,
\begin{equation*}
\max\left\{\|\widehat{T} - T\|_n^2, \|\widehat{T} - T\|_{\rm F}^2\right\} \leq \frac{c_2 \mu_{\max}}{\mu_{\min}}s \lambda^2,
\end{equation*}
when $n$ is sufficiently large, where $\widehat{T}$ is the regularized least squares estimators defined by (\ref{EqnGeneral}) with regularizer given by (\ref{eq:varpen}). In addition,
\begin{equation*}
\min_{\widetilde{T}} \max_{T \in  \mathcal{T}_3} \|\widetilde{T} - T\|_{\rm F}^2 \geq  c_3\mu_{\min} \sigma^2 \frac{s \max\{p, \log (m/\sqrt{s})\}}{n},
\end{equation*}
for some constant $c_3>0$, with probability at least $1/2$, where the minimum is taken over all estimators $\widetilde{T}$ based on data $\{X^{(t)}: t=0,\ldots,n+p\}$.
\etheos

Theorem \ref{ThmUpperVAR} provides, to our best knowledge, the only lower bound result for multivariate time series. The upper bound is also novel and is
different from Proposition 4.1 in \cite{BasuMichail15} since we impose sparsity only on the large $m$ directions and not over the $p$ lags, whereas \cite{BasuMichail15} impose sparsity through vectorization. Note that Proposition 4.1 in \cite{BasuMichail15} follows directly from Lemma~\ref{LemSparsity} with $d_1 = p$ and $d_2 = d_3 = m$. Using the sparsity regularizer \cite{BasuMichail15} vectorize the problem and prove restricted strong convexity whereas since we leave the problem as a multi-response problem, we requried the more refined technique used for proving Lemma~\ref{LemEmpLowerBound}. 

\subsection{Pairwise interaction tensor models}

\label{SecTheoryPairwise}

Finally, we consider the tensor regression (\ref{eq:model}) where $T$ follows a pairwise interaction model. More specifically, $(X^{(i)},Y^{(i)})$, $i=1,2,\ldots, n$ are independent copies of a random couple $X\in \RR^{d_1\times d_2\times d_3}$ and $Y\in \RR$ such that
$$
Y=\langle X, T\rangle +\epsilon
$$
and
$$
T_{j_1j_2j_3}=A^{(12)}_{j_1j_2}+A^{(13)}_{j_1j_3}+A^{(23)}_{j_2j_3}.
$$
Here $A^{(k_1,k_2)}\in \RR^{d_{k_1}\times d_{k_2}}$ such that
$$
A^{(k_1,k_2)}\one =\zero, \qquad {\rm and}\qquad  (A^{(k_1,k_2)})^\top\one =\zero.
$$
The pairwise interaction was used originally by \cite{RendleSchmidt09,RendleSchmidt10} for personalized tag recommendation, and later analyzed in \cite{ChenLyu13}. \cite{Hoff2014} briefly introduced a single index additive model (amongst other tensor models) which is a sub-class of the pairwise interaction model. The regularizer we consider is:
\begin{equation}
\label{eq:pairwise0}
\calR(A)= \|A^{(12)}\|_\ast+\|A^{(13)}\|_\ast+\|A^{(23)}\|_\ast.
\end{equation}
It is not hard to see that $\calR$ defined above is decomposable with respect to $\calA(P_1,P_2,P_3)$ for any projection matrices.

Let
\begin{eqnarray*}
\calT_4=\{A\in \RR^{d_1\times d_2\times d_3}: A_{j_1j_2j_3}=A^{(12)}_{j_1j_2}+A^{(13)}_{j_1j_3}+A^{(23)}_{j_2j_3}, A^{(k_1,k_2)}\in \RR^{d_{k_1}\times d_{k_2}},\\
 A^{(k_1,k_2)}\one =\zero, \qquad {\rm and}\qquad  (A^{(k_1,k_2)})^\top\one =\zero\\
 \max_{k_1,k_2} {\rm rank}(A^{(k_1,k_2)})\le r\}.
\end{eqnarray*}
For simplicity, we assume i.i.d. Gaussian design so $c_\ell^2 = c_u^2 = 1$.
\btheos
\label{th:pairwise}
Under the pairwise interaction model with $T\in \calT_4$, if
$$
\lambda \geq 3\sigma \sqrt{\max\{d_1,d_2,d_3\}\over n},
$$
such that $\sqrt{r}\lambda$ converges to zero as $n$ increases, then there exist constants $c_1,c_2>0$ such that with probability at least $1-\min\{d_1,d_2,d_3\}^{-c_1}$,
\begin{equation*}
\max\left\{\|\widehat{T} - T\|_n^2, \|\widehat{T} - T\|_{\rm F}^2\right\}  \leq c_2 r \lambda^2,
\end{equation*}
when $n$ is sufficiently large, where $\widehat{T}$ is the regularized least squares estimate defined by (\ref{EqnGeneral}) with regularizer given by (\ref{eq:pairwise0}). In addition,
\begin{equation*}
\min_{\widetilde{T}} \max_{T \in  \mathcal{T}_4} \|\widetilde{T} - T\|_{\rm F}^2 \geq  \frac{c_3 \sigma^2 r\max\{d_1,d_2,d_3\}}{n},
\end{equation*}
for some constant $c_3>0$, with probability at least $1/2$, where the minimum is taken over all estimate $\widetilde{T}$ based on data $\{(X^{(i)},Y^{(i)}): 1\le i\le n\}$.
\etheos

As in the other settings, Theorem \ref{th:pairwise} establishes the minimax optimality of the regularized least squares estimate (\ref{EqnGeneral}) when using an appropriate convex decomposable regularizer. Since this is single response and the norm involves matricization, this result is a straightforward extension to earlier results. 

\section{Numerical Experiments}

\label{SecSim}

In this section, we provide a series of numerical experiments that both support our theoretical results and display the flexibility of our general framework. In particular, we consider several different models including: third-order tensor regression with a scalar response (Section \ref{SecSimThird}); fourth-order tensor regression (Section \ref{SecSimFourth}); matrix-response regression with both group sparsity and low-rankness regularizers (Section \ref{SecSimMat});  multi-variate sparse auto-regressive models (Section \ref{SecSimVAR}); and pairwise interaction models (Section \ref{SecSimPairwise}). To perform the simulations in a computationally tractable way, we adapt the block coordinate descent approaches in multi-response case developed by \cite{SimonFriedmanHastie13}, and those developed by \cite{QinGoldfarb12} for univariate response settings, to capture group sparsity and low-rankness regularizers.

To fix ideas, in all numerical experiments, the covariate tensors $X^{(i)}$s were independent standard Gaussian ensembles (except for the multivariate auto-regressive models); and the noise $\epsilon^{(i)}$s are i.i.d. random tensors with elements following $N(0,\sigma^2)$ independently. As to the choice of tuning parameter, we adopt grid search on $\lambda$ to find the one with the least estimation error (in terms of mean squared error) in all our numerical examples.

\subsection{Third-order tensor regression}
\label{SecSimThird}
First we consider a third-order tensor regression model:
$$
Y^{(i)} = \langle B, X^{(i)}  \rangle + \epsilon^{(i)}
$$
where
$B\in \mathbb{R}^{d\times d \times d}, \quad Y^{(i)}, \epsilon^{(i)} \in \mathbb{R}, \quad X^{(i)}\in \mathbb{R}^{d\times d\times d}$. The regression coefficient tensor $B$ was generated as follows: the first $s$ slices $B_{\cdot\cdot 1},\ldots B_{\cdot\cdot s}$ are i.i.d standard normal ensembles; and the remaining slices $B_{\cdot\cdot s+1},\ldots B_{\cdot\cdot d_3}$ are set to be zero. Naturally, we consider here the group-sparsity regularizer:
$$
\min_{A\in \mathbb{R}^{d\times d\times d}} \left\{ \frac{1}{2n} \sum_{i=1}^n \|Y^{(i)}- \langle A, X^{(i)}  \rangle\|_{\rm F}^2 +\lambda \sum_{j_3 = 1}^{d} \|A_{\cdot \cdot j_3}\|_{\text{F}}  \right\}.
$$
\begin{figure}[!htb]\centering
   \begin{minipage}{0.32\textwidth}
     \includegraphics[width=\linewidth]{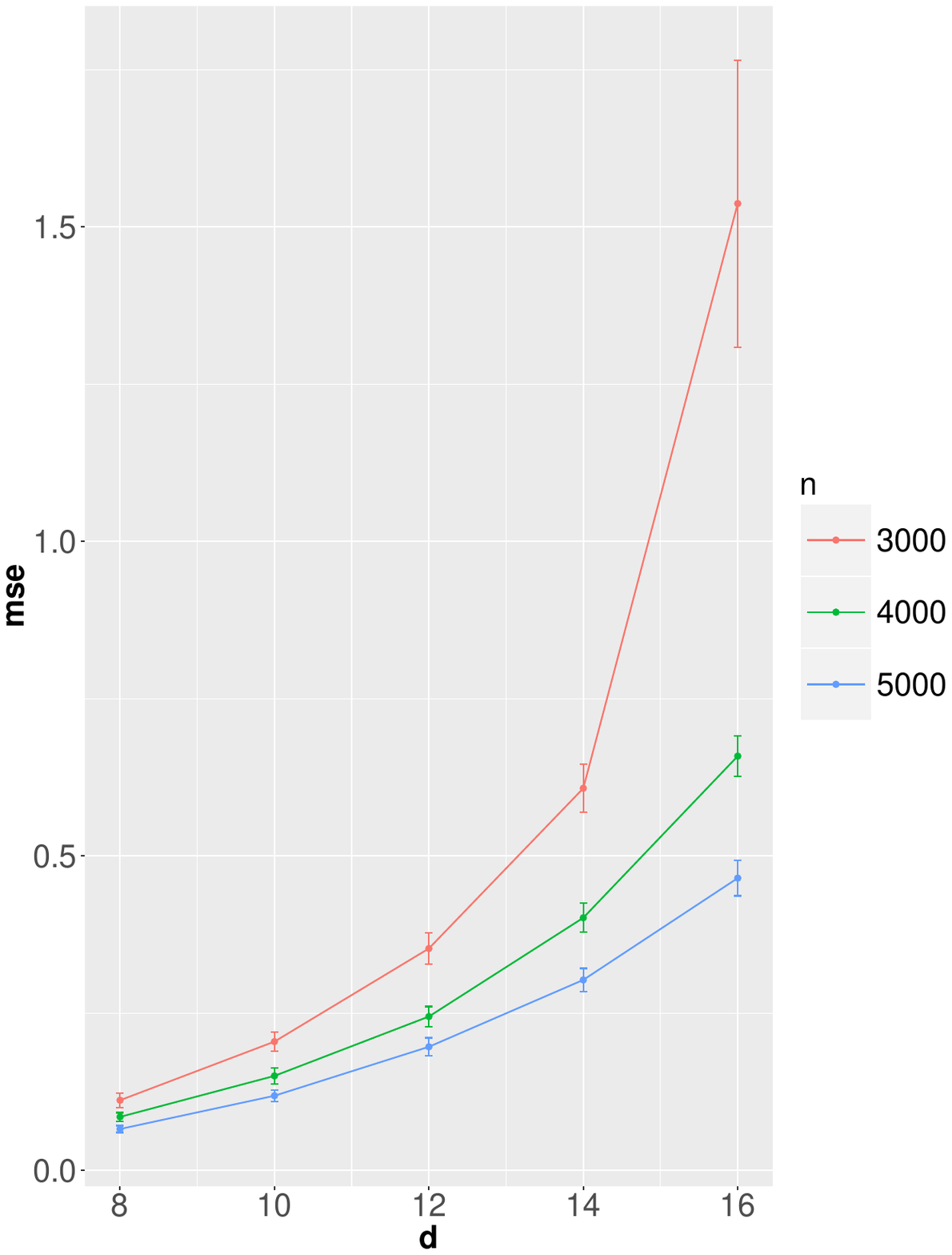}
   \end{minipage}
   \begin {minipage}{0.32\textwidth}
     \includegraphics[width=\linewidth]{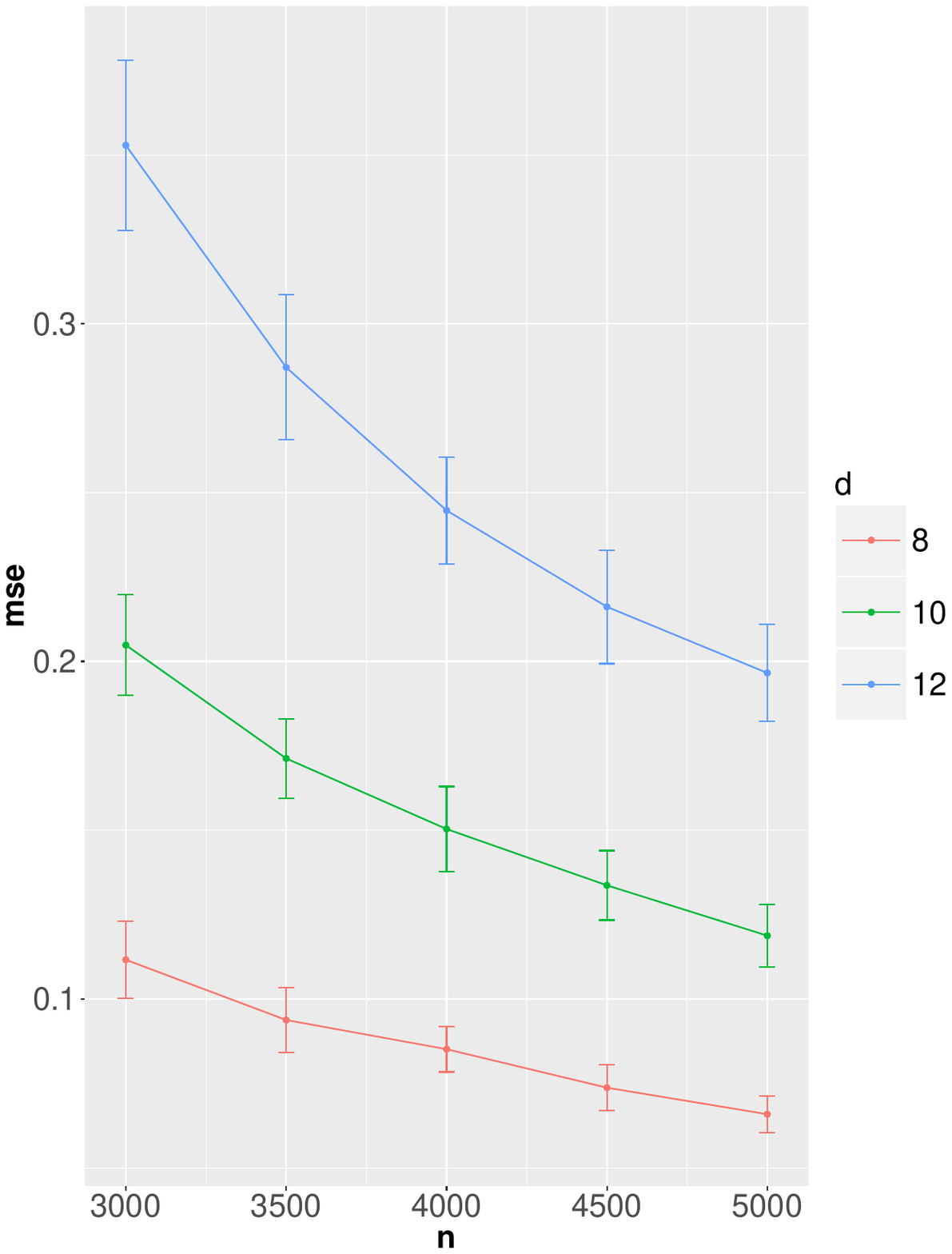}
   \end{minipage}
    \begin {minipage}{0.32\textwidth}
     \includegraphics[width=\linewidth]{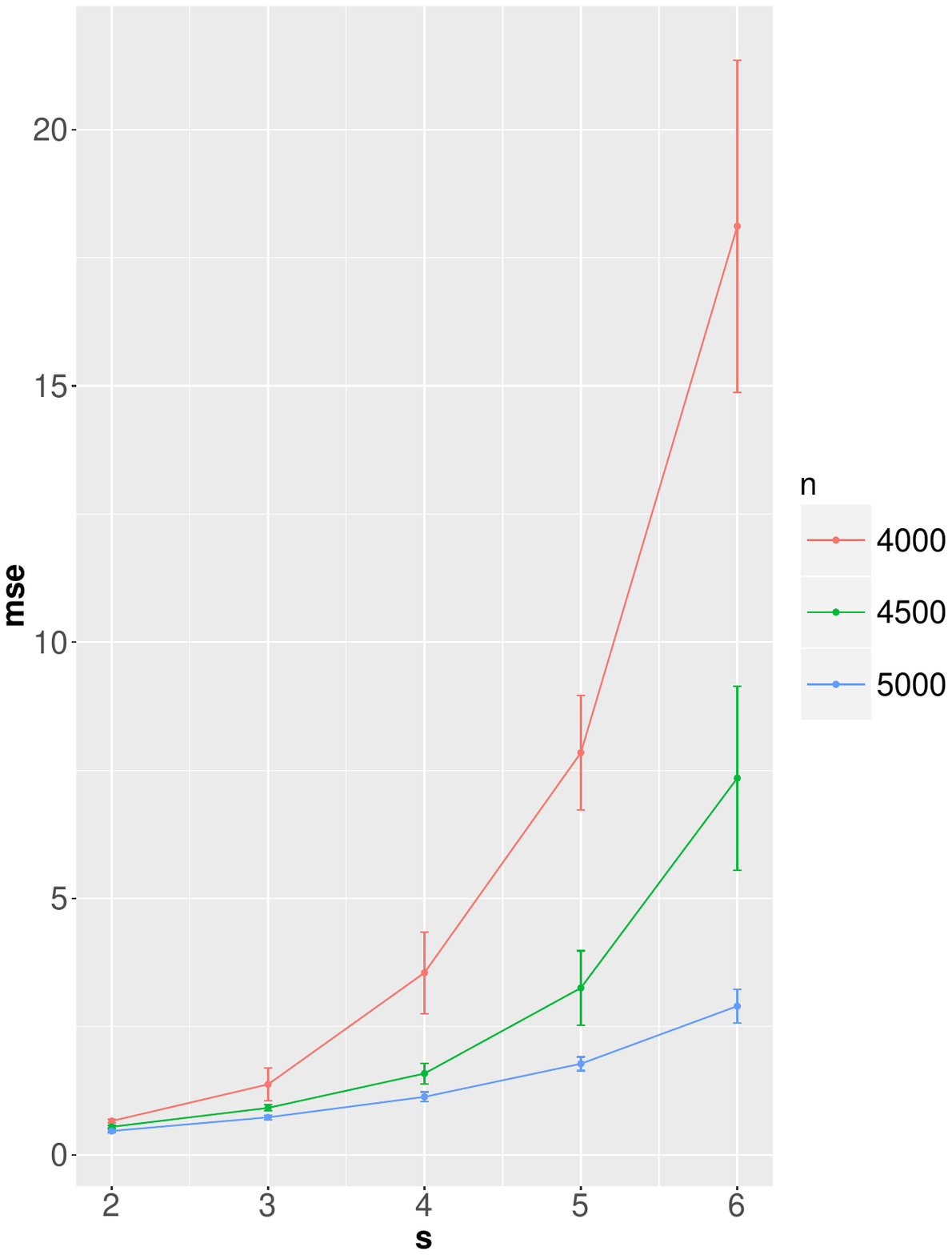}
   \end{minipage}
   \caption{Mean squared error of the group-sparsity regularization for third order tensor regression. The plot was based on 50 simulation runs and the error bars in each panel represent $\pm$ one standard deviation.}
   \label{figure3}
\end{figure}
Figure \ref{figure3} shows the mean squared error of the estimate averaged over 50 runs (with standard deviation) versus $d$, $n$ and $s$ respectively. In the left and middle panels, we set $s = 2$, whereas in the right panel, we fixed $d = 16$. As we can observe, the mean squared error increases approximately according to $d^2$, $s$, and $1/n$ which agrees with the risk bound given in Lemma~\ref{LemSparsityslice}.

We also considered a setting where $B$ is slice-wise low-rank. More specifically, the $s$ nonzero slices $B_{\cdot\cdot 1},\ldots B_{\cdot\cdot s}$ were random rank-$r$ matrices. In this case, the slice-wise low-rankness regularizer can be employed:
$$
\min_{A\in \mathbb{R}^{d_1\times d_2\times d_3}} \left\{ \frac{1}{2n} \sum_{i=1}^n \|Y^{(i)}- \langle A, X^{(i)}  \rangle\|_{\rm F}^2 +\lambda \sum_{j_3 = 1}^{d_3} \|A_{\cdot \cdot j_3}\|_*  \right\}.
$$
\begin{figure}[!htb]\centering
   \begin{minipage}{0.32\textwidth}
     \includegraphics[width=\linewidth]{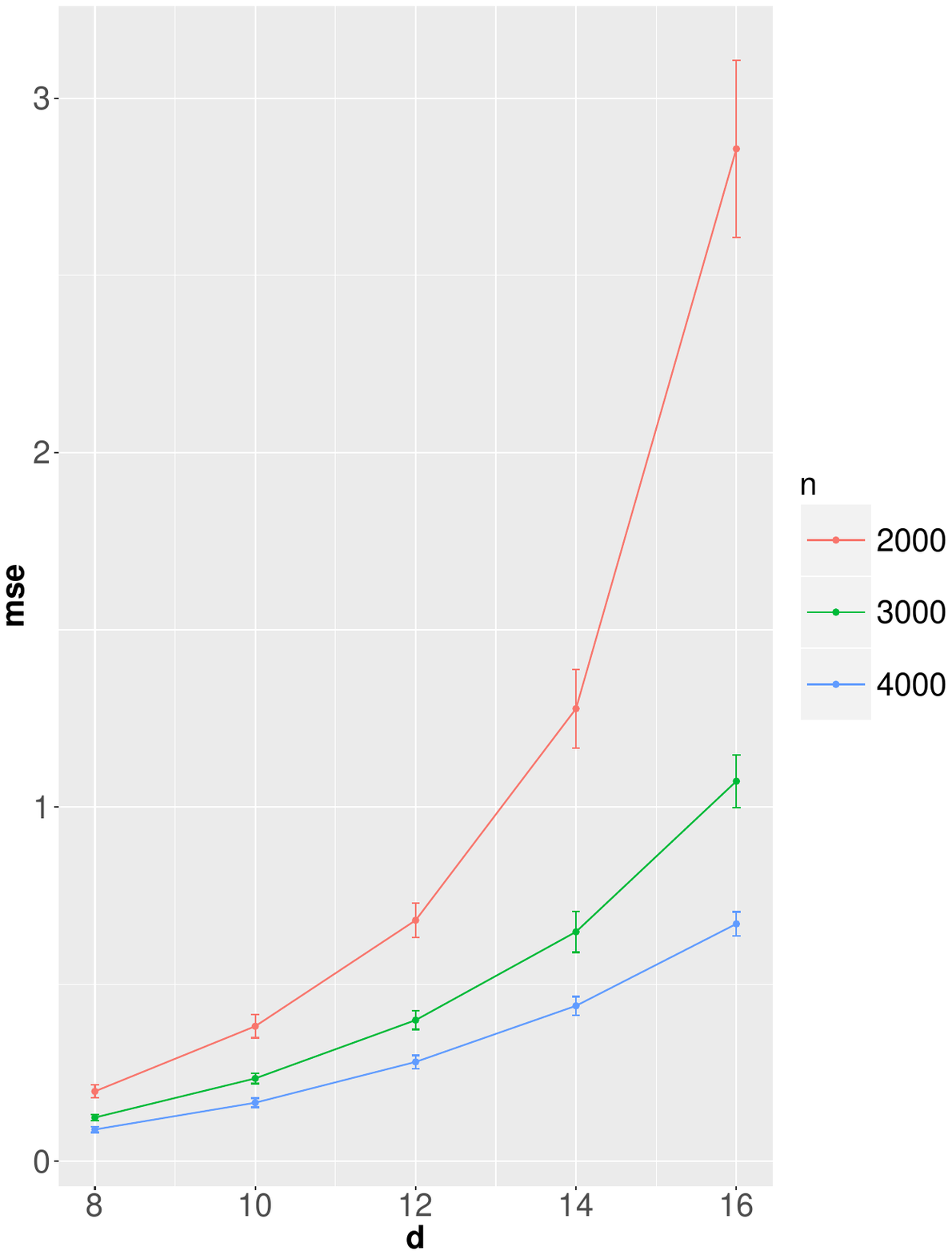}
   \end{minipage}
   \begin {minipage}{0.32\textwidth}
     \includegraphics[width=\linewidth]{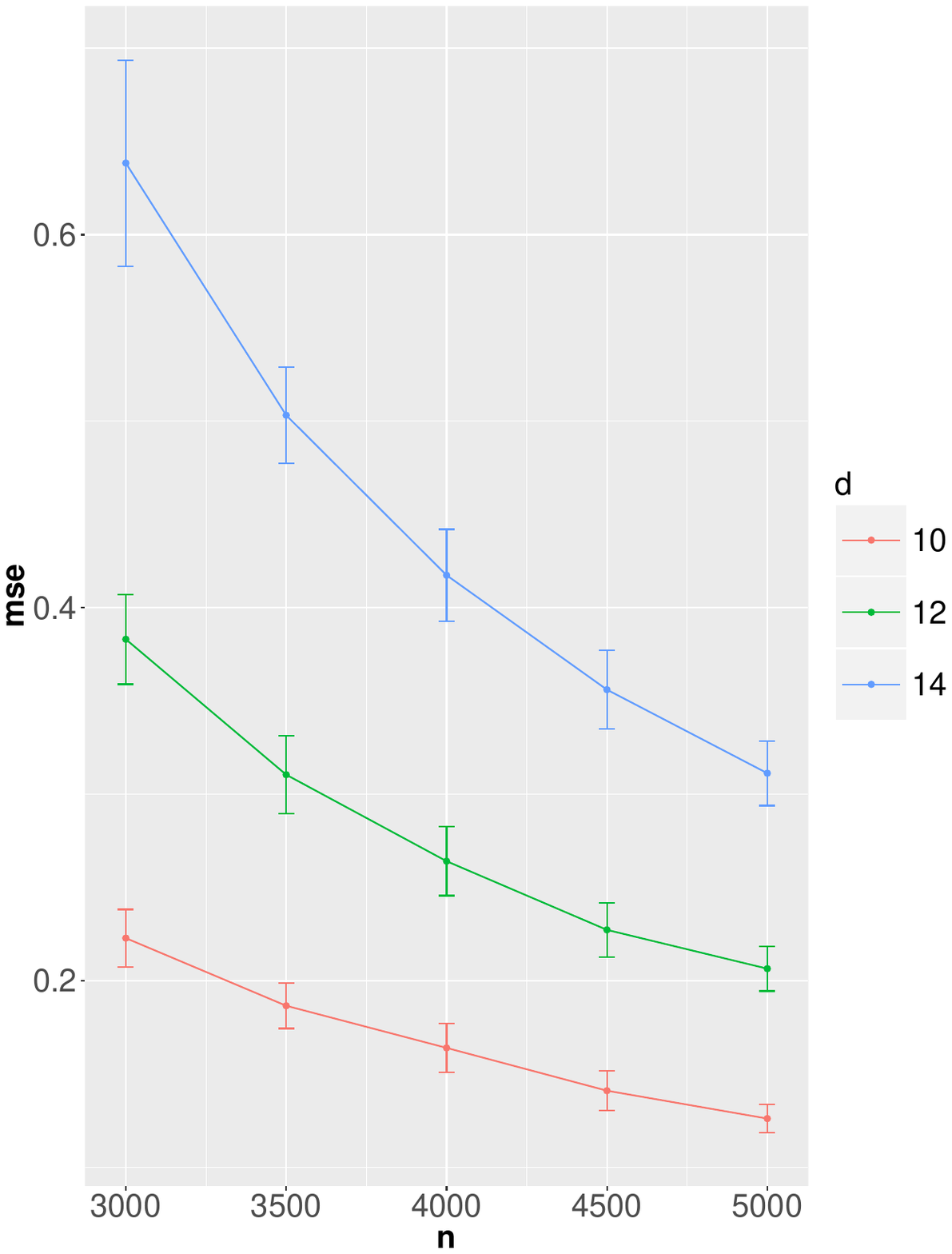}
   \end{minipage}
    \begin {minipage}{0.32\textwidth}
     \includegraphics[width=\linewidth]{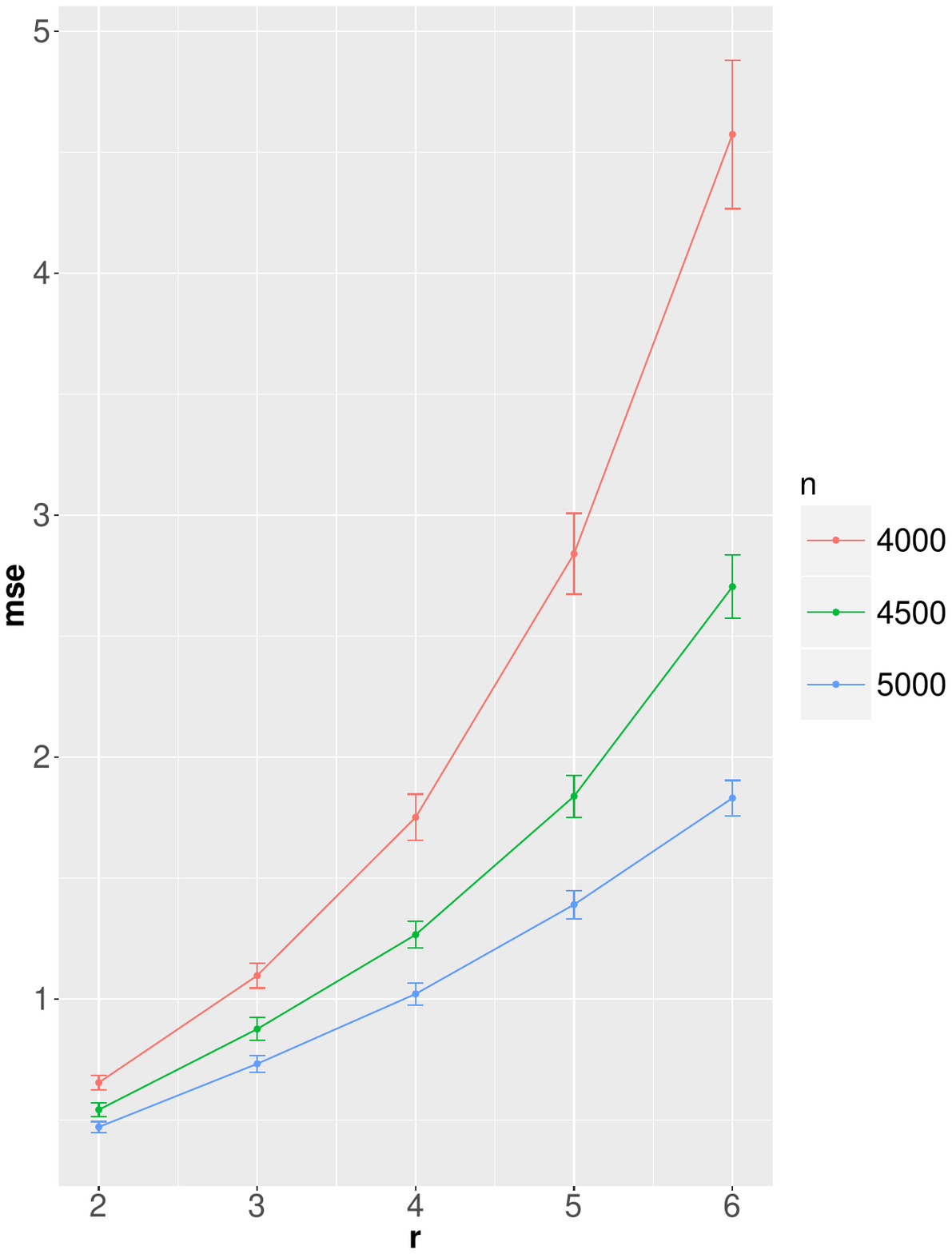}
   \end{minipage}
   \caption{Mean squared error for third order tensor regression with low-rank slices tensor coefficients. The plot was based on 50 simulation runs and the error bars in each panel represent $\pm$ one standard deviation.}\label{figure4}
\end{figure}
The performance of the estimate, averaged over 50 simulation runs, is summarized by Figure \ref{figure4} where in the left and middle panels $r = 2$, and in the right panel, $d = 16$. Once again, our results are consistent with our theoretical results.


\subsection{Fourth-order tensor regression}
\label{SecSimFourth}
Although we have focused on third order tensors for brevity, our treatment applies to higher order tensors as well. As an illustration, we now consider fourth order models where $ B\in \mathbb{R}^{d\times d \times d\times d}, \quad Y^{(i)}, \epsilon^{(i)} \in \mathbb{R}, \quad X^{(i)}\in \mathbb{R}^{d\times d\times d\times d}.$

To generate low-rank fourth-order tensors, we impose low CP rank as follows: generate four independent groups of $r$ independent random vectors of unit length, $\{u_{k,1}\}_{k=1}^r$, $\{u_{k,2}\}_{k=1}^r$, $\{u_{k,3}\}_{k=1}^r$ and $\{u_{k,4}\}_{k=1}^r$ via performing an SVD of Gaussian random matrix two times and keeping the $r$ pairs of leading singular vectors, and then compute the outer-product yielding a rank-$r$ tensor $B =\sum_{k=1}^r u_{k,1} \otimes u_{k,2} \otimes u_{k,3} \otimes u_{k,4}$.

We consider two different regularization schemes. First we impose low-rank structure through mode-1 matricization:
$$
\min_{A\in \mathbb{R}^{d\times d\times d\times d}} \left\{ \frac{1}{2n} \sum_{i=1}^n \|Y^{(i)}- \langle A, X^{(i)}  \rangle\|_{\rm F}^2 +\lambda  \|\mathcal{M}_1(A)\|_*  \right\}.
$$
Secondly we use the square matricization as follows:
$$
\min_{A\in \mathbb{R}^{d\times d\times d\times d}} \left\{ \frac{1}{2n} \sum_{i=1}^n \|Y^{(i)}- \langle A, X^{(i)}  \rangle\|_{\rm F}^2 +\lambda  \| \mathcal{M}_{12}(A)\|_*  \right\},
$$
where $\mathcal{M}_{12}(\cdot)$ reshape a fourth order tensor into a $d^2\times d^2$ matrix by collapsing its first two indices, and last two indices respectively. Table~\ref{table1} shows the average root-mean-square error (RMSE, for short) for both approaches. As we can see, the 2-by-2 approach appears superior to the 1-by-3 approach which is also predicted by the theory. 

\begin{table}[htbp]
\centering
\begin{tabular}{|c|c|c|c|c|c|c|c|c|}
	\hline
	$n$  & $d$ & $r$ & $\sigma$ & SNR & RMSE ({\small Mode-1 Matricization}) & RMSE ({\small Square Matricization})  \\
	\hline
	2000 & 7 & 5 & 10 & 3.0 (0.1) & 0.53 (0.01) & 0.51 (0.01)    \\
	\hline
	2000 & 7 & 3 & 10 & 1.5 (0.1) & 0.58 (0.02) & 0.49 (0.02)   \\
	\hline
	4000 & 10 & 3 & 10 & 1.7 (0.1) & 0.67 (0.01) & 0.51 (0.02)   \\
	\hline
\end{tabular}
\caption{Tensor regression with fourth order tensor covariates and scale response based on matricization: RMSE were computed based on 50 simulations runs. Numbers in parentheses are standard errors.}
\label{table1}
\end{table}

\subsection{Matrix-response regression}

\label{SecSimMat}
Our general framework can handle matrix-responses in a seamless fashion. For demonstration, we consider here matrix-response regression with both group sparsity and low-rankness regularizer. More specifically, the following model was considered:
$$
Y^{(i)} = \langle B, X^{(i)}\rangle + \epsilon^{(i)}
$$
where $B\in \mathbb{R}^{d\times d \times d}, \quad Y^{(i)}, \epsilon^{(i)} \in \mathbb{R}^{d\times d}, \quad X^{(i)}\in \mathbb{R}^{d}.$ As before, to impose group sparsity, the first $s$ slices of $B$ were generated as Gaussian ensembles and the remaining slices were set to zero.

For both the group sparsity and low-rankness regularizers, we used the matrix-version algorithm for group-penalized multi-response regression in \cite{SimonFriedmanHastie13}.  For each block of the coordinate descent, the sub-problem with both $\ell_1$ and nuclear norm penalty have closed-form solutions.

\begin{figure}[!htb]\centering
   \begin{minipage}{0.32\textwidth}
     \includegraphics[width=\linewidth]{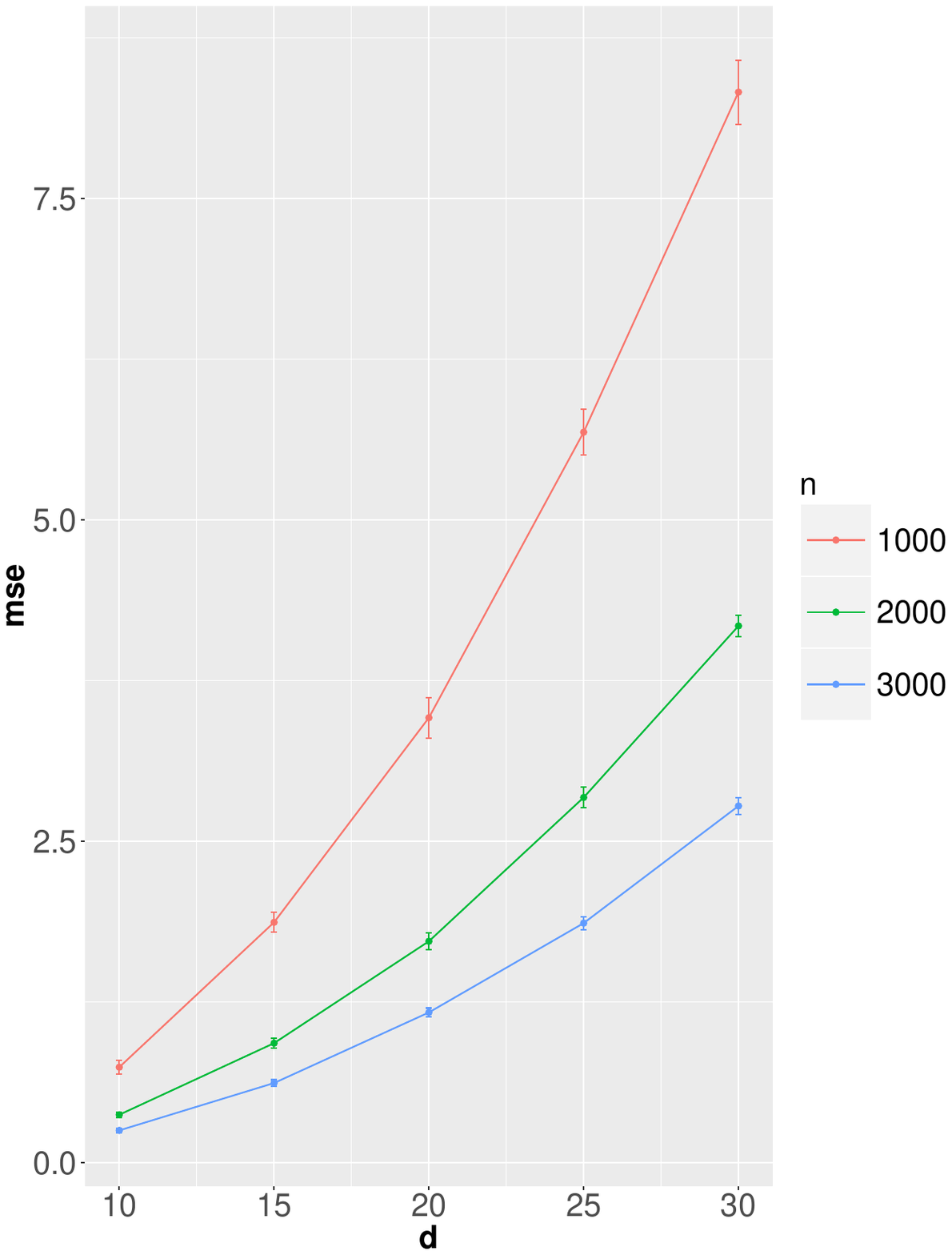}
   \end{minipage}
   \begin {minipage}{0.32\textwidth}
     \includegraphics[width=\linewidth]{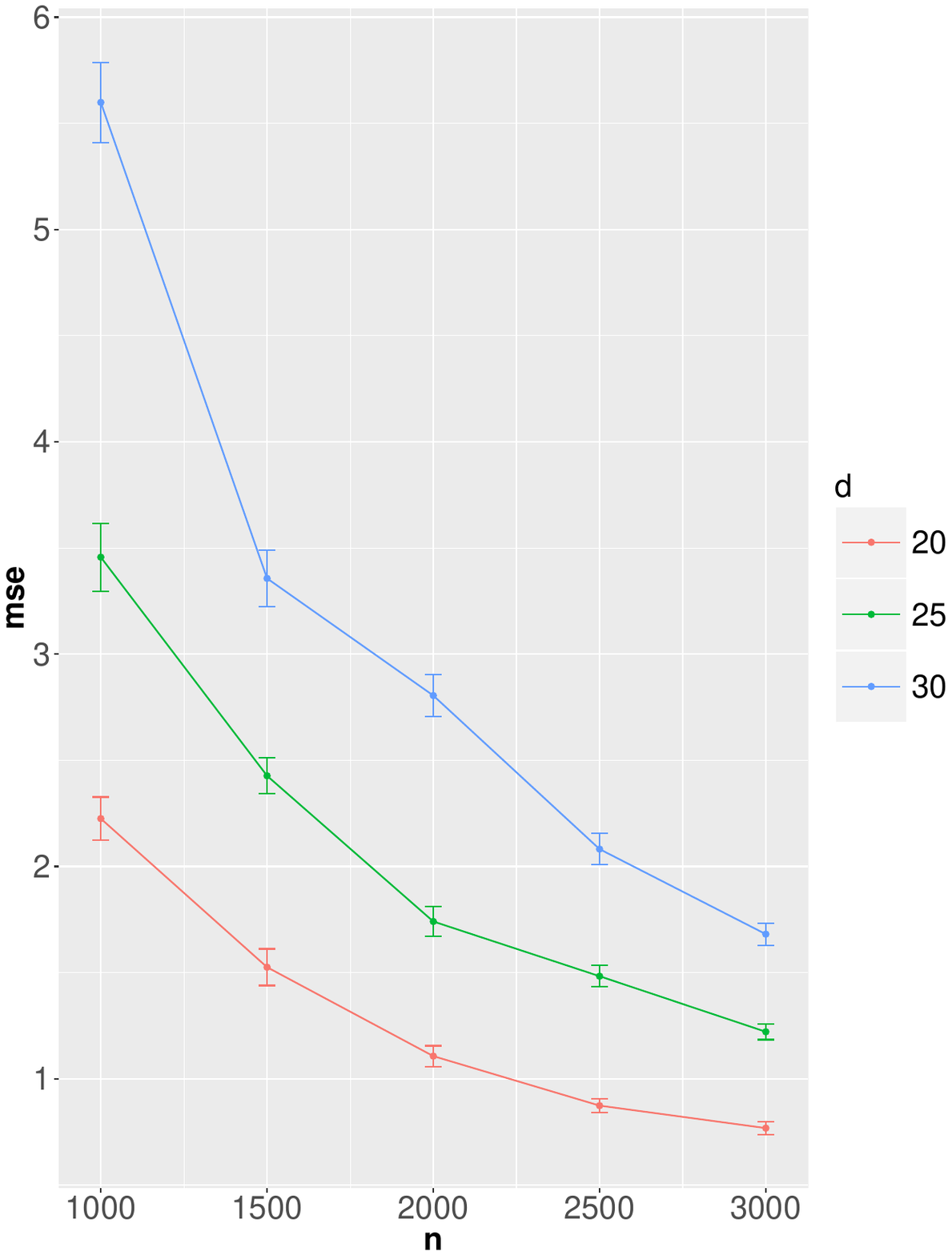}
   \end{minipage}
    \begin {minipage}{0.32\textwidth}
     \includegraphics[width=\linewidth]{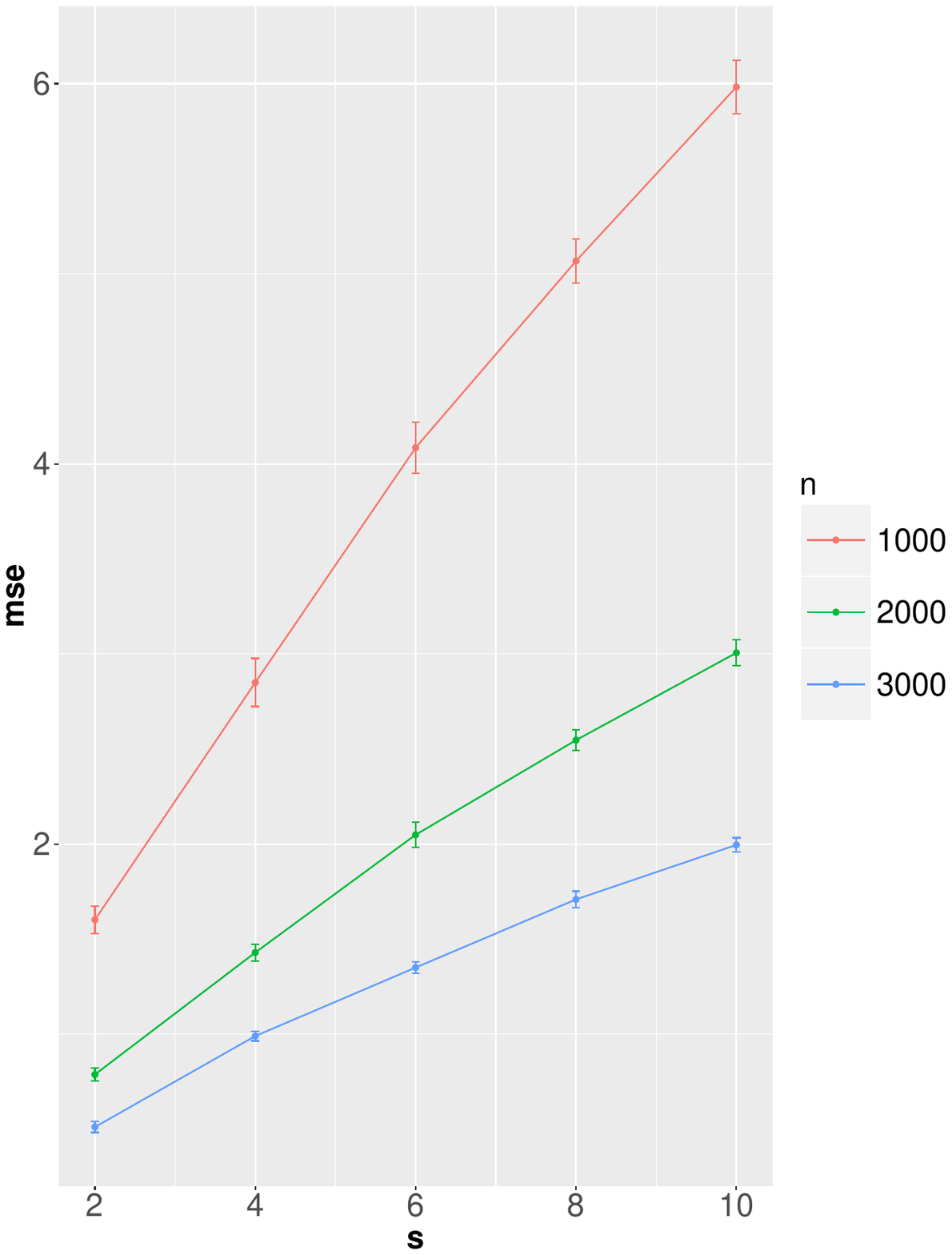}
   \end{minipage}
   \caption{Matrix response regression with sparse slices tensor coefficients. The plot was based on 50 simulation runs and the error bars in each panel represent $\pm$ one standard deviation.}\label{figure1}
\end{figure}

Figure~\ref{figure1} shows the average (with standard deviation) mean squared error over 50 runs versus the $d$, $n$ and $s$ parameter. (Here $d_1 = d_2 = d_3 = d$). As we observe, the mean-squared error increase approximately according to $\log d$, $s$, and $1/n$ which supports our upper bound in Theorem~\ref{ThmUpperMultiReg}.

\begin{figure}[!htb]\centering
   \begin{minipage}{0.32\textwidth}
     \includegraphics[width=\linewidth]{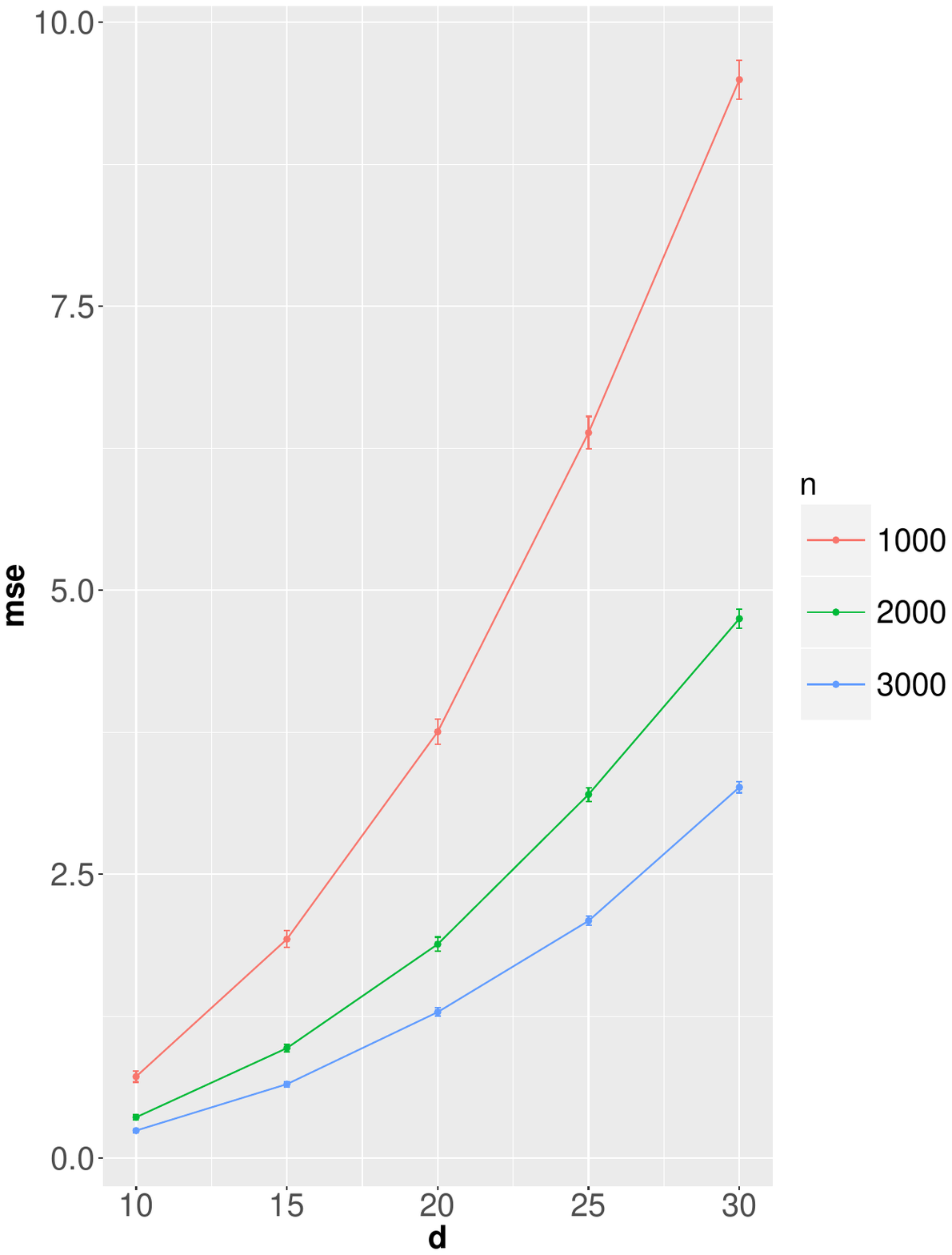}
   \end{minipage}
   \begin {minipage}{0.32\textwidth}
     \includegraphics[width=\linewidth]{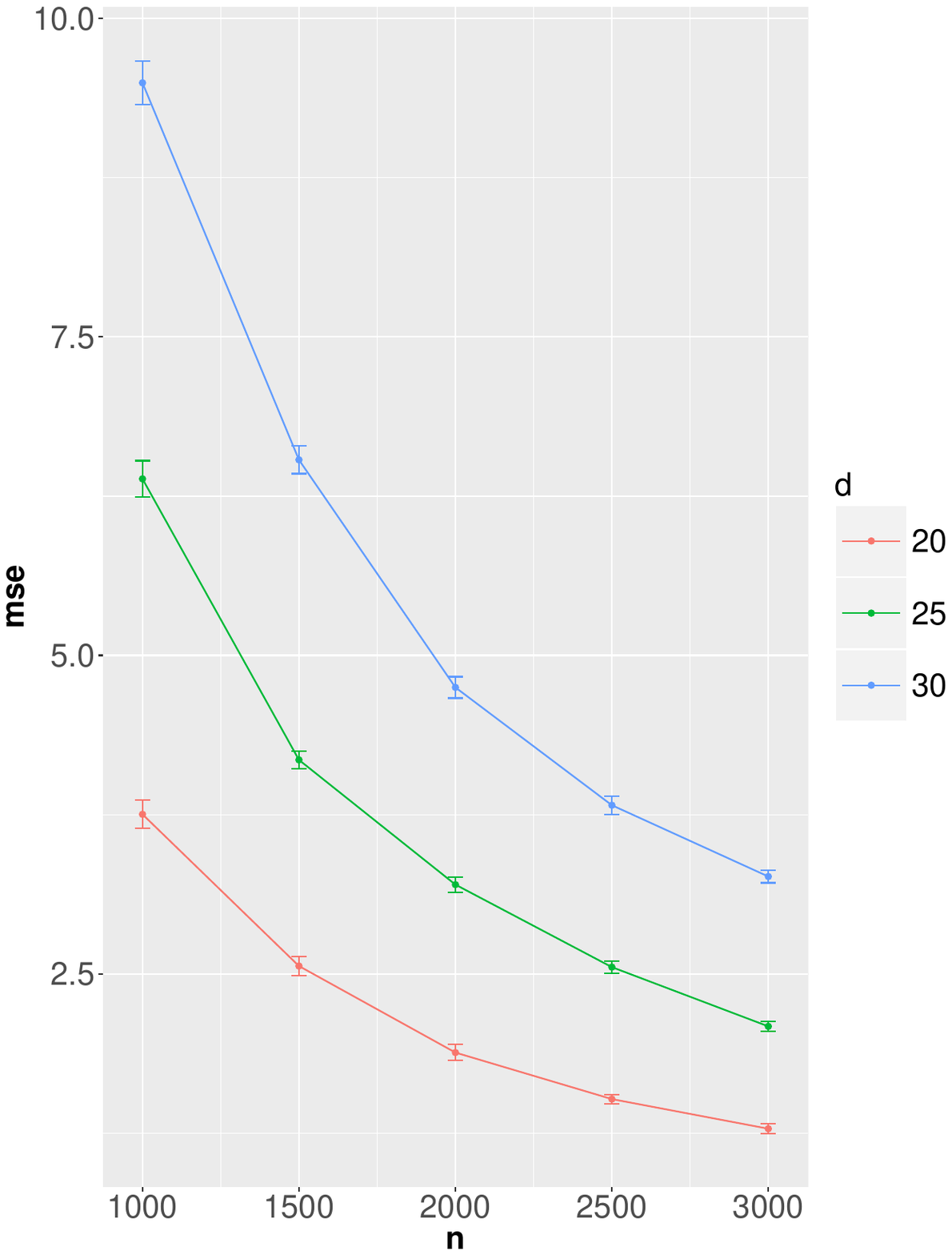}
   \end{minipage}
    \begin {minipage}{0.32\textwidth}
     \includegraphics[width=\linewidth]{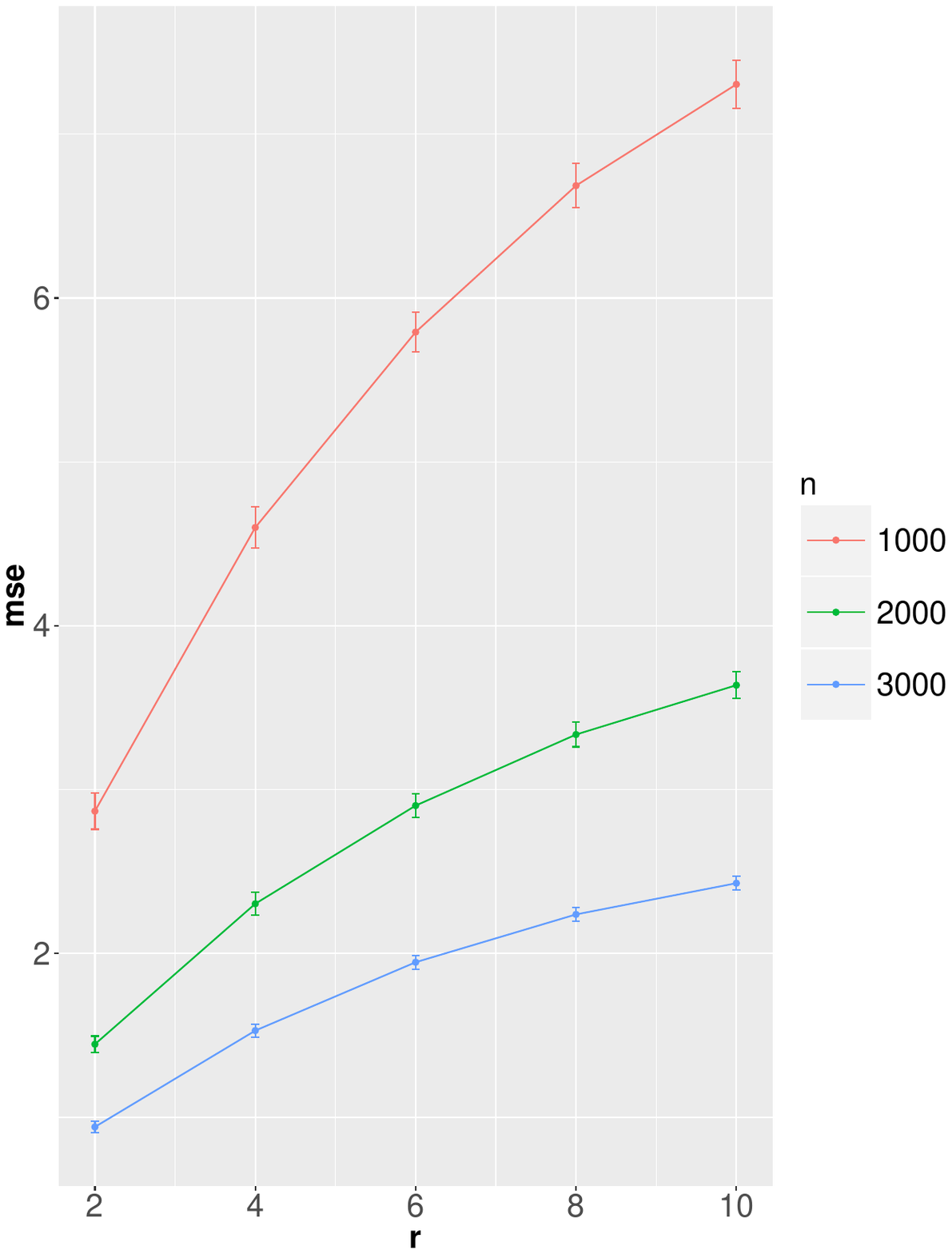}
   \end{minipage}
   \caption{Matrix response regression with low-rank slices tensor coefficients. The plot was based on 50 simulation runs and the error bars in each panel represent $\pm$ one standard deviation.}\label{figure2}
\end{figure}

We also generated low-rank $B$ in the same fashion as before. Figure \ref{figure2} plots the average (with standard deviation)  mean squared error against $d$, $n$ and $r$ respectively. These results are consistent with the main result in Theorem~\ref{ThmUpperMultiReg1}.

\subsection{Multivariate sparse auto-regressive models}
\label{SecSimVAR}
Now we consider dependent covariates and responses through the multivariate auto-regressive model. Recall that the generative model is:
$$
X^{(t+p)} = \sum_{j=1}^p B_{\cdot\cdot j} X^{(t+p-j)} +\epsilon^{(t)},
$$
where $1\leq t\leq n$ represents the time index, $1\leq j\leq p$ represents the lag index, $\{X^{(t)}\}_{t=0}^{(n+p)}$ is an $m$-dimensional vector, $\epsilon^{(t)}\sim \sigma \mathcal{N}(0, I_{m\times m})$ represents the additive noise.

We consider four different low-dimensional structures for $B$ and we choose the entries of $B$ to be sufficiently small to ensure the time series is stable.
\begin{itemize}
\item Slice-wise sparsity: $B_{\cdot\cdot 1},\ldots B_{\cdot\cdot s}$ are $s$ non-zero slices of diagonal matrix, where diagonal elements are constants $\rho$ with $\rho = 2$. $B_{\cdot\cdot s+1},\ldots B_{\cdot\cdot d_3}$ are zero slices.
\item Sparse low-rank slices: $B_{\cdot\cdot 1},\ldots B_{\cdot\cdot s}$ are $s$ non-zero slices, which are independent random rank-$r$ matrix (truncated matrix with i.i.d elements from $N(0, \tau^2)$). $B_{\cdot\cdot s+1},\ldots B_{\cdot\cdot d_3}$ are zero slices. Here $\tau = 1/150$ for $m = 10$ and $\tau = 1/500$ for $m = 20$.
\item Group Sparsity by lag (sparse normal slices):$B_{\cdot\cdot 1},\ldots B_{\cdot\cdot s}$ are $s$ non-zero slices, where elements follow i.i.d $N(0, \tau^2)$ with $\tau = 0.05$. $B_{\cdot\cdot s+1},\ldots B_{\cdot\cdot d_3}$ are zero slices.
\item Group sparsity by coordinate (sparse normal fibers): \\
$B_{s_1s_2\cdot}$ is a vector of i.i.d normal elements following $N(0, \tau^2)$ ($\tau = 0.1$) when $(s_1,s_2)\in \mathcal{S}$, which is a random sample of size $s$ from $\{ 1,\ldots, m\} \times \{1,\ldots, m\} $, and zero otherwise.
\end{itemize}

\begin{sidewaystable}
\centering
\begin{tabular}{|c|c|c|c|c|c|c|c|c|}
	\hline
	 Regularizer & coefficient tensor & $n$  & $m$ & $p$ & $s / r$  &$\sigma$ & SNR (sd) & rmse (sd )   \\
	\hline
	Vectorized  & $s$ diagonal slices & 2000 & 10 & 10 & $s = 5$  & 2 & 5.5 (1.4) &  0.43 (0.15)   \\
	\cline{2-9}
	Sparsity & $s$ diagonal slices & 2000 & 20 & 20 & $s = 5$ & 2 & 2.9 (0.7) &  0.25 (0.04)   \\
	\hline
	Low-rank & $s$ rank-$r$ slices & 2000 & 10 & 10 & $s = 10, r = 5$  & 0.05  & 1.0 (0.1) &  0.37 (0.01)   \\
	\cline{2-9}
	Slices & $s$ rank-$r$ slices & 2000 & 20 & 20 & $s = 10, r = 5$  & 0.05 & 1.2 (0.1) &  0.82 (0.02)   \\
	\hline
	Group Sparsity & $s$ Gaussian slices & 2000 & 10 & 10 & $s = 5$ & 0.2  & 0.42 (0.02) &  0.40 (0.01)   \\
	\cline{2-9}
	(by lag) & $s$ Gaussian slices & 2000 & 20 & 20 & $s = 5$ & 0.2 & 0.42 (0.01) &  0.64 (0.02)   \\
	\hline
	Group Sparsity & $s$ Gaussian fibers &  2000 & 10 & 10 & $s = 10$  & 0.02  & 0.93 (0.2) &  0.35 (0.04)   \\
	\cline{2-9}
	(by coordinate) & $s$ Gaussian fibers & 2000 & 20 & 20 & $s = 10$  & 0.02 & 1.2 (0.1) &  0.72 (0.07)   \\
	\hline
\end{tabular}
\caption{Multivariate auto-regressive model with various sparsity/low-rankness: RMSE were computed based on 50 simulations runs. Numbers in parentheses are standard errors.}
\label{table2}
\end{sidewaystable}

Table \ref{table2} shows the average rmse for 50 runs of each case as a function of $m$, $p$, $s$ and $r$. In general, the smaller the $n$ is, or the larger the $m$(or $p$) is, the harder it is to recover the coefficient $B$. These findings are consistent with our theoretical developments. 

\subsection{Pairwise interaction tensor models}
\label{SecSimPairwise}

Finally we consider the so-called pairwise interaction tensor models as described in Section~\ref{SecTheoryPairwise}. To implement this regularization scheme we kept iterating among the matrix slices $A_{1,2}$, $A_{1,3}$ and $A_{2,3}$ and updating one of the three at a time while assuming the other two components are fixed. For the  update of $A_{k_1,k_2}$, we conducted an approximated projection onto the zero-row-sum/zero-column-sum subspace after each generalized gradient descent (soft thresholding) step. 
$$A_{k_1, k_2}^{(i+1)} = \hat P (  \hat P_{\lambda \eta} (A_{k_1, k_2}^{(i)} -\eta \nabla f))$$
where $\eta$ is the step size for the gradient step, $\nabla f$ is the gradient of the least square objective function, $\hat P_{\lambda\eta}$ is the singular space soft-thresholding operator with threshold $\lambda\eta$, and $\hat P$ is the approximated projection operator that make any given matrix have zero row sums (by shifting rows) and zero column sums (by shifting columns). We simulated independent random low-rank matrix $C_{k_1, k_2}$s and make them have zero column sums and row sums by $B_{k_1,k_2} = \hat P (C_{k_1, k_2})$.

\begin{table}[htbp]
\centering
\begin{tabular}{|c|c|c|c|c|c|c|c|c|}
	\hline
	$n$  & $d_1$ & $d_2$ & $d_3$ & s & $\sigma$ & RMSE & SNR  \\
	\hline
	2000 & 40 & 40 & 40 & 5 & 10 & 0.54 (0.02) & 2.4 (0.1)   \\
	\hline
	2000 & 40 & 40 & 40 & 10 & 10 & 0.70 (0.01)  &  3.4 (0.1)   \\
	\hline
	2000 & 20 & 20 & 20 & 5 & 10 & 0.39 (0.01) & 1.7 (0.1)  \\
	\hline
	2000 & 20 & 20 & 20 & 10 & 10 & 0.37 (0.01) &  2.3 (0.1) \\
	\hline
	1000 & 20 & 20 & 20 & 5 & 10 & 0.58 (0.02) & 1.6 (0.1)   \\
	\hline
	1000 & 20 & 20 & 20 & 10 & 10 & 0.63 (0.02) &  2.3 (0.1)  \\
	\hline
\end{tabular}
\caption{Pairwise interaction model: RMSE were computed based on 50 simulations runs. Numbers in parentheses are standard errors.}
\label{table3}
\end{table}

Table~\ref{table3} shows the average (with standard deviation) rmse under different $r$, $d$, $n$ combinations under 50 runs. In general, the rmse in estimating the tensor coefficient increases as $s$ and $d$ increases.

\section{Proofs}

\label{SecProofs}
In this section, we present the proofs to our main results. We begin with the proof of Theorem~\ref{ThmUpper}.
\subsection{Proof of Theorem~\ref{ThmUpper}}
Our proof involves the following main steps:
\begin{itemize}
\item In the initial step, we use an argument similar to those developed in \cite{Neg10} to exploit weak decomposability and properties of the empirical risk minimizer and convex duals to upper bound $\|\widehat{T} - T\|_n$ in terms of $\mathcal{R}(\widehat{T} - T)$, and $\lambda$.
\item Next, we use properties of Gaussian random variables and supremum of Gaussian processes to express the lower bound on $\lambda$ in terms of the Gaussian width $\mathbb{E}[\mathcal{R}^*(G)]$ (Lemma~\ref{LemGauss} below).
\item The final and most challenging step involves proving a one-sided uniform law relating $\|\widehat{T} - T\|_n$ to $\|\widehat{T} - T\|_{\rm F}$ (Lemma~\ref{LemEmpLowerBound} below) which is analogous to restricted strong convexity. The proof for Lemma~\ref{LemEmpLowerBound} uses a novel truncation argument and is similar in spirit to that of Lemma 4 in ~\cite{RasWaiYu12}. Lemma~\ref{LemEmpLowerBound} is necessary to incorporate multiple responses as existing results relating the $\|\cdot\|_n$ to the population $\|\cdot\|_{\rm F}$ norm \citep[e.g.,][]{Mendelson10,RasWaiYu10b,vandeGeer} only apply to univariate functions.
\end{itemize} 
Throughout $\mathcal{R}(A)$ refers to the weakly decomposable regularizer over the tensor $A$. For a tensor $A$, we shall write $A_0$ and $A^\perp$ as its projections onto $\calA_0$ and $\calA^\perp$ with respect to the Frobenius norm, respectively.
 
Since $\widehat{T}$ is the empirical minimizer,
\begin{eqnarray*}
\frac{1}{2n} \sum_{i=1}^n \| Y^{(i)}- \langle X^{(i)} , \widehat{T} \rangle \|_{\rm F}^2 + \lambda \mathcal{R}( \widehat{T}) & \leq & \frac{1}{2n} \sum_{i=1}^n \| Y^{(i)}- \langle X^{(i)} , T \rangle \|_{\rm F}^2 + \lambda \mathcal{R}(T).
\end{eqnarray*}
Substituting $Y^{(i)}  = \langle X^{(i)} , T \rangle + \epsilon^{(i)}$ and $\Delta = \widehat{T} - T$,
\begin{eqnarray*}
\frac{1}{2n}\sum_{i=1}^n  \| \langle X^{(i)}, \Delta\rangle \|_{\rm F}^2 & \leq & \frac{1}{n} \left|\sum_{i=1}^{n} \langle \epsilon^{(i)} \otimes X^{(i)}, \Delta \rangle\right|  + \lambda (\mathcal{R}(T)  - \mathcal{R}(\widehat{T}))\\
& \le & \mathcal{R}^*\left(\frac{1}{n} \sum_{i=1}^n {\epsilon^{(i)} \otimes X^{(i)}} \right) \mathcal{R}(\Delta)  + \lambda (\mathcal{R}(T)  - \mathcal{R}(\widehat{T}_0) - c_\calR\mathcal{R}(\widehat{T}^\perp))\\
& \leq & \mathcal{R}^*\left(\frac{1}{n} \sum_{i=1}^n {\epsilon^{(i)} \otimes X^{(i)}} \right) \mathcal{R}(\Delta)  + \lambda (\mathcal{R}(\Delta_0) - c_\calR\mathcal{R}(\Delta^\perp)),
\end{eqnarray*}
where the second inequality follows from the decomposability and the last one follows from triangular inequality. 

Let $G \in \mathbb{R}^{d_1 \times d_2 \times ...\times d_N}$ be an tensor where each entry is i.i.d. $\mathcal{N}(0,1)$. Recall the definition of Gaussian width:
\begin{equation*}
w_G[\mathbb{B}_\calR(1)] = \mathbb{E}[\mathcal{R}^*(G)].
\end{equation*}
For simplicity let $$\eta_\calR = \frac{3+c_\calR}{2c_\calR}$$
and recall that $\lambda \geq {2 c_u \eta_{\calR}}n^{-1/2} \mathbb{E}[\mathcal{R}^*(G)]$. We have the following Lemma:

\blems
\label{LemGauss}
If $\lambda\geq {2 \sigma c_u \eta_{\calR}}n^{-1/2} \mathbb{E}[\mathcal{R}^*(G)]$, then 
\begin{equation*}
\lambda \geq \eta_{\calR} \mathcal{R}^*\left(\frac{1}{n} \sum_{i=1}^n {\epsilon^{(i)} \otimes X^{(i)}} \right),
\end{equation*}
with probability at least $1 - \exp\{-{\eta_{\calR}^2 \mathbb{E}[\mathcal{R}^*(G)]^2}/{4}\}$
\elems

The proof relies on Gaussian comparison inequalities and concentration inequalities.

\begin{proof}[Proof of Lemma~\ref{LemGauss}]

Recall that we have set: 
$$\lambda \geq \frac{2 \sigma \eta_\calR c_u }{\sqrt{n}} \mathbb{E}[\mathcal{R}^*(G)].$$
First we show that $\lambda \geq {\sigma c_u \eta_{\calR}}n^{-1/2} \mathcal{R}^*(G)$ with high probability using concentration of Lipschitz functions for Gaussian random variables (see Theorem~\ref{ThmConcGaussLip} in Appendix~\ref{AppTailBounds}). First we prove that $f(G) = \mathcal{R}^*(G) = \sup_{A \in \mathbb{B}_\calR(1)}\langle G, A \rangle$ is a $1$-Lipschitz function in terms of $G$. In particular note that:
\begin{eqnarray*}
f(G) - f(G') =  \sup_{A:\mathcal{R}(A) \leq 1} \langle G, A \rangle - \sup_{A: \mathcal{R}(A) \leq 1} \langle G', A \rangle.
\end{eqnarray*}
Let $\widetilde{A} := \argmax_{A: \mathcal{R}(A) \leq 1} \langle G, A \rangle$. Then 
\begin{eqnarray*}
\sup_{A:\mathcal{R}(A) \leq 1} \langle G, A \rangle - \sup_{A: \mathcal{R}(A) \leq 1} \langle G', A \rangle & = & \langle G, \widetilde{A} \rangle - \sup_{\mathcal{R}(A) \leq 1} \langle G', A \rangle\\
& \leq & \langle G, \widetilde{A} \rangle - \langle G', \widetilde{A} \rangle\\
& \leq & \langle G - G', \widetilde{A} \rangle \\
& \leq & \sup_{A:\mathcal{R}(A) \leq 1} \langle G - G', A \rangle \\
& \leq & \sup_{A:\|A\|_{\rm F} \leq 1} \langle G - G', A \rangle \\
& \leq & \|G - G'\|_{\rm F},
\end{eqnarray*}
where recall that $\|A\|_{\rm F} \leq \mathcal{R}(A)$ which implies the second last inequality. Therefore $f(G)$ is a $1$-Lipschitz function with respect to the Frobenius norm. Therefore, by applying Theorem~\ref{ThmConcGaussLip} in Appendix~\ref{AppTailBounds}, 
\begin{equation*}
\mathbb{P}\left\{\left|\sup_{A \in \mathbb{B}_\calR(1)}\langle G, A \rangle - \mathbb{E}[\sup_{A \in \mathbb{B}_\calR(1)}\langle G, A \rangle]\right| > w_G(\mathbb{B}_\calR(1)) \right\} \leq 2 \exp\left(-\frac{1}{2} w_G^2[\mathbb{B}_\calR(1)]\right).
\end{equation*}
Therefore $$\lambda \geq \frac{\eta_\calR \sigma c_u }{\sqrt{n}} \mathcal{R}^*(G)$$ with probability at least $1 - 2\exp\{-w_G^2[\mathbb{B}_\calR(1)]/2\}$.

To complete the proof, we use a Gaussian comparison inequality between the supremum of the process ${c_u }n^{-1/2} \langle G, A \rangle$ and $n^{-1} \sum_{i=1}^n \langle \epsilon^{(i)} \otimes X^{(i)}, A\rangle$ over the set $\mathbb{B}_\calR(1)$. Recall that:
\begin{equation*}
\mathcal{R}^*\left(\frac{1}{n} \sum_{i=1}^n {\epsilon^{(i)} \otimes X^{(i)}} \right) = \sup_{A\in \mathbb{B}_\calR(1)}\left\langle A, \frac{1}{n} \sum_{i=1}^n {\epsilon^{(i)} \otimes X^{(i)}} \right\rangle.
\end{equation*}

Recall that each $\epsilon^{(i)} \in \mathbb{R}^{d_{M+1} \times d_{M+2} \times ... \times d_N}$ is an i.i.d. standard centered Gaussian tensor with each entry having variance $\sigma^2$ and $\mbox{vec}(X) \in \mathbb{R}^{nd_1d_2\cdots d_M}$ is a Gaussian vector covariance $\Sigma \in \mathbb{R}^{(nD_M) \times (nD_M)}$. First we condition on all the $\epsilon^{(i)}$'s, which are indepedendent of the $X^{(i)}$'s. Further let $\{w^{(i)}: i=1,\ldots,n\}$ be i.i.d. standard normal Gaussian tensors where $w^{(i)} \in \mathbb{R}^{d_{1} \times d_{2} \times ... \times d_M}$. First we condition on all the $\epsilon^{(i)}$'s, which are indepedendent of the $X^{(i)}$'s. Assuming (\ref{AssCov}) and using a standard Gaussian comparison inequality due to Lemma~\ref{LemAnderson} in Appendix~\ref{AppTailBounds} proven earlier in~\cite{Anderson55}, if we condition on the $\epsilon^{(i)}$'s we get
\begin{equation*}
\mathbb{P}\left\{\sup_{A: \mathcal{R}(A) \leq 1}\frac{1}{n} \sum_{i=1}^n {\langle \epsilon^{(i)} \otimes X^{(i)}, A\rangle} > x \right\} \leq \mathbb{P}\left\{\sup_{A: \mathcal{R}(A) \leq 1}\frac{1}{n} \sum_{i=1}^n {\langle \epsilon^{(i)} \otimes w^{(i)}, A\rangle} > \frac{x}{c_u} \right\},
\end{equation*}
since
$${\rm Cov}({\rm vec}(X)) = \Sigma \preceq c_u^2 I_{(nD_M) \times (nD_M)}.$$
Now we define the $w_j \in \mathbb{R}^n$ as the standard random vector where $1 \leq j \leq D_M$ and $w_j = (w_j^{(1)},...,w_j^{(n)})$. Conditioning on the $w^{(i)}$'s and dealing with the randomness in the $\epsilon^{(i)}$'s, 
$$
\frac{1}{n} \sum_{i=1}^n {\langle \epsilon^{(i)} \otimes w^{(i)}, A\rangle}  \leq \max_{1 \leq j \leq D_M}\frac{\|w_j\|_{\ell_2}}{\sqrt{n}} \frac{\sigma}{\sqrt{n}}\langle G, A \rangle, 
$$
where $G \in \mathbb{R}^{d_1 \times d_2 \times ... \times d_N}$ is an i.i.d. standard normal tensor, since the $\epsilon^{(i)}$'s are i.i.d. standard normal. Now we upper bound 
$$
 \max_{1 \leq j \leq D_M}\frac{\|w_j\|_{\ell_2}}{\sqrt{n}},
$$
using standard $\chi^2$ tail bounds. Since $\|w_j\|_{\ell_2}^2$ is a $\chi^2$ random variable with $n$ degrees of freedom, for each $j$,
$$
\mathbb{P}\left\{\frac{\|w_j\|_{\ell_2}^2}{n} \geq 4\right\} \leq \exp(-n),
$$
using $\chi^2$ tail bounds provided in Appendix~\ref{AppTailBounds} presented in~\cite{LauMas98}. Now taking the union bound over $D_M$,
$$
\mathbb{P}\left\{ \max_{1 \leq j \leq D_M}\frac{\|w_j\|_{\ell_2}}{\sqrt{n}} \geq 2\right\} \leq \exp(\log D_M -n),
$$ 
and provided $n \geq 2 \log D_M$, it follows that with probability greater than $1 - \exp(-n/2)$, 
$
 \max_{1 \leq j \leq D_M}\frac{\|w_j\|_{\ell_2}}{\sqrt{n}} \leq 2. 
$
Therefore, with probability at least $1 - \exp(-n/2)$, 
$$
\frac{1}{n} \sum_{i=1}^n {\langle \epsilon^{(i)} \otimes w^{(i)}, A\rangle}  \leq \frac{2\sigma}{\sqrt{n}}\langle G, A \rangle.
$$

Now we apply Slepian's lemma \citep{Slepian62} to complete the proof. Slepian's lemma is stated  in Appendix~\ref{AppTailBounds}. Applying Slepian's lemma (Lemma~\ref{LemSlepian} in Appendix~\ref{AppTailBounds}),
\begin{equation*}
\mathbb{P}\left\{\sup_{\mathcal{R}(A)\leq 1}\frac{1}{n} \sum_{i=1}^n {\langle \epsilon^{(i)} \otimes w^{(i)}, A\rangle} > x \right\} \leq \mathbb{P}\left\{\sup_{\mathcal{R}(A)\leq 1}\frac{2\sigma}{\sqrt{n}} \langle G, A\rangle > x \right\},
\end{equation*}
for all $x > 0$. Substituting $x$ by $x/c_u$ means that
\begin{equation*}
\mathbb{P}\left\{\mathcal{R}^*\left(\frac{1}{n} \sum_{i=1}^n {\langle \epsilon^{(i)} \otimes w^{(i)}, A \rangle}\right)  > x \right\} \leq \mathbb{P}\left\{ \frac{2 \sigma c_u}{\sqrt{n}} \mathcal{R}^*(G) > x \right\},
\end{equation*}
for any $x > 0$. This completes the proof.
\end{proof}
\vskip 25pt

In light of Lemma~\ref{LemGauss}, for the remainder of the proof, we can condition on the event that
$$\lambda \geq \eta_{\calR} \mathcal{R}^*\left(\frac{1}{n} \sum_{i=1}^n {\epsilon^{(i)} \otimes X^{(i)}} \right).$$
Under this event,
\begin{eqnarray*}
\frac{1}{2n}\sum_{i=1}^n  \| \langle X^{(i)}, \Delta \rangle \|_{\rm F}^2 & \leq &  \frac{1}{\eta_\calR}\lambda \mathcal{R}(\Delta)  +  \lambda (\mathcal{R}(\Delta_0) - c_\calR\mathcal{R}(\Delta^\perp))\\
&\leq & \left(1+\frac{1}{\eta_\calR}\right)\lambda \mathcal{R}(\Delta_0) - \left(c_\calR - \frac{1}{\eta_\calR}\right) \lambda \mathcal{R}(\Delta^\perp).
\end{eqnarray*}
Since
$$\frac{1}{2n}\sum_{i=1}^n  \| \langle \Delta, X^{(i)} \rangle \|_{\rm F}^2 \geq 0,$$
we get
\begin{equation*}
\mathcal{R}(\Delta^\perp) \leq \frac{3}{c_\calR} \mathcal{R}(\Delta_{0}).
\end{equation*}

Hence we define the cone
$$\mathcal{C} = \left\{ \Delta\;|\; \mathcal{R}(\Delta^\perp) \leq 3c_\calR^{-1} \mathcal{R}(\Delta_{0}) \right\},$$
and know that $\Delta \in \mathcal{C}$. Hence
$$
\frac{1}{2n}\sum_{i=1}^n  \| \langle X^{(i)},\Delta \rangle \|_{\rm F}^2\leq \frac{3(1+c_\calR)}{3+c_\calR} \lambda \mathcal{R}(\Delta_{0})\leq \frac{3(1+c_\calR)}{3+c_\calR} \sqrt{s(\calA)} \lambda \|\Delta\|_{\rm F}.
$$
Recall that
$$\frac{1}{n}\sum_{i=1}^n  \| \langle X^{(i)}, \Delta\rangle \|_{\rm F}^2 = \|\Delta\|_n^2.$$
Thus,
\begin{eqnarray*}
\|\Delta\|_n^2  & \leq & \frac{6(1+c_\calR)}{3+c_\calR}\sqrt{s(\mathcal{A})} \lambda \|\Delta\|_{\rm F}.
\end{eqnarray*}
For convenience, in the remainder of this proof let $$\delta_n:= \frac{6(1+c_\calR)}{3+c_\calR}\sqrt{s(\mathcal{A})} \lambda.$$

Now we split into three cases. (i) If $\|\Delta\|_n \geq \| \Delta \|_{\rm F}$, then $$\max\{\| \Delta\|_n, \| \Delta\|_{\rm F}\} \leq \delta_n.$$
On the other hand if (ii) $\|\Delta\|_n \leq \| \Delta \|_{\rm F}$ and $\|\Delta\|_{\rm F} \leq \frac{c_u}{c_{\ell}} \delta_n$, then
$$\max\{\| \Delta\|_n, \| \Delta\|_{\rm F}\} \leq \frac{c_u}{c_{\ell}} \delta_n.$$
Hence the only case we need to consider is (iii) $\|\Delta\|_n \leq \| \Delta \|_{\rm F}$ and $\|\Delta\|_{\rm F} \geq c_uc_{\ell}^{-1} \delta_n$. Now we follow a similar proof technique to the proof for Theorem 1 in \cite{RasWaiYu12}.

Let us define the following set:
\begin{eqnarray*}
\mathcal{C}(\delta_n) := \left\{ \Delta \in \mathbb{R}^{d_1 \times d_2 \times\cdots\times d_N}\;|\;\mathcal{R}(\Delta^\perp) \leq 3 c_\calR^{-1} \mathcal{R}(\Delta_{0}),\;\|\Delta\|_n \leq \|\Delta\|_{\rm F} \right\}.
\end{eqnarray*}

Further, let us define the event:
\begin{equation*}
\mathcal{E}(\delta_n) := \left\{ \|\Delta\|_n^2 \geq \frac{1}{4}\|\Delta\|_{\rm F}^2 \;|\; \Delta \in \mathcal{C}(\delta_n),\; \|\Delta\|_{\rm F} \geq \frac{c_u}{c_\ell}\delta_n \right\}.
\end{equation*}
Let us define the alternative event:
\begin{equation*}
\mathcal{E}'(\delta_n) := \{ \|\Delta\|_{n}^2 \geq \frac{1}{4} \|\Delta\|_{\rm F}^2 \;|\; \Delta \in \mathcal{C}(\delta_n),\; \|\Delta\|_{\rm F} =  \frac{c_u}{c_\ell} \delta_n \}.
\end{equation*}
We claim that it suffices to show that $\mathcal{E}'(\delta_n)$ holds with probability at least $1 - \exp(-c n )$ for some constant $c > 0$. In particular, given an arbitrary non-zero $\Delta \in \mathcal{C}(\delta_n)$, consider the re-scaled tensor
$$\widetilde{\Delta} = \frac{c_u\delta_n}{c_\ell} \frac{\Delta}{\|\Delta\|_{\rm F}}.$$
Since $\Delta \in \mathcal{C}(\delta_n)$ and $\mathcal{C}(\delta_n)$ is star-shaped, we have $\widetilde{\Delta} \in \mathcal{C}(\delta_n)$ and $\|\widetilde{\Delta}\|_{\rm F} = c_uc_\ell^{-1}\delta_n$ by construction. Consequently, it is sufficient to prove that $\mathcal{E}'(\delta_n)$ holds with high probability.

\blems
\label{LemEmpLowerBound}
Assume that for any $c> 0$, there exists an $n$ such that $\sqrt{s} \lambda \leq c$. Then there exists a $\tilde{c} > 0$ such that
\begin{equation*}
\mathbb{P}\big(\mathcal{E}'(\delta_n)\big) \geq 1 - \exp(-\tilde{c}n).
\end{equation*}
\elems

\begin{proof}[Proof of Lemma~\ref{LemEmpLowerBound}]
Denote by $D_N=d_1d_2\cdots d_N$ and $D_M=d_1d_2\cdots d_M$. Now we define the random variable
$$Z_n(\mathcal{C}(\delta_n)) = \sup_{\Delta \in \mathcal{C}(\delta_n)} \left\{\frac{c_u^2}{c_\ell^2}\delta_n^2 - \frac{1}{n} \sum_{i=1}^n \|\langle \Delta, X^{(i)} \rangle\|_{F}^2\right\},$$
then it suffices to show that
$$Z_n(\mathcal{C}(\delta_n)) \leq \frac{c_u^2\delta_n^2}{2 c_\ell^2}.$$
Recall that the norm
$$\|\Delta\|_n^2 = \frac{1}{n} \sum_{i=1}^n \|\langle \Delta, X^{(i)} \rangle\|_{F}^2.$$ 
To expand this out, reacall $[M] = \{1,2,...,M\}$ and define an extension of the standard matricization
$$
\tilde{\Delta} := \mathcal{M}_{[M]}(\Delta) \in \mathbb{R}^{D_M \times D_N/D_M},
$$
which groups together the first $M$ modes. With a slight abuse of notation, it follows that
$$\|\Delta\|_n^2 = \frac{1}{n} \sum_{i=1}^n \sum_{m=1}^{D_N/D_M}\langle \tilde{\Delta}_m, \mbox{vec}(X^{(i)})\rangle^2,$$
where $\tilde{\Delta}_m \in \mathbb{R}^{D_M}$, and clearly $\mbox{vec}(X^{(i)}) \in \mathbb{R}^{D_M}$.
In order to complete the proof we make use of a truncation argument. For a constant $\tau > 0$ to be chosen later, consider the truncated quadratic function
$$
\phi_{\tau}(u) = \min\{u^2,\tau^2\},
$$ 
and define 
$$\Delta_{m, \tau}(X) = \mbox{sign}\left(\langle \tilde{\Delta}_m, \mbox{vec}(X) \rangle)\sqrt{\phi_{\tau}(\langle \tilde{\Delta}_m, \mbox{vec}(X) \rangle)}\right)$$
where $X \in \mathbb{R}^{d_1 \times d_2 \times...\times d_M}$ is an input tensor. Further let
$$\Delta_{\tau}(X) = (\Delta_{1,\tau}(X),\Delta_{2,\tau}(X),...,\Delta_{D_N/D_M, \tau}(X))$$
and
$$
\|\Delta_{m, \tau}\|_{\rm F}^2 = \mathbb{E}[\Delta^{2}_{m, \tau}(X)],
$$
$\|\Delta_{\tau}\|_{\rm F}^2 = \sum_{m=1}^{D_N/D_M}{\|\Delta_{m, \tau}\|_{\rm F}^2}$, and similarly
$$
\|\Delta_{m, \tau}\|_n^2 = \frac{1}{n}\sum_{i=1}^{n}{\Delta^{2}_{m, \tau}(X^{(i)})},
$$
$\|\Delta_{\tau}\|_n^2 = \sum_{m=1}^{D_N/D_M}{\|\Delta_{m, \tau}\|_n^2}$.
By construction, for any $\Delta \in  \mathcal{C}(\delta_n)$, $\|\Delta\|_n^2 \geq \|\Delta_{\tau}\|_n^2$ and hence
$$
\|\Delta\|_n^2 \geq \|\Delta_{\tau}\|_{\rm F}^2 - \sup_{\Delta \in \mathcal{C}(\delta_n)}|\|\Delta_{\tau}\|_n^2 - \|\Delta_{\tau}\|_{\rm F}^2|
$$
The remainder of the proof consists of showing that for a suitable of $\tau$,
$$
\|\Delta_{\tau}\|_{\rm F}^2 \geq \frac{3}{4}\|\Delta\|_{\rm F}^2,\;\;\mbox{for all}\; \Delta \in  \mathcal{C}(\delta_n)
$$
and
$$
\mathbb{P}\left\{Z_n \geq \frac{c_u^2}{4c_{\ell}^2}\delta_n^2  \right\} \leq c_1 e^{-c_2 n \delta_n^2},
$$
where $Z_n :=\sup_{\Delta \in \mathcal{C}(\delta_n)}|\|\Delta_{\tau}\|_n^2 - \|\Delta_{\tau}\|_{\rm F}^2|$. By definition
\begin{eqnarray*}
\|\Delta_m\|_{\rm F}^2 - \|\Delta_{m,\tau}\|_{\rm F}^2 & \leq & \mathbb{E}\big[\langle \tilde{\Delta}_m, \mbox{vec}(X) \rangle^2 \mathbf{1}[|\tilde{\Delta}_m, \mbox{vec}(X) \rangle| \geq \tau]  \big] \\ &\leq&  \sqrt{\mathbb{E}[\langle \tilde{\Delta}_m, \mbox{vec}(X) \rangle^4]}\sqrt{\mathbb{P}[|\tilde{\Delta}_m, \mbox{vec}(X) \rangle| \geq \tau]} \\
& \leq & \frac{\mathbb{E}[\langle \tilde{\Delta}_m, \mbox{vec}(X) \rangle^4]}{\tau^2},
\end{eqnarray*}
where the second last inequality follows from the Cauchy-Schwartz inequality and the final inequality follows from Markov's inequality. Since $\langle \tilde{\Delta}_m, \mbox{vec}(X) \rangle$ is a Gaussian random variable,
$$
\mathbb{E}[\langle \tilde{\Delta}_m, \mbox{vec}(X) \rangle^4] \leq 3 \mathbb{E}[\langle \tilde{\Delta}_m, \mbox{vec}(X) \rangle^2] = 3 \|\Delta_m\|_{\rm F}^2.
$$
Therefore
$$
\|\Delta_m\|_{\rm F}^2 - \|\Delta_{m,\tau}\|_{\rm F}^2 \leq 9 \frac{ \|\Delta_m\|_{\rm F}^2}{\tau^2}.
$$
Setting $\tau^2 = 36$, $\|\Delta_m\|_{\rm F}^2 - \|\Delta_{m,\tau}\|_{\rm F}^2 \leq \frac{1}{4}\|\Delta_m\|_{\rm F}^2$. Summing over $m$ implies
$$
\|\Delta\|_{\rm F}^2 - \|\Delta_{\tau}\|_{\rm F}^2 \leq \frac{1}{4}\|\Delta\|_{\rm F}^2,
$$
which implies $\|\Delta_{\tau}\|_{\rm F}^2 \geq \frac{3}{4} \|\Delta\|_{\rm F}^2$. Now to prove the high probability bound on $Z_n$ by first upper bounding $E[Z_n]$. A standard symmetrization argument \citep[see e.g.,][]{Pollard84}, shows that
$$
\mathbb{E}_X[Z_n] \leq 2 \mathbb{E}_{X, z}\left[\sup_{\Delta \in \mathcal{C}(\delta_n)}\left|\frac{1}{n}\sum_{i=1}^n \sum_{m=1}^{D_N/D_M}z_m^{(i)}\phi_{\tau}(\langle\tilde{\Delta}_m, \mbox{vec}(X^{(i)}) \rangle)\right|\right],
$$
where $z_m^{(i)}$ are i.i.d. Rademacher random variables (that is $P(z_m^{(i)}=+1) = P(z_m^{(i)}=-1)  = 1/2$). Since $\phi_{\tau}(\langle\tilde{\Delta}_m, \mbox{vec}(X^{(i)}) \rangle)$ is a Lipschitz function with Lipschitz constant $2 \tau$,  the Ledoux-Talagrand contrtaction inequality \citep{LedTal91} implies that ,
\begin{eqnarray*}
\mathbb{E}_X[Z_n] &\leq& 2 \mathbb{E}_{X, z}\left[\sup_{\Delta \in \mathcal{C}(\delta_n)}\left|\frac{1}{n}\sum_{i=1}^n \sum_{m=1}^{D_N/D_M}z_m^{(i)}\phi_{\tau}(\langle\tilde{\Delta}_m, \mbox{vec}(X^{(i)}) \rangle)\right|\right] \\
&\leq& 8\tau \mathbb{E}_{X, z}\left[\sup_{\Delta \in \mathcal{C}(\delta_n)}\left|\frac{1}{n}\sum_{i=1}^n \sum_{m=1}^{D_N/D_M}z_m^{(i)}\langle\tilde{\Delta}_m, \mbox{vec}(X^{(i)}) \rangle\right|\right].  
\end{eqnarray*}
Using standard comparisons between Rademacher and Guassian complexities \citep[see, e.g., Lemma 4 of][]{Bartlett02}, there exists a $C > 0$ such that
\begin{eqnarray*}
&&\mathbb{E}_{X, z}\left[\sup_{\Delta \in \mathcal{C}(\delta_n)}\left|\frac{1}{n}\sum_{i=1}^n \sum_{m=1}^{D_N/D_M}z_m^{(i)}\langle\tilde{\Delta}_m, \mbox{vec}(X^{(i)}) \rangle\right|\right] \\
&\leq& C\mathbb{E}_{X, w}\left[\sup_{\Delta \in \mathcal{C}(\delta_n)}\left|\frac{1}{n}\sum_{i=1}^n \sum_{m=1}^{D_N/D_M}w_m^{(i)}\langle\tilde{\Delta}_m, \mbox{vec}(X^{(i)}) \rangle\right|\right],
\end{eqnarray*}
where $w_m^{(i)}$s ($1 \leq i \leq n$ and $1 \leq m \leq D_N/D_M$) are independent standard normal random variables. 

Next we upper bound the Gaussian complexity
$$
\mathbb{E}_{w} \left(\sup_{\Delta \in \mathcal{C}(\delta_n)} \frac{1}{n} \sum_{i=1}^n \langle w^{(i)} \otimes X^{(i)}, \Delta \rangle \right).$$
Clearly,
$$
\frac{1}{n} \sum_{i=1}^n \langle w^{(i)} \otimes X^{(i)}, \Delta \rangle \leq \mathcal{R}^*\left(\frac{1}{n} \sum_{i=1}^n w^{(i)} \otimes X^{(i)} \right) \mathcal{R}(\Delta) \leq \frac{\lambda}{\eta_{\mathcal{R}}} \mathcal{R}(\Delta).
$$
by the definition of $\lambda$ and our earlier argument. Since $\Delta \in \mathcal{C}(\delta_n)$,
$$
\frac{\lambda}{\eta_{\mathcal{R}}} \mathcal{R}(\Delta) \leq \frac{\lambda(1 + 3c^{-1}_{\mathcal{R}})}{\eta_{\mathcal{R}}} \mathcal{R}(\Delta_0) \leq \frac{\lambda(1 + 3c^{-1}_{\mathcal{R}})}{\eta_{\mathcal{R}}} \sqrt{s(\mathcal{A})}\|\Delta_0\|_{\rm F} \leq \frac{c_u(1 + 3c^{-1}_{\mathcal{R}})}{c_\ell \eta_{\mathcal{R}}}\delta_n \sqrt{s(\mathcal{A})}\lambda.
$$
Therefore,
$$\mathbb{E}_{w} \left(\sup_{\Delta \in \mathcal{C}(\delta_n)} \frac{1}{n} \sum_{i=1}^n \langle w^{(i)} \otimes X^{(i)}, \Delta \rangle \right) \leq \frac{c_u(1 + 3c^{-1}_{\mathcal{R}})}{c_\ell \eta_{\mathcal{R}}}\delta_n \sqrt{s(\mathcal{A})}\lambda.
$$
Since
$$\frac{c_u(1 + 3c^{-1}_{\mathcal{R}})}{c_\ell \eta_{\mathcal{R}}} \sqrt{s(\mathcal{A})}\lambda = \frac{2 c_u}{c_{\ell}} \delta_n,$$
we have
$$
\mathbb{E}_{w} \left(\sup_{\Delta \in \mathcal{C}(\delta_n)} \frac{1}{n} \sum_{i=1}^n \langle w^{(i)} \otimes X^{(i)}, \Delta \rangle \right) \leq \frac{c_u}{c_\ell \eta_{\mathcal{R}}}\delta_n^2.
$$
Finally we need a concentration bound to show that
$$
\mathbb{P}\left\{Z_n \geq \frac{c_u^2}{4c_{\ell}^2}\delta_n^2  \right\} \leq c_1 e^{-c_2 n \delta_n^2}.
$$
In particular using Talagrand's theorem for empirical processes \citep{Tal96}. By construction $\phi_{\tau}(.) \leq \tau^2 = 36$ and
\begin{eqnarray*}
\mbox{Var}\left(\sum_{m=1}^{D_N/D_M}\phi_{\tau}(\langle \tilde{\Delta}_m, \mbox{vec}(X) \rangle)\right) &\leq& \sum_{m=1}^{D_N/D_M}\mathbb{E}[\phi_{\tau}^4(\langle \tilde{\Delta}_m, \mbox{vec}(X) \rangle)]\\
&\leq& \tau^2\sum_{m=1}^{D_N/D_M}{\|\tilde{\Delta}_m\|_{\ell_2}^2} \leq \frac{36 c_u^2 \delta_n^2}{c_{\ell}^2}. 
\end{eqnarray*}
Consequently Talagrand's inequality implies that
$$
\mathbb{P}\big(Z_n \geq E[Z_n] + u\big) \leq c_1 \exp\left(-\frac{c_2 n u^2}{3\delta_n^2 + 9u^2} \right).
$$
Since $\mathbb{E}[Z_n] \leq \frac{c_u}{c_\ell \eta_{\mathcal{R}}}\delta_n^2$, the claim follows by setting $u =  \frac{c_u}{c_\ell \eta_{\mathcal{R}}}\delta_n^2$. 
\end{proof}
\vskip 25pt

Finally we return to the main proof. On the event $\mathcal{E}(\delta_n)$, it now follows easily that,
\begin{equation*}
\max\{\|\Delta\|_{\ell_2}^2, \|\Delta\|_n^2\} \leq  \frac{\eta_{\mathcal{R}} c_u^2}{c_\ell^2}  s(\mathcal{A})\lambda^2.
\end{equation*}
This completes the proof for Theorem~\ref{ThmUpper}.

\subsection{Proof of other results in Section~\ref{SecBounds}}

In this section we present proofs for the other main results from Section \ref{SecBounds}, deferring the more technical parts to the appendix. 

\begin{proof}[Proof of Lemmas ~\ref{LemSparsity}, \ref{LemSparsityfiber} and \ref{LemSparsityslice}]
We prove these three lemmas together since the proofs follow a very similar argument. First let $S \subset \{1,2,3\}$ denote the directions in which sparsity is applied and $D_S = \prod_{k \in S}{d_k}$ denote the total dimension in all these directions. For example, in Lemma~\ref{LemSparsity} $S = \{1,2,3\}$ and $D_S = d_1d_2d_3$, for Lemma~\ref{LemSparsityfiber}, $S = \{2,3\}$ and $D_S = d_2d_3$ and for Lemma~\ref{LemSparsityslice}, $S =\{1\}$ and $D_S = d_1$. Recall $N = \{1,2,3\}$ and $D_N  = d_1d_2d_3$.

Note that $\mathcal{R}^{*}(G)$ can be represented by the variational form:
\begin{equation*}
\mathcal{R}^{*}(G) = \sup_{\|\mbox{vec}(u)\|_{\ell_1} \leq 1, \|v\|_{\rm F} \leq 1}  \langle G, u \otimes v \rangle,
\end{equation*}
where $u \in \mathbb{R}^{d_{S_1} \times ...\times d_{S_{|S|}}}$ and $v \in \mathbb{R}^{d_{S^c_1} \times ...\times d_{S^c_{N-|S|}}}$. Now we express the supremum of this Gaussian process as:
\begin{equation*}
\sup_{(u,v) \in V} \mbox{vec}(u)^\top \mathcal{M}_S(G) \mbox{vec}(v),
\end{equation*}
where recall $\mathcal{M}_S$ is the matricization involving either slice or fiber $S$. The remainder of the proof follows from Lemma~\ref{LemSupGauss1} in Appendix~\ref{AppSupGauss}.
\fpro
\vskip 25pt

\begin{proof}[Proof of Lemma~\ref{LemSparsityslice2}]
Recall that 
\begin{equation*}
\mathcal{R}^{*}(G) := \max_{1 \leq j_3 \leq d_3} \left\|G_{..j_3}\right\|_s.
\end{equation*}
For each $1 \leq j_3 \leq d_3$, Lemma~\ref{LemGaussTens} in Appendix~\ref{AppSupGauss} with $N = 2$ satisfies the concentration inequality 
\begin{equation*}
\mathbb{E}[\|G_{..j_3}\|_s]  \leq \sqrt{6(d_1 + d_2)}.
\end{equation*}
Applying standard bounds on the maximum of functions of independent Gaussian random variables,
\begin{equation*}
\mathbb{E}[\max_{1 \leq j_3 \leq d_3}\|G..j_3\|_s]  \leq \sqrt{6(d_1 + d_2 + \log d_3)}.
\end{equation*}
This completes the proof.
\end{proof}
\vskip 25pt

\spro[Proof of Lemma \ref{le:sglasso2}]
Using the standard nuclear norm upper bound for a matrix in terms of rank and Frobenius norm:
\begin{eqnarray*}
\mathcal{R}_4^2(A) & = & \left(\sum_{j_3=1}^{d_3}\|A_{\cdot \cdot j_3}\|_\ast\right)^2 \\
& \leq &  \left(\sum_{j_3=1}^{d_3}\sqrt{{\rm rank}(A_{\cdot \cdot j_3})}\|A_{\cdot \cdot j_3}\|_{\rm F}\right)^2\\
& \leq & \sum_{j_3=1}^{d_3}{\rm rank}(A_{\cdot \cdot j_3}) \sum_{j_3=1}^{d_3}\|A_{\cdot \cdot j_3}\|_{\rm F}^2 = \sum_{j_3=1}^{d_3}{\rm rank}(A_{\cdot \cdot j_3}) \|A\|_{\rm F}^2,  
\end{eqnarray*}
where the final inequality follows from the Cauchy-Schwarz inequality. Finally, note that for any $A \in \Theta_4(r)/\{0\}$,
$$\sum_{j_3=1}^{d_3}{\rm rank}(A_{\cdot \cdot j_3})\le r,$$
which completes the proof.
\fpro
\vskip 25pt

\begin{proof}[Proof of Lemma \ref{LemNuclearNorm}]
Note that $\mathcal{R}^*(G) = \|G\|_s$, we can directly apply Lemma~\ref{LemGaussTens} with $N=3$ from Appendix~\ref{AppSupGauss}.
\end{proof}
\vskip 25pt

\spro[Proof of Lemma \ref{LemNuclearNorm1}]
From Tucker decomposition (\ref{eq:tucker}), it is clear that for any $A\in \Theta_5(r)$, we can find sets of vectors $\{u_k: k=1,\ldots, r^2\}$, $\{v_k: k=1,\ldots, r^2\}$ and $\{w_k: k=1,\ldots, r^2\}$ such that
$$
A= \sum_{k=1}^{r^2}{u_k \otimes v_k \otimes w_k},
$$
and in addition,
$$
u_{k}^\top u_{k'}=(v_{k}^\top v_{k'})(w_{k}^\top w_{k'})=0
$$
for any $k\neq k'$. It is not hard to see that
$$
\|A\|_{\rm F}^2=\sum_{k=1}^{r^2}\left(\|u_k\|_{\ell_2}^2\|v_k\|_{\ell_2}^2\|w_k\|_{\ell_2}^2\right).
$$
On the other hand, as shown by \cite{YuanZhang14},
$$
\|A\|_{\ast}=\sum_{k=1}^{r^2}\left(\|u_k\|_{\ell_2}\|v_k\|_{\ell_2}\|w_k\|_{\ell_2}\right).
$$
The claim then follows from an application of Cauchy-Schwartz inequality.
\fpro
\vskip 25pt

\begin{proof}[Proof of Lemma~\ref{LemLowRank}]
Recall that we are considering the regularizer
$$
\calR_6^\ast(A)=3\max\left\{\|\calM_1(A)\|_s, \|\calM_2(A)\|_s, \|\calM_3(A)\|_s\right\},
$$
and our goal is to upper bound 
$$
\mathcal{R}^*_6(G)=3\max_{1\le k\le 3}\|\calM_k(G)\|_s.
$$
Once again apply Lemma~\ref{LemGaussTens} in Appendix~\ref{AppSupGauss} with $N=2$ for each matricization implies
$$
\mathbb{E}[\mathcal{R}^*_6(G)] \leq 4 \max(\sqrt{d_1},\sqrt{d_2},\sqrt{d_3}).
$$
\end{proof}	
\vskip 25pt

\spro[Proof of Lemma \ref{LemLowRank1}]
It is not hard to see that
\begin{eqnarray*}
\calR_6(A)^2 & = & {1\over 9}\left(\|\calM_1(A)\|_\ast+\|\calM_2(A)\|_\ast+\|\calM_3(A)\|_\ast\right)^2 \\
& \leq & {1\over 9}(\sqrt{r_1} + \sqrt{r_2} + \sqrt{r_3})^2 \|A\|_{\rm F}^2\\
& \leq & \max\{r_1(A), r_2(A), r_3(A)\}\|A\|_{\rm F}^2,
\end{eqnarray*}
which completes the proof.
\fpro

\subsection{Proof of results in Section~\ref{SecExamples}}

In this section we prove the results in Section~\ref{SecExamples}. First we provide a general minimax lower result that we apply to our main results. Let $\mathcal{T} \subset \mathbb{R}^{d_1 \times d_2 \times\cdots\times d_N}$ be an arbitrary subspace of order-$N$ tensors. 

\btheos
\label{ThmLower}
Assume that (\ref{AssCov}) holds and there exists a finite set $\{A^1, A^2,\ldots,A^{m} \} \in \mathcal{T}$ of tensors such that $\log m \geq 128 n \delta^2$, such that 
\begin{equation*}
n \sigma^2 c_u^{-2} \delta^2 \leq \|A^{\ell_1} - A^{\ell_2}\|_{\rm F}^2 \leq 8 n \sigma^2 c_u^{-2} \delta^2 ,
\end{equation*}
for all $\ell_1\neq \ell_2 \in [m]$ and all $\delta > 0$. Then
\begin{equation*}
\min_{\widetilde{T}} \max_{T \in \mathcal{T}} \| \widetilde{T} - T \|_{\rm F}^2  \geq c \sigma^2 c_u^{-2} \delta^2,
\end{equation*}
with probability at least $1/2$ for some $c >0$.
\etheos

\spro
We use standard information-theoretic techniques developed in \cite{IbrHas81} and extended in ~\cite{YanBar99}. Let $\{A^1, A^2,\ldots,A^m \}$ be a set such that
$$\|A^{\ell_1} - A^{\ell_2}\|_{\rm F}^2 \geq n \sigma^2 c_u^{-2} \delta^2$$
for all $\ell_1 \neq \ell_2$, and let $\widetilde{m}$ be a random variable uniformly distributed over the index set $[m] = \{1,2,\ldots,m\}$. 

Now we use a standard argument which allows us to provide a minimax lower bound in terms of the probability of error in a multiple hypothesis testing problem \citep[see, e.g.,][]{YanBar99,Yu} then yields the lower bound. 
\begin{eqnarray*}
\inf_{\widetilde{T}} \sup_{T \in \mathcal{T}} \mathbb{P}\left\{\|\widetilde{T} - T\|_{\rm F}^2 \geq \frac{\sigma^2 c_u^{-2} \delta^2}{2} \right\} \geq \inf_{\widetilde{T}} \mathbb{P}(\widetilde{T} \neq A^{\widetilde{m}})
\end{eqnarray*}
where the infimum is taken over all estimators $\widetilde{T}$ that are measurable functions of $X$ and $Y$. 

Let $X = \{ X^{(i)}: i=1,\ldots,n\}$, $Y = \{ Y^{(i)}: i=1,\ldots,n \}$ and $E = \{ \epsilon^{(i)}: i=1,\ldots,n \}$. Using Fano's inequality \citep[see, e.g.,][]{Cover}, for any estimator $\widetilde{T}$, we have:
\begin{equation*}
\mathbb{P}[\widetilde{T} \neq A^{\widetilde{m}} | X ] \geq 1 - \frac{I_X(A^{\widetilde{m}}; Y)+\log 2}{\log m}.
\end{equation*}
Taking expectations over $X$ on both sides, we have
\begin{equation*}
\mathbb{P}[\widetilde{T} \neq A^{\widetilde{m}}] \geq 1 - \frac{\mathbb{E}_X[I_X(A^{\widetilde{m}}; Y)]+\log 2}{\log m}.
\end{equation*}

For $\ell = 1,2,\ldots,m$, let $\mathbb{Q}^{\ell}$ denote the condition distribution of $Y$ conditioned on $X$ and the event $\{T = A^{\ell}\}$, and $D_{\rm KL}(\mathbb{Q}^{\ell_1} ||\mathbb{Q}^{\ell_2} )$ denote the Kullback-Leibler divergence between $\mathbb{Q}^{\ell_1}$ and $\mathbb{Q}^{\ell_2}$. From the convexity of mutual information \citep[see, e.g.,][]{Cover}, we have the upper bound
$$I_X(T; Y) \leq \frac{1}{{m \choose 2}} \sum_{\ell_1,\ell_2 = 1}^m {D_{\rm KL}(\mathbb{Q}^{\ell_1} ||\mathbb{Q}^{\ell_2} )}.$$
Given our linear Gaussian observation model~\eqref{eq:model},
\begin{equation*}
D_{\rm KL}(\mathbb{Q}^{\ell_1} ||\mathbb{Q}^{\ell_2} ) = \frac{1}{2\sigma^2} \sum_{i=1}^n \left(\langle A^{\ell_1}, X^{(i)} \rangle - \langle A^{\ell_2}, X^{(i)} \rangle \right)^2  = \frac{n\| A^{\ell_1} - A^{\ell_2}\|_n^2}{2\sigma^2}.
\end{equation*}
Further if (\ref{AssCov}) holds, then
\begin{equation*}
\mathbb{E}_X[I_X(T; Y)] \leq \frac{n}{2\sigma^2 {m \choose 2} } \sum_{\ell_1 \neq \ell_2} {\mathbb{E}_X[\| A^{\ell_1} - A^{\ell_2}\|_n^2] } \leq c_u^2\frac{n}{2\sigma^2 {m \choose 2} } \sum_{\ell_1 \neq \ell_2} {\| A^{\ell_1} - A^{\ell_2}\|_{\rm F}^2 }.
\end{equation*}

Based on our construction, there exists a set $\{A^1, A^2,\ldots,A^{m} \}$ where each $A^\ell \in \mathcal{T}$ such that $\log m \geq C n \delta^2$ and 
$$c_u^{-1} \delta \leq \|A^{\ell_1} - A^{\ell_2} \|_{\rm F} \leq 8 c_u^{-1} \delta$$
for all $\ell_1 \neq \ell_2 \in \{1,2,\ldots,m\}$.
If (\ref{AssCov}) holds, then
$$\mathbb{E}_X\left(\| A^{\ell_1} - A^{\ell_2}\|_n^2\right) \leq c_u^2\| A^{\ell_1} - A^{\ell_2}\|_{\rm F}^2$$
and we can conclude that
\begin{equation*}
\mathbb{E}_X[I_X(T; Y)] \leq 32 c_u^2 \sigma^{-2} n \delta^2,
\end{equation*}
and from the earlier bound due to Fano's inequality, for and $\delta > 0$ such that
$$\frac{32 c_u^2 \sigma^2 n \delta^2 + \log 2}{\log m} \leq \frac{1}{2},$$
we are guaranteed that
$$\mathbb{P}\left\{\widetilde{T} \neq A^{\widetilde{m}}\right\} \geq \frac{1}{2}.$$
The proof is now completed because $\log m \geq 128 n \delta^2$  and $32 n \delta^2 \geq \log 2$.
\fpro
\vskip 25pt

\spro[Proof of Theorem \ref{ThmUpperMultiReg}]
The proof for the upper bound follows directly from Lemma~\ref{LemSparsityslice} with $d_1 = d_2 = m$ and $d_3 = p$ and noting that the overall covariance $\Sigma \in \mathbb{R}^{(nD_M) \times (nD_M)}$ is block-structured with blocks $\widetilde{\Sigma}$ since each of the samples is independent. Hence
$$c_\ell^2 \leq \lambda_{\min}(\Sigma) \leq \lambda_{\max}(\Sigma) \leq c_u^2.$$

To prove the lower bound, we use Theorem~\ref{ThmLower} and construct a suitable packing set for $\mathcal{T}_1$. The way we construct this packing is to construct two separate packing sets and select the set with the higher packing number using a similar argument to that used in \cite{RasWaiYu12} which also uses two separate packing sets. The first packing set we consider involves selecting the $s$-dimensional slice $A_{..S}$ where $A \subset[j_3]$ and $S = \{1,2,...,s\}$. Consider vectorizing each slice so $v = \mbox{vec}(A_{..S}) \in \mathbb{R}^{sm^2}$. Hence in order to apply Theorem~\ref{ThmLower},  we define the set $\mathcal{T}$ to be slices which is isomorphic to the vector space $\mathbb{R}^{sm^2}$. Using Lemma~\ref{LemFullHuperCube} in Appendix~\ref{SecHypercube}, there exists a packing set $\{v^1, v^2,...,v^N\} \in \mathbb{R}^{sm^2}$ such that $\log N \geq c s m^2$ and for all $v^{\ell_1}, v^{\ell_2}$ where $\ell_1 \neq \ell_2 $,
$$\frac{\delta^2}{4} \leq \|v^{\ell_1} - v^{\ell_2}\|_{\rm F}^2 \leq \delta^2$$
for any $\delta > 0$. If we choose $\delta = c\sqrt{s}m/\sqrt{n}$, then Theorem~\ref{ThmLower} implies the lower bound
\begin{equation*}
\min_{\widetilde{T}} \max_{T \in \mathcal{T}_1} \| \widetilde{T} - T \|_{\rm F}^2  \geq cc_u^{-2}\frac{sm^2}{n},
\end{equation*}
with probability greater than $1/2$.

The second packing set we construct is for the slice $A_{11\cdot} \in \mathbb{R}^p$. Since in the third direction only $s$ of the $p$ co-ordinates are non-zero, the packing number for any slice is analogous to the packing number for $s$-sparse vectors with ambient dimension $p$. Letting $v = A_{11\cdot}$, we need to construct a packing set for
$$\{v \in \mathbb{R}^p\;|\;\|v\|_{\ell_0} \leq s \}.$$
Using Lemma~\ref{LemSparseHuperCube} in Appendix~\ref{SecHypercube}, there exists a discrete set $\{v^1, v^2,...,v^N\}$ such that $\log N \geq cs \log(p/s)$ for some $c >0$ and $$\frac{\delta^2}{8} \leq \|v^k - v^{\ell}\|_{\ell_2}^2 \leq \delta^2$$
for $k \neq \ell$ for any $\delta > 0$. Setting $\delta^2 = s n^{-1}\log (p/s)$,
\begin{equation*}
\min_{\widetilde{T}} \max_{T \in \mathcal{T}_1} \| \widetilde{T} - T \|_{\rm F}^2  \geq cc_u^{-2} \frac{s \log(p/s)}{n},
\end{equation*}
with probability greater than $1/2$.

Taking a maximum over lower bounds involving both packing sets completes the proof of the lower bound in Theorem~\ref{ThmUpperMultiReg}.
\fpro
\vskip 25pt

\spro[Proof of Theorem \ref{ThmUpperMultiReg1}]
The upper bound follows directly from Lemma~\ref{LemSparsityslice2} with $d_1 = d_2 = m$ and $d_3 = p$ and noting that the overall covariance $\Sigma \in \mathbb{R}^{(nD_M) \times (nD_M)}$ is block-structured with blocks $\widetilde{\Sigma}$ since each of the samples is independent.

To prove the lower bound, we use Theorem~\ref{ThmLower} and construct a suitable packing set for $\mathcal{T}_2$. Once again we construct two separate packings and choose the set that leads to the larger minimax lower bound. For our first packing set, we construct a packing along one slice. Let us assume $A = (A_{\cdot\cdot1},...,A_{\cdot\cdot p})$, where $\mbox{rank}(A_{\cdot\cdot1}) \leq r$ and
$$A_{\cdot\cdot 2} =\cdots= A_{\cdot\cdot p} = 0.$$
If we let $A_{\cdot\cdot 1} = M$ where $M \in \mathbb{R}^{m \times m}$ then $A = (M, 0,..,0) \in \mathbb{R}^{m \times m \times p}$. Using Lemma~\ref{LemLowRankHypercube} in Appendix~\ref{SecHypercube}, there exists a set $\{A^1, A^2,...,A^{N} \}$ such that $\log N \geq c r m$ and
$$\frac{\delta^2}{4} \leq \|A^{\ell_1} - A^{\ell_2}\|_{\rm F}^2 \leq \delta^2$$ 
for all $\ell_1 \neq \ell_2$ and any $\delta > 0$. Here we set $\delta = \sqrt{rm/n}$. Therefore using Theorem~\ref{ThmLower}
\begin{equation*}
\min_{\widetilde{T}} \max_{T \in \mathcal{T}_2} \| \widetilde{T} - T \|_{\rm F}^2  \geq c\sigma^2 c_u^{-2} \frac{rm}{n},
\end{equation*}
with probability greater than $1/2$. 

The second packing set for $\mathcal{T}_2$ involves a packing in the space of singular values since
$$\sum_{j=1}^p {\mbox{rank}(A_{\cdot\cdot j})} \leq r.$$ 
Let $\{\sigma_{jk}: k=1,\ldots,m\}$ be the singular values of the matrix $A_{\cdot\cdot j}$. Under our rank constraint, we have
$$\sum_{j=1}^p \sum_{k=1}^m {\mathbb{I}(\sigma_{jk} \neq 0)} \leq s.$$
Let $v \in \mathbb{R}^{mp}$ where
$$v = \mbox{vec}((\sigma_{jk})_{1 \leq j \leq p, 1 \leq k \leq m}).$$
Note that
$$\sum_{j=1}^p \sum_{k=1}^m {\mathbb{I}(\sigma_{jk} \neq 0)} \leq r$$
implies $\|v\|_{\ell_0} \leq r$. Using Lemma~\ref{LemSparseHuperCube}, there exists a set $\{v^1, v^2,...,v^N \}$, such that $\log N \geq c r \log(mp/r)$ and for all $\ell_1\neq \ell_2$,
$$\frac{\delta^2}{4} \leq \|v^{\ell_1}- v^{\ell_2}\|_{\ell_2}^2 \leq \delta^2$$
for any $\delta > 0$. If we set $\delta^2 = rn^{-1}\log(mp/r)$. Therefore using Theorem~\ref{ThmLower},
\begin{equation*}
\min_{\widetilde{T}} \max_{T \in \mathcal{T}_2} \| \widetilde{T} - T \|_{\rm F}^2  \geq c\sigma^2 c_u^{-2} \frac{r \log(mp/r)}{n},
\end{equation*}
with probability greater than $1/2$. Hence taking a maximum over both bounds,
\begin{equation*}
\min_{\widetilde{T}} \max_{T \in \mathcal{T}_2} \| \widetilde{T} - T \|_{\rm F}^2  \geq c \sigma^2 c_u^{-2} \frac{r \max\{m, \log(p/r), \log m \}}{n} = c\sigma^2 c_u^{-2} \frac{r \max\{m, \log(p/r)\}}{n},
\end{equation*}
with probability greater than $1/2$.
\fpro
\vskip 25pt

\spro[Proof of Theorem \ref{ThmUpperVAR}]
The upper bound with 
$$\lambda \ge 3 \sqrt{\frac{\max\{p,2\log m\}}{\mu_{\min} n}}$$
follows directly from Lemma \ref{LemSparsityfiber} with $d_1=p$ and $d_2 = d_3 = m$ and (\ref{AssCov}) is satisfied with $c_u^2 = 1/\mu_{\min}$ and $c_{\ell}^2 =1/\mu_{\max}$ according to (\ref{EqnCondVAR}).

To prove the lower bound is similar to the proof for the lower bound in Theorem~\ref{ThmUpperMultiReg}. Once again we use Theorem~\ref{ThmLower} and construct a two suitable packing sets for $\mathcal{T}_3$. The first packing set we consider involves selecting an arbitrary subspace
$$\widetilde{\mathcal{T}} := \{ A = (A_{j_1,j_2,j_3})_{j_1, j_2, j_3}\;|\; 1 \leq j_1 \leq \sqrt{s},\; 1 \leq j_2 \leq \sqrt{s},\; 1 \leq j_3 \leq p \}.$$
Now if we let $v = \mbox{vec}(A)$, then $v$ comes from an $sp$-dimensional vector space for any $A\in \widetilde{\calT}$. Using Lemma~\ref{LemFullHuperCube} in Appendix~\ref{SecHypercube}, there exists a packing set $\{v^1, v^2,...,v^N\} \in \mathbb{R}^{sp}$ such that $\log N \geq c s p$ and for all $v^{\ell_1}, v^{\ell_2}$ where $\ell_1 \neq \ell_2$,
$$\frac{\delta^2}{4} \leq \|v^{\ell_1} - v^{\ell_2}\|_{\rm F}^2 \leq \delta^2$$
for any $\delta > 0$. If we choose $\delta = \sqrt{sp/n}$, then Theorem~\ref{ThmLower} implies the lower bound
\begin{equation*}
\min_{\widetilde{T}} \max_{T \in \mathcal{T}_3} \| \widetilde{T} - T \|_{\rm F}^2  \geq c\sigma^2 c_u^{-2} \frac{sp}{n},
\end{equation*}
with probability greater than $1/2$. Further $c_u^2 = 1/\mu_{\min}$. 

The second packing set we construct is for the slice $A_{1,j_2,j_3}$ for any $1 \leq j_2, j_3 \leq m$. Since in the second and third direction only $s$ of the co-ordinates are non-zero, we consider the vector space
$$\{v \in \mathbb{R}^{m^2}\;|\; \|v\|_{\ell_0} \leq s\}.$$
Once again using the standard standard hypercube construction in Lemma~\ref{LemSparseHuperCube} in Appendix~\ref{SecHypercube}, there exists a discrete set $\{v^1, v^2,...,v^N\}$ such that $\log N \geq cs \log(m^2/s)$ for some $c >0$ and
$$\frac{\delta^2}{8} \leq \|v^{\ell_1} - v^{\ell_2}\|_{\ell_2}^2 \leq \delta^2$$
for $\ell_1 \neq \ell_2$ for any $\delta > 0$. Setting $\delta = sn^{-1} \log (m^2/s)$ yields
\begin{equation*}
\min_{\widetilde{T}} \max_{T \in \mathcal{T}_3} \| \widetilde{T} - T \|_{\rm F}^2  \geq c\sigma^2 c_u^{-2} \frac{s \log(m/\sqrt{s})}{n},
\end{equation*}
with probability greater than $1/2$. Taking a maximum over lower bounds involving both packing sets completes the proof of of our lower bound.
\fpro
\vskip 25pt

\spro[Proof of Theorem \ref{th:pairwise}]
The upper bound follows from a slight modification of the statement in Lemma~\ref{LemSparsityslice2}. In particular since $\calR(A)= \|A^{(12)}\|_\ast+\|A^{(13)}\|_\ast+\|A^{(23)}\|_\ast$, the dual norm is 
\begin{equation}
\label{eq:pairwise}
\calR^*(A)=\max_{1\le k_1<k_2\le 3}\|A^{(k_1 k_2)}\|_s.
\end{equation} 
Hence, following the same technique as used in Lemma~\ref{LemSparsityslice2}
\begin{equation}
\label{eq:pairwise}
\mathbb{E}[\calR^*(G)] \leq c\max_{1\le k_1<k_2\le 3}\sqrt{\frac{\max\{d_{k_1}, d_{k_2}\}}{n}} = c\sqrt{\frac{\max\{d_1, d_2, d_3\}}{n}}.
\end{equation}
It is also straightforward to see that $s(\calT_4) \leq r$.

To prove the lower bound, we construct three packing sets and select the one with the largest packing number. Recall that 
\begin{eqnarray*}
\calT_4=\{A\in \RR^{d_1\times d_2\times d_3}: A_{j_1j_2j_3}=A^{(12)}_{j_1j_2}+A^{(13)}_{j_1j_3}+A^{(23)}_{j_2j_3}, A^{(k_1,k_2)}\in \RR^{d_{k_1}\times d_{k_2}},\\
 A^{(k_1,k_2)}\one =\zero, \qquad {\rm and}\qquad  (A^{(k_1,k_2)})^\top\one =\zero\\
 \max_{k_1,k_2} {\rm rank}(A^{(k_1,k_2)})\le r\}.
\end{eqnarray*}
Therefore our three packings are for $A^{(12)} \in \mathbb{R}^{d_1 \times d_2}$, $A^{(13)} \in \mathbb{R}^{d_1 \times d_3}$, and $A^{(23)} \in \mathbb{R}^{d_2 \times d_3}$ assuming each has rank $r$. We focus on packing in $A^{(12)} \in \mathbb{R}^{d_1 \times d_2}$ since the approach is similar in the other two cases. Using Lemma~\ref{LemGaussTens} from Appendix~\ref{AppSupGauss} in combination with Theorem~\ref{ThmLower},
\begin{equation*}
\min_{\widetilde{T}} \max_{T \in \calT_4} \| \widetilde{T} - T \|_{\rm F}^2  \geq c\sigma^2 c_u^{-2} \frac{r \min\{d_1, d_2\} }{n},
\end{equation*}
with probability greater than $1/2$. Repeating this process for packings in $A^{(13)} \in \mathbb{R}^{d_1 \times d_3}$, and $A^{(23)} \in \mathbb{R}^{d_2 \times d_3}$ assuming each has rank $r$ and taking a maximum over all three bounds yields the overall minimax lower bound
\begin{equation*}
\min_{\widetilde{T}} \max_{T \in \calT_4} \| \widetilde{T} - T \|_{\rm F}^2  \geq cc_u^{-2} \frac{r \max\{d_1, d_2, d_3\}}{n},
\end{equation*}
with probability greater than $1/2$. 
\fpro

\bibliographystyle{plainnat}

\bibliography{Biblio_GroupLasso}

\appendix

\section{Results for Gaussian random variables}

\label{AppTailBounds}

In this section we provide some standard concentration bounds that we use throughout this paper. First, we provide standard $\chi^2$ tail bounds. 
due to Laurent and Massart~\cite{LauMas98}:
\blems 
\label{LemChi}
Let $Z$ be a centralized $\chi^2$ random variable with $m$ degrees of freedom. Then for all $x \geq 0$,
\begin{eqnarray*}
\mathbb{P}[Z - m \geq 2 \sqrt{mx} + 2x] & \leq exp(-x),\;\;and\\
\mathbb{P}[Z - m \leq -2 \sqrt{mx}] & \leq exp(-x).
\end{eqnarray*}
\elems

\subsection{Gaussian comparison inequalities}
\label{AppGaussComp}

The first result is a classical result from \cite{Anderson55}.

\blems[Anderson's comparison inequality]
\label{LemAnderson}
Let $X$ and $Y$ be zero-mean Gaussian random vectors with covariance $\Sigma_X$ and $\Sigma_Y$ respectively. If $\Sigma_X - \Sigma_Y$ is positive semi-definite then for any convex symmetric set $C$,
\begin{equation*}
\mathbb{P}(X \in C) \leq \mathbb{P}(Y \in C).
\end{equation*}
\elems

The following Lemma is Slepian's inequality~\cite{Slepian62} which allows to upper bound the supremum of one Gaussian process by the supremum of another Gaussian process.

\blems[Slepian's Lemma] 
\label{LemSlepian}
Let $\{G_s, s \in S\}$ and $\{H_s, s \in S\}$ be two centered Gaussian processes defined over the same index set $S$. Suppose that both processes are almost surely bounded. For each $s, t \in S$, if $\mathbb{E}(G_s - G_t)^2 \leq \mathbb{E}(H_s - H_t)^2$, then $\mathbb{E}[\sup_{s \in S} G_s] \leq \mathbb{E}[\sup_{s \in S} H_s]$. Further if $\mathbb{E}(G_s^2) = \mathbb{E}(H_s^2)$ for all $s \in S$, then
\begin{equation*}
\mathbb{P}\left\{\sup_{s \in S} G_s > x\right\} \leq \mathbb{P}\left\{\sup_{s \in S} H_s > x\right\},
\end{equation*}
for all $x > 0$.
\elems

Finally, we require a standard result on the concentration of Lipschitz functions over Gaussian random variables.

\btheos[Theorem 3.8 from ~\cite{Massart03}]
\label{ThmConcGaussLip}
Let $g \sim \mathcal{N}(0, I_{d \times d})$ be a $d$-dimensional Gaussian random variable. Then for any function $F : \mathbb{R}^d \rightarrow \mathbb{R}$ such that $|F(x) - F(y)| \leq L\|x-y\|_{\ell_2}$ for all $x, y \in \mathbb{R}^d$, we have
\begin{equation*} 
\mathbb{P}\big[|F(g) - \mathbb{E}[F(g)]| \geq t \big] \leq 2\exp \left(-\frac{t^2}{2L^2} \right),
\end{equation*}
for all $t > 0$.
\etheos

\section{Suprema for i.i.d. Gaussian tensors}

In this section we provide important results on suprema of i.i.d. Gaussian tensors over different sets.

\label{AppSupGauss}

\subsection{The group $\ell_2$-$\ell_\infty$ norm}

Let $G \in \mathbb{R}^{d_1 \times d_2}$ be an i.i.d. Gaussian matrix and define the set
$$V := \{(u, v) \in \mathbb{R}^{d_1} \times \mathbb{R}^{d_2}\;|\; \|u\|_{\ell_2} \leq 1\;, \|v\|_{\ell_1} \leq 1\}.$$
Using this notation, let us define the define the random quantity:
\begin{equation*}
M(G, V) := \sup_{(u,v) \in V} u^\top G v.
\end{equation*}
Then we have the following overall bound.
\blems
\label{LemSupGauss1}
\begin{equation*}
\mathbb{E}[M(G, V)] \leq 3(\sqrt{d_1} + \sqrt{\log d_2}).
\end{equation*}
\elems
\spro
Our proof user similar ideas to the proof of Theorem 1 in ~\cite{RasWaiYu10b}. We need to upper bound $\mathbb{E}[M(G, V)]$. We are taking the supremum of the Gaussian process
$$\sup_{\|u\|_{\ell_2} \leq 1,\; \|v\|_{\ell_2} \leq 1} u^\top G v.$$
We now construct a second Gaussian process $\widetilde{G}_{u, v}$ over the set $V$ and apply Slepian's inequality (see Lemma~\ref{LemSlepian} in Appendix~\ref{AppGaussComp}) to upper bound
$$\sup_{\|u\|_{\ell_2} \leq 1,\; \|v\|_{\ell_2} \leq 1} u^\top G v$$
by the supremum over our second Gaussian process. $\widetilde{G}_{u, v}$. In particular, let us define the process as:
\begin{equation*}
\widetilde{G}_{u, v} = g^\top u + h^\top v,
\end{equation*}
where the vectors $(g, h) \in \mathbb{R}^{d_1} \times \mathbb{R}^{d_2}$ are i.i.d. standard normals (also independent of each other). It is straightforward to show that both $u^\top G v$ and $g^\top u + h^\top v$ are zero-mean. Further it is straightforward to show that
\begin{equation*}
\mbox{Var}(\widetilde{G}_{u, v} - \widetilde{G}_{u', v'}) = \|u - u'\|_{\ell_2}^2 + \|v - v'\|_{\ell_2}^2.
\end{equation*}
Now we show that
$$\mbox{Var}(u^\top G v - u'^\top G v') \leq \|u - u'\|_{\ell_2}^2 + \|v - v'\|_{\ell_2}^2.$$
To this end, observe that
\begin{eqnarray*}
\mbox{Var}(u^\top G v - u'^\top G v') & = & \|u v^\top - u' v'^\top\|_{\rm F}^2 \\
& = & \|(u-u')v^\top + u'(v-v')^\top\|_{\rm F}^2\\
& = & \|v\|_{\ell_2}^2 \|u - u'\|_{\ell_2^2} + \|u'\|_{\ell_2}^2 \|v - v'\|_{\ell_2}^2 \\
&&+ 2 (u^\top u' - \|u'\|_{\ell_2}\|u\|_{\ell_2})(v^\top v' - \|v'\|_{\ell_2}\|v\|_{\ell_2}).   
\end{eqnarray*}
First note that $\|v\|_{\ell_2}^2 \leq \|v\|_{\ell_1}^2 \leq 1$ for all $v \in V$ and $\|u'\|_{\ell_2}^2 \leq 1$. By the Cauchy-Schwarz inequality, $v^\top v' - \|v'\|_{\ell_2}\|v\|_{\ell_2} \leq 0$ and $u^\top u' - \|u'\|_{\ell_2}\|u\|_{\ell_2} \leq 0$. Therefore 
\begin{eqnarray*}
\mbox{Var}(u^\top G v - u'^\top G v') & \leq & \|u - u'\|_{\ell_2^2} + \|v - v'\|_{\ell_2}^2.   
\end{eqnarray*}
Consequently using Lemma~\ref{LemSlepian}
\begin{equation*}
\mathbb{E}[M(G, V)] \leq \mathbb{E}[\sup_{\|u\|_{\ell_2} \leq 1}g^\top u + \sup_{\|v\|_{\ell_1} \leq 1} h^\top v].
\end{equation*}
Therefore:
\begin{eqnarray*}
\mathbb{E}[M(G, V)] & \leq & \mathbb{E}[\sup_{\|u\|_{\ell_2} \leq 1}g^\top u + \sup_{\|v\|_{\ell_1} \leq 1} h^\top v]\\
& = & \mathbb{E}[\sup_{\|u\|_{\ell_2} \leq 1}g^\top u] + \mathbb{E}[\sup_{\|v\|_{\ell_1} \leq 1} h^\top v]\\
& = & \mathbb{E}[\|g\|_{\ell_2}] + \mathbb{E}[\|h\|_{\ell_\infty}]. 
\end{eqnarray*}
By known results on Gaussian maxima \citep[see e.g.][]{LedTal91},
$$\mathbb{E}[\|h\|_{\ell_\infty}] \leq 3 \sqrt{\log d_2}$$
and
$$\mathbb{E}[\|g\|_{\ell_2}] \leq \sqrt{d_1} + o(\sqrt{d_1}) \leq \frac{3}{2} \sqrt{D_j}.$$
Therefore
\begin{eqnarray*}
\mathbb{E}[M(G, V)] & \leq & \frac{3}{2} \sqrt{d_1} + 3 \sqrt{\log d_2}. 
\end{eqnarray*}
\fpro

\subsection{Spectral norm of tensors}

Our proof is based on an extension of the proof techniques used for the proof of Proposition 1 in \cite{NegWain11}.
\blems 
\label{LemGaussTens}
Let $G \in \mathbb{R}^{d_1 \times d_2 \times\cdots\times d_N}$ be a random sample from an i.i.d. Gaussian tensor ensemble. Then we have
\begin{equation*}
\mathbb{E}[\|G\|_s] \leq 4 \log(4N) \sum_{k=1}^{N}{\sqrt{d_k}}.
\end{equation*}
\elems
\spro
Recall the definition of $\|G\|_s$:
\begin{equation*}
\|G\|_s = \sup_{(u_1, u_2,\ldots,u_N) \in S^{d_1 - 1}\times S^{d_2 - 1}\times\cdots\times S^{d_N-1}} \langle u_1 \otimes u_2 \otimes \cdots\otimes u_N, G\rangle.
\end{equation*}
Since each entry $\langle u_1 \otimes u_2 \otimes \cdots\otimes u_N, G\rangle$ is a zero-mean Gaussian random variable, $\|G\|_s$ is the supremum of a Gaussian process and therefore the concentration bound follows from Theorem 7.1 in~\cite{Ledoux01}.

We use a standard covering argument to upper bound $\mathbb{E}[\|G\|_s]$. Let $\{u_1^1,u_1^2,\ldots,u_1^{M_1} \}$ be a $1/2N$ covering number of the sphere $S^{d_1 - 1}$ in terms of vector $\ell_2$-norm. Similarly for all $2 \leq k \leq N$, let $\{u_k^1,u_k^2,\ldots,u_k^{M_k} \}$ be a $1/2N$ covering number of the sphere $S^{d_k - 1}$. Therefore
\begin{eqnarray*}
&&\langle u_1 \otimes u_2 \otimes \cdots \otimes u_{N-1} \otimes u_N, G\rangle\\
& \leq& \langle u_1 \otimes u_2 \otimes \cdots \otimes u_{N-1} \otimes u_N^j, G\rangle + \langle u_1 \otimes u_2 \otimes \cdots\otimes u_{N-1} \otimes (u_N - u_N^j), G\rangle.
\end{eqnarray*}
Taking a supremum over both sides, 
\begin{equation*}
\|G\|_s \leq \max_{j = 1,\ldots,M_N} \langle u_1 \otimes u_2 \otimes \cdots \otimes u_{N-1} \otimes u_N^j, G\rangle + \frac{1}{2N} \|G\|_s.
\end{equation*}
Repeating this argument over all $N$ directions,
\begin{equation*}
\|G\|_s \leq 2 \max_{j_1 = 1,2,\ldots,M_1,\ldots,j_N = 1,2,\ldots,M_N} \langle u_1^{j_1} \otimes u_2^{j_2} \otimes \cdots \otimes u_N^{j_N}, G\rangle.
\end{equation*}
By construction, each variable $ \langle u_1^{j_1} \otimes u_2^{j_2} \otimes \cdots \otimes u_N^{j_N}, G\rangle$ is a zero-mean Gaussian with variance at most $1$, so by standard bounds on Gaussian maxima,
\begin{equation*}
\mathbb{E}[\|G\|_s] \leq 4\sqrt{\log(M_1 \times M_2 \times\cdots\times M_N)} \leq 4 [\sqrt{\log M_1} +\cdots+ \sqrt{\log M_N}].
\end{equation*}
There exist a $1/2N$-coverings of $S^{d_k - 1}$ with $\log M_k \leq d_k \log(4N)$ which completes the proof. 
\fpro

\section{Hypercube packing sets}
\label{SecHypercube}

In this section, we provide important results for the lower bound results. One key concept is the so-called \emph{Hamming distance}.The Hamming distance is between two vectors $v \in \mathbb{R}^d$ and $v' \in \mathbb{R}^d$ is defined by:
\begin{equation*}
d_H(v, v') = \sum_{j=1}^{d}{\mathbb{I}(v_j \neq v_j')}.
\end{equation*}

\blems
\label{LemFullHuperCube}
Let $\mathcal{C} = [-1, +1]^{d}$ where $d \geq 6$. Then there exists a discrete subset $\{v^{1}, v^{2},...,v^{m}\} \subset \mathcal{C}$, such that $\log m \geq c d$ for some constant $c > 0$, and for all $\ell_1 \neq \ell_2$,
\begin{equation*}
\frac{\delta^2}{4} \leq \|v^{\ell_1} - v^{\ell_2}\|_{\ell_2}^2 \leq \delta^2,
\end{equation*}
for any $\delta > 0$.
\elems
\spro
Let
$$v^{\ell} \in \left\{-\frac{\delta}{\sqrt{d}}, \frac{\delta}{\sqrt{d}}\right\}^{d},$$
i.e. a member of the $d$-dimensional hypercube re-scaled by $\sqrt{3}\delta/(2 \sqrt{d})$. Recall the definition of Hamming distance provided above. In this case amounts to the places either $v_j$ or $v'_j$ is negative, but both or not negative. Then according to Lemma 4 in ~\cite{Yu}, there exists a subset  re-scaled of this hypercube $v^{1}, v^{2},...,v^{m}$, such that
$$d_H(v^{\ell_1}, v^{\ell_2}) \geq \frac{d}{3}$$
and $\log m \geq c d$. Clearly, 
\begin{equation*}
\|v^{\ell_1} - v^{\ell_2}\|_{\ell_2}^2 = \frac{3\delta^2}{4 d} d_H(v^{\ell_1}, v^{\ell_2}) \geq \frac{\delta^2}{4}.
\end{equation*}
Further, 
\begin{equation*}
\|v^{\ell_1} - v^{\ell_2}\|_{\ell_2}^2 \leq \frac{3\delta^2}{4 d} \times d \leq \frac{3\delta^2}{4} \leq \delta^2.
\end{equation*}
This completes the proof.
\fpro

Next we provide a hupercube packing set for the sparse subset of vectors. That is the set
$$V := \{v \in \mathbb{R}^d\;|\; \|v\|_{\ell_0} \leq s \}.$$
This follows from Lemma 4 in \cite{RasWaiYu11} which we state here for completeness.

\blems
\label{LemSparseHuperCube}
Let $\mathcal{C} = [-1, +1]^{d}$ where $d \geq 6$. Then there exists a discrete subset $\{v^{1}, v^{2},...,v^{m}\} \subset V\cap \mathcal{C}$, such that $\log m \geq cs\log(d/s)$ for some $c > 0$, and for all $\ell_1 \neq \ell_2$,
\begin{equation*}
\frac{\delta^2}{8} \leq \|v^{\ell_1} - v^{\ell_2}\|_{\ell_2}^2 \leq \delta^2,
\end{equation*}
for any $\delta > 0$.
\elems

Finally we present a packing set result from Lemma 6 in \cite{AgarNegWainwright} that packs into the set of rank-$r$ $d_1 \times d_2$ matrices.

\blems
\label{LemLowRankHypercube}
Let $\min\{d_1, d_2\} \geq 10$, and let $\delta > 0$. Then for each $1 \leq r \leq \min\{d_1, d_2\}$, there exists a set of $d_1 \times d_2$ matrices $\{A^1, A^2,...,A^m\}$ with rank-$r$ with cardinality $\log m \geq c r \min\{d_1, d_2\}$ for some constant $c > 0$ such that
\begin{equation*}
\frac{\delta^2}{4} \leq \|A^{\ell_1} - A^{\ell_2}\|_{\rm F}^2 \leq \delta^2,	
\end{equation*} 
for all $\ell_1\neq \ell_2$.
\elems

\end{document}